\documentclass[10pt, reqno]{amsart}
\usepackage{amssymb}
\usepackage{amsthm}
\usepackage{graphicx}
\usepackage{mathabx}

\DeclareMathOperator{\grad}{grad}
\newcommand{\Laplacian}{\Delta}
\DeclareMathOperator{\diver}{div}
\DeclareMathOperator{\arcsinh}{arcsinh}
\newcommand{\supnorm}[1]{\vvvert #1\vvvert}
\newcommand{\llangle}{\langle\!\langle}
\newcommand{\rrangle}{\rangle\!\rangle}
\newcommand{\termi}{\text{I}}
\newcommand{\termii}{\text{II}}
\newcommand{\fordiff}{\nabla_+}
\newcommand{\backdiff}{\nabla_-}
\newcommand{\fordiffone}{\nabla_{+,1}}
\newcommand{\backdiffone}{\nabla_{-,1}}
\newcommand{\etadiffs}{\eta_s}
\newcommand{\etadifft}{\eta_t}
\newcommand{\sigmadiffs}{\sigma_s}
\newcommand{\sigmadifft}{\sigma_t}
\newcommand{\Gdiffs}{G_s}
\newcommand{\thetadiffs}{\theta_s}
\newcommand{\thetadifft}{\theta_t}
\newcommand{\epsilondifft}{\varepsilon_t}
\newcommand{\etadiffss}{\eta_{ss}}
\newcommand{\sigmadiffss}{\sigma_{ss}}
\newcommand{\sigmadifftss}{\sigma_{tss}}
\newcommand{\etadiffst}{\eta_{st}}
\newcommand{\thetadiffss}{\theta_{ss}}
\newcommand{\Gdiffss}{G_{ss}}
\newcommand{\etadifftt}{\eta_{tt}}
\newcommand{\thetadifftt}{\theta_{tt}}
\newcommand{\etadifftx}{\eta_{tx}}
\newcommand{\sigmadiffst}{\sigma_{st}}
\newcommand{\phidifftt}{\varphi_{tt}}
\newcommand{\epsilondiffst}{\varepsilon_{st}}
\newcommand{\epsilondiffss}{\varepsilon_{ss}}
\newcommand{\epsilondifftt}{\varepsilon_{tt}}
\newcommand{\epsilondiffs}{\varepsilon_s}
\newcommand{\deltas}{\delta_{s}}
\newcommand{\deltass}{\delta_{ss}}
\newcommand{\etameandiffss}{\overline{\eta}_{ss}}
\newcommand{\etameandiffst}{\overline{\eta}_{st}}
\newcommand{\etameandiffs}{\overline{\eta}_s}

\newcommand{\sigmameandiffss}{\overline{\sigma}_{ss}}
\newcommand{\etameandiffsss}{\overline{\eta}_{sss}}
\newcommand{\etameandifft}{\overline{\eta}_t}

\newcommand{\alphadiffr}{\alpha_r}
\newcommand{\phidiffr}{\varphi_r}
\newcommand{\phidifft}{\varphi_t}
\newcommand{\phidiffrr}{\varphi_{rr}}
\newcommand{\alphadiffrr}{\alpha_{rr}}

\newtheorem{theorem}{Theorem}[section]
\newtheorem{proposition}[theorem]{Proposition}
\newtheorem{lemma}[theorem]{Lemma}
\newtheorem{corollary}[theorem]{Corollary}
\newtheorem*{prop}{Proposition}
\newtheorem*{theo}{Theorem}
\newtheorem*{lemm}{Lemma}

\theoremstyle{definition}
\newtheorem{defn}[theorem]{Definition}
\newtheorem{exmp}[theorem]{Example}

\newtheorem{rem}[theorem]{Remark}

\numberwithin{equation}{section}

\begin{document}

\title{The motion of whips and chains}
\author{Stephen C. Preston}
\address{Department of Mathematics, University of Colorado, Boulder,
CO 80309-0395} \email{Stephen.Preston@colorado.edu}
\date{\today}

\maketitle

\section*{Abstract}

We study the motion of an inextensible string (a whip) fixed at one point in
the absence of gravity, satisfying the equations
$$ \etadifftt = \partial_s(\sigma \etadiffs), \qquad \sigmadiffss -
\lvert \etadiffss\rvert^2 \sigma = -\lvert \etadiffst\rvert^2, \qquad
\lvert \etadiffs\rvert^2 \equiv 1$$ with boundary conditions
$\eta(t,1)=0$ and $\sigma(t,0)=0$. We prove local existence and
uniqueness in the space defined by the weighted Sobolev energy
$$ \sum_{\ell=0}^m \int_0^1 s^{\ell} \lvert \partial_s^{\ell} \etadifft\rvert^2 \,ds + \int_0^1 s^{\ell+1} \lvert \partial_s^{\ell+1}\eta\rvert^2 \, ds,$$
when $m\ge 3$. In addition we show persistence of smooth solutions as long 
as the energy for $m=3$ remains bounded. 
We do this via the method of lines, 
approximating with a discrete system of coupled pendula (a chain)
for which the same estimates hold.

\tableofcontents

\section{Introduction and background}

\subsection{Introduction}

In this paper, we explore the motion of a whip, modeled as an
inextensible string. We prove that the partial differential equation
describing this motion is locally well-posed in certain weighted
Sobolev spaces. In addition, we are interested in the motion of a
chain, modeled as a coupled system of $n$ pendula, in the limit as
$n$ approaches infinity. We show that the motion of the chain
converges to that of the whip.

Although the equations of motion are well-known and have been studied by many authors, there are
few results known about the general existence and uniqueness
problem. Reeken~\cite{reeken2}~\cite{reeken3} proved local existence
and uniqueness for the \emph{infinite} string in $\mathbb{R}^3$ with
gravity and initial data sufficiently close (in $H^{26}$) to the vertical
solution, but aside from this, we know of no other existence result.
In the current paper we prove a local well-posedness
theorem for arbitrary initial data for the finite string.

One reason this problem is somewhat complicated is that the equation
of motion is hyperbolic, nonlinear, nonlocal, degenerate on
a spatial boundary, and possibly even elliptic under certain
conditions.

If $\eta \colon \mathbb{R}\times [0,1] \to \mathbb{R}^d$ describes
the position $\eta(t,s)$ of the whip, then one can derive that the equation of
motion in the absence of gravity and under the inextensibility
constraint $\langle \etadiffs, \etadiffs\rangle \equiv 1$ is
\begin{equation}\label{basicPDE}
\etadifftt(t,s) = \partial_s\big(\sigma(t,s) \etadiffs(t,s)\big).
\end{equation}
Incorporating gravity introduces some complications; to keep things
as simple as possible, we will neglect it.

Equation \eqref{basicPDE} is a standard wave equation; however, the tension $\sigma$ is
determined nonlocally, as a consequence of the inextensibility
constraint, by the ordinary differential equation
\begin{equation}\label{tensionODE}
\sigmadiffss(t,s) - \lvert \etadiffss(t,s)\rvert^2 \sigma(t,s) =
-\lvert \etadiffst(t,s)\rvert^2.
\end{equation}
With one end fixed and one end free, the boundary conditions are
$\eta(t,1)\equiv 0$ and $\sigma(t,0)\equiv 0$, along with
the compatibility condition $\partial_s\sigma(t,1)\equiv 0$.

We use the energy
$$ E_m = \sum_{\ell=0}^m \int_0^1 s^{\ell} \lvert \partial_s^{\ell} \etadifft(t,s)\rvert^2 \, ds
+ \int_0^1 s^{\ell+1} \lvert \partial_s^{\ell+1}\eta(t,s)\rvert^2 \, ds,$$
and show that for small time we have local existence and uniqueness
in the space for which the energy $E_3$ is bounded. Precisely, for any nonnegative integer $m$, define $N_m[0,1]$ to be the space 
of functions $f\colon [0,1]\to\mathbb{R}^d$ such that 
\begin{equation}\label{Nnorm}
\lVert f\rVert^2_{N_m} = \sum_{\ell=0}^m \int_0^1 s^{\ell} \Big\lvert \frac{d^{\ell}f}{ds^{\ell}}\Big\rvert^2 \, ds
\end{equation}
is bounded; then $E_m = \lVert \etadifft\rVert^2_{N_m} + \lVert \eta\rVert^2_{N_{m+1}}$.

We prove the following result:
\begin{theorem}\label{theorem1}
Suppose $\gamma\colon [0,1]\to \mathbb{R}^d$ and $w\colon [0,1]\to \mathbb{R}^d$
are restrictions of functions on $[0,2]$ satisfying the oddness condition $\gamma(2-s)=-\gamma(s)$ 
and $w(2-s)=-w(s)$, and that they have bounded weighted Sobolev norms,
$\lVert \gamma\rVert_{N_4} < \infty$
and
$ \lVert w\rVert_{N_3} < \infty$.
Suppose that in addition we have
$$ \lvert \gamma'(s)\rvert^2 \equiv 1 \qquad \text{and}\quad \langle \gamma'(s), w'(s)\rangle \equiv 0 \text{ for all $s\in[0,1]$.}$$

Then there is a $T>0$ such that there is a unique solution $\eta$ of
the equation \eqref{basicPDE} in $L^{\infty}([0,T], N_4[0,1])\cap
W^{1,\infty}([0,T], N_3[0,1])$ satisfying $\eta(0,s)=\gamma(s)$ and $\etadifft(0,s) = w(s)$. 
\end{theorem}
We prove this by showing that the corresponding discrete
energy $e_3$ for the chain with $n$ links is uniformly bounded for
small time, independently of $n$. The solution is then a weak-* limit of 
the chain solutions in $N_4$, which converges strongly in $N_3$ and
hence in $C^2$. One could prove this more directly using a Galerkin method,
but the present technique allows us to simultaneously discuss convergence
of the discrete approximation.

All the higher energies $E_m(t)$ can be bounded in terms of $E_3(t)$, so that $C^{\infty}$ 
initial conditions yield $C^{\infty}$ solutions for short time.
As a consequence, we derive a simple global existence criterion: if the initial conditions are 
$C^{\infty}$ functions, then a $C^{\infty}$ solution exists on $[0,T]$ if and only if
$E_3(t)$ is uniformly bounded on $[0,T]$. 
Of course, one expects blowup of the whip equation, at least for
some initial data, since the whole purpose of a whip is to construct
the initial condition so that the velocity of the free end
approaches infinity after a short time. See McMillen and
Goriely~\cite{mg} for a discussion of such issues; although our
model neglects some of the phenomena they consider, one expects that
the situations are similar in many ways. For the heuristics of
blowup in our situation, see Thess et al.~\cite{tzn}. The simplest
blowup mechanism appears to be the closing off of a loop along the whip; as a loop
shrinks, there appears a kink in the whip, representing blowup of
both the curvature and the angular velocity.

The paper is organized as follows. In Section \ref{basicequationssection}, we discuss
the equations for the whip and derive the corresponding equations for a chain with $n$ links,
in terms of difference operators, emphasizing the role of odd and even extensions in order 
to get the fixed endpoint conditions satisfied automatically.
In Section \ref{greenfunctionsection} we discuss the solution of the tension equation \eqref{tensionODE}
in terms of  a Green function, showing that the 
tension is positive except at $s=0$ and deriving a similar result for the chain.
We also derive sharp upper and lower estimates for the Green function.
In Section \ref{weightedsobolevsection} we explain why we need weighted energies, and we derive the 
analogues of the Sobolev and Poincar\'e inequalities for weighted norms, which are used throughout the rest 
of the paper. 

In Section \ref{tensionboundssection}, we give estimates for the tension $\sigma$ in terms of $\eta$ and $\etadifft$. For 
the $C^1$ norms of $\sigma$ we use the bounds on the Green function; for higher derivatives we bound the 
weighted Sobolev norms of $\sigma$ in terms of those of $\eta$ and $\etadifft$. These bounds are used in Section 
\ref{energyestimatesection} to derive the main energy estimate, to bound the time derivative of one energy in terms of 
another energy. 
Section \ref{existenceuniquenesssection} contains the proof of Theorem \ref{theorem1}. 
 Uniqueness is proved using a low-order 
 estimate for the difference of two solutions. Finally in Section \ref{futureresearchsection} we discuss 
 related open problems. To make the paper a bit easier to read, we have moved all of the longer proofs into 
 an Appendix.

Victor Yudovich found several results on this problem, although he did not
publish anything on it to my knowledge. I learned of this problem from Alexander Shnirelman, and I would
like to thank him for many useful discussions about it.

\subsection{Background}

The study of the inextensible string is one of the oldest
applications of calculus, going back to Galileo, and yet it is still
being studied to this day. One is especially concerned about kinks in the solution
and what the appropriate jump conditions should be; authors such as
O'Reilly and Varadi~\cite{ov}, Serre~\cite{serre}, and
Reeken~\cite{reeken1} have discussed these issues in detail from
differing points of view.

The first problem to be studied was finding the shape of a hanging chain, 
first solved incorrectly by Galileo and then correctly by Leibniz and
Bernoulli, one of the first major applications of the calculus of
variations. The shape of small-magnitude vibrations of a chain hanging straight down
(in a linear approximation) goes back to the Bernoullis and
Euler~\cite{truesdell}, and is taught in textbooks today as an
example of Bessel functions; see Johnson~\cite{johnson} and
Schagerl-Berger~\cite{sb} for related problems.
Kolodner~\cite{kolodner}, Dickey~\cite{dickey1},
Luning-Perry~\cite{lp}, and Allen-Schmidt~\cite{as} studied the
problem of a uniformly rotating inextensible string, one of the few
other problems that can be solved more or less exactly.

Burchard and Thomas~\cite{bt} obtained a local well-posedness result for the related problem of
inextensible elastica, in which there is a potential energy term
reflecting a resistance to bending; however it is not clear whether
the solutions are preserved in the limit as the potential term goes
to zero, so this result does not help in the present situation.

Many authors have studied the problem of a vertically folded chain
falling from rest; this is a classical problem that appears in
several textbooks (\cite{antman}, \cite{dickey2}, \cite{hamel}, and
\cite{rosenberg}). In recent years the problem has been debated in
the physics literature, in particular the issue of whether energy is
conserved and whether the tip of the chain falls at an acceleration
equal to gravity or faster (\cite{cal} \cite{calmar}
\cite{capmaz} \cite{dSR} \cite{hhr} \cite{ih} \cite{ov} \cite{ssst}
\cite{tp} \cite{tpg} \cite{st}). See Wong-Yasui~\cite{wy} or
McMillen~\cite{mcmillen} for a good survey of the literature.

McMillen and Goriely (\cite{gm} and \cite{mg}) studied a tapered whip theoretically,
numerically, and experimentally, showing that the crack comes not from the tip
but rather from a loop that straightens itself out. They use a different
model, however, in which the tension depends locally on the configuration.
Thess et al.~\cite{tzn} studied the blowup problem for the
closed inextensible string, especially as a model of the blowup
problem for the Euler equations for a 3D ideal fluid. They found
evidence of blowup from loops closing off, showing numerically that
$\sup_s \lvert \etadiffst\rvert \simeq \frac{1}{T-t}$ and $\sup_s
\lvert \etadiffss\rvert \simeq \frac{1}{(T-t)^{3/2}}$, where $T$ is
the blowup time.

\section{The basic equations}\label{basicequationssection}

In this section, we present the equations for both whips and chains,
assuming no external forces. Our boundary conditions come from the 
assumption that one end of the whip or chain is held fixed at the origin, 
while the other end is free. We describe the whip as a function 
$\eta\colon [0,T]\times [0,1] \to \mathbb{R}^d$, and describe the chain as 
a sequence of functions $\eta_k\colon [0,T]\to \mathbb{R}^d$ for $1\le k\le n+1$.
Our formulas simplify if we assume the fixed point occurs at $s=1$, i.e., $\eta(t,1)=0$
for all $t$; for the chain, we assume $\eta_{n+1}(t)=0$ for all $t$.

\subsection{The whip equations}\label{whipsection}

We will just present the equations here with a sketch of the derivation; the
reader may refer to \cite{whipcurvature} for a detailed derivation
and discussion. Schagerl et al.~\cite{ssst} and Thess et
al.~\cite{tzn} also present derivations from minimum principles: the
basic idea is to find a critical point of the action $\int_0^T \int_0^1 \lvert
\etadifft(t,s)\rvert^2
\, ds \,dt$ subject to the
constraint $\lvert \etadiffs(t,s)
\rvert^2 \equiv 1$. 

A variation $\zeta$ must satisfy $\zeta(t,1)=0$ and
 $\langle \etadiffs(t,s), \partial_s\zeta(t,s)\rangle \equiv 0$,
and if $\eta$ is a critical point of the action, then 
$\int_0^T \int_0^1 \langle \etadifft, \partial_t\zeta\rangle \, ds \, dt = 0$ for all such $\zeta$.
Integrating by parts, we conclude that a critical $\eta$ must satisfy
the evolution equation 
\begin{equation}\label{generalstringevol}
\etadifftt(t,s) = \partial_s\big( \sigma(t,s) \etadiffs(t,s)\big)
\qquad \eta(t,1) = 0,
\end{equation}
for some function $\sigma$, where $\sigma(t,0)=0$.

Differentiating $\lvert \etadiffs\rvert^2\equiv 1$ twice with respect
to $t$, we find that $\sigma$ is determined by the following
boundary-value problem for an ordinary differential equation (for
each fixed $t$):
\begin{equation}\label{generalstringconst}
\begin{split}
\sigmadiffss(t,s) - \lvert \etadiffss(t,s)\rvert^2 \sigma(t,s) = -\lvert \etadiffst(t,s)\rvert^2,\qquad
\sigma(t,0) = 0, \quad
\sigmadiffs(t,1)=0.
\end{split}
\end{equation}

\begin{figure}[ht]
\begin{center}
\includegraphics[scale=0.3]{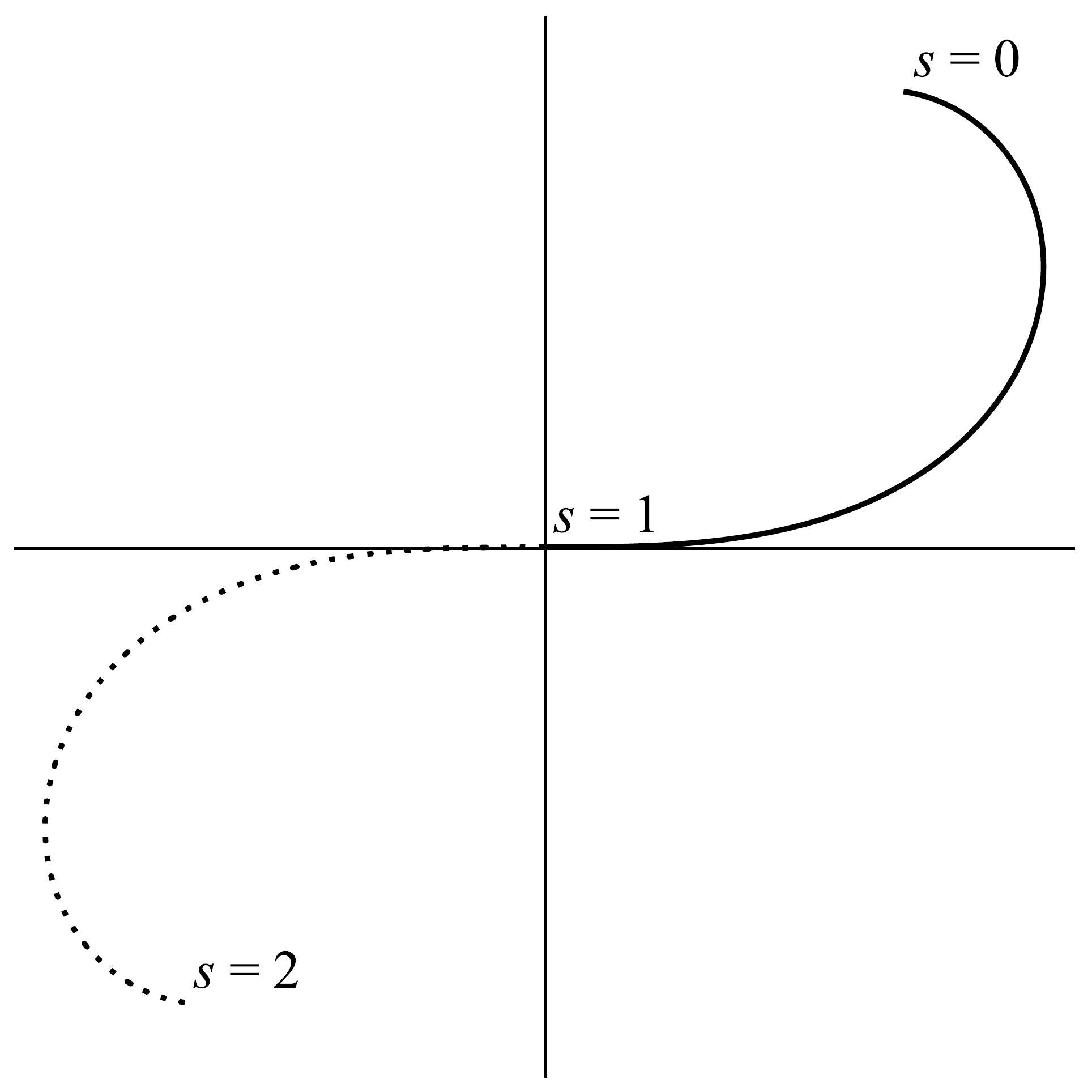}
\caption{The free end of the curve is at $s=0$, while the fixed end is at $s=1$. We imagine the curve
extending smoothly through the origin to $s=2$ through an odd reflection such that $\eta(s)=-\eta(1-s)$. Under such an extension, the tension extends to a smooth function satisfying $\sigma(s)=\sigma(1-s)$.}\label{oddcurve}
\end{center}
\end{figure}

The boundary conditions are compatible with the evolution equation as long as $\eta$ can be
extended to an odd function through $s=1$; in that case $\sigma$ can
be extended to an even function through $s=1$, which is where we get the extra boundary condition $\sigmadiffs(t,1)=0$. See Figure \ref{oddcurve}. Oddness and evenness give us the correct
boundary conditions for all higher derivatives of $\eta$ and $\sigma$ at $s=1$, which is crucial 
for the a priori estimates. Furthermore there is a discrete analogue of oddness and evenness 
for the chain which both simplifies the equations and helps greatly in defining the higher 
discrete energies. 

A sometimes simpler way of dealing with the constraint $\lvert \etadiffs\rvert \equiv 1$ is to consider $\etadiffs$ 
as a curve on the unit sphere in $S^{d-1}$. For simplicity we will assume $d=2$ when doing this, although 
the technique works in spherical coordinates in any dimension. We write 
\begin{equation}\label{spherical}
\etadiffs(t,s) = \big( \cos{\theta(t,s)}, \sin{\theta(t,s)}\big);
\end{equation}
a straightforward computation verifies that \eqref{generalstringevol} becomes
\begin{equation}\label{thetaevol}
\thetadifftt(t,s) = \sigma(t,s) \thetadiffss(t,s) + 2\sigmadiffs(t,s) \thetadiffs(t,s),
\end{equation}
while \eqref{generalstringconst} becomes
\begin{equation}\label{thetaconst}
\sigmadiffss(t,s) - \thetadiffs(t,s)^2 \sigma(t,s) = -\thetadifft(t,s)^2.
\end{equation}
The fact that $\eta$ is odd through $s=1$ forces us to have $\theta$ even through $s=1$, so the boundary condition on \eqref{thetaevol} is $\thetadiffs(t,1)=0$. We could work out all the estimates directly in terms of the system \eqref{thetaevol}--\eqref{thetaconst}, but the discrete versions of these equations are substantially more complicated than the discrete versions of \eqref{generalstringevol}--\eqref{generalstringconst}, even when $d=2$.

If $\sigma(t,s)$ is strictly positive for $0<s\le 1$, then equation
\eqref{generalstringevol} is a hyperbolic equation with a parabolic
degeneracy at $s=0$ (since we must have $\sigma(t,0)=0$). As such,
the only condition necessary to impose at $s=0$ is that $\eta(t,0)$
remain finite.

We point out that equation \eqref{generalstringevol} cannot be an ordinary 
differential equation on any infinite-dimensional Sobolev manifold: the right 
side is obviously an unbounded operator even in the simplest case. Hence we
cannot hope to prove existence and uniqueness using the techniques of
Picard iteration on an infinite-dimensional space, as in Ebin-Marsden~\cite{em}.  Instead
we will work directly with the partial differential equation using energy estimates.

\subsection{The chain equations}

We now derive the equations for the finite model, consisting of $(n+1)$
particles in $\mathbb{R}^d$, each of mass $\frac{1}{n}$, one of which is held fixed. The
particles are assumed to be joined by rigid links of length
$\frac{1}{n}$, whose mass is negligible. The position of the
$k^{\text{th}}$ particle is $\eta_{k}(t)$ for $1\le k\le n+1$; we assume the fixed
end is the $(n+1)^{\text{st}}$ particle, so that
$\eta_{n+1}(t)\equiv {0}$ for all time.\footnote{It might seem more natural to assume $\eta_0(t)=0$, but our choice makes the tensions $\sigma_k$ proportional to $k/n$ rather than $(n-k)/n$, simplifying some formulas.}
The configuration space is thus homeomorphic to $(S^{d-1})^n$, and is naturally embedded in $\mathbb{R}^{dn}$.

The kinetic energy in $\mathbb{R}^{dn}$ is
\begin{equation}\label{finitekinetic} K = \frac{1}{2n} \sum_{k=1}^n
\lvert\dot{\eta}_k \rvert^2.\end{equation} In addition the
constraints are given by
\begin{equation}\label{constrainerizer}
h_k(\eta_1, \ldots, \eta_n) = \frac{1}{2} \lvert \eta_{k+1} - \eta_k\rvert^2 = \frac{1}{2n^2},
\qquad 1\le k\le n.
\end{equation}
Stationary points of the constrained action satisfy the equations of motion $\ddot{\eta}_k = -\sum_{j=1}^n n^2 \sigma_j \, \partial_{\eta_k} h_j
$ for some Lagrange multipliers $\sigma_j$.
More explicitly, we have
\begin{equation}\label{pendulaevolution}
\ddot{\eta}_k = n^2\sigma_k (\eta_{k+1} - \eta_k) - n^2 \sigma_{k-1}
(\eta_k - \eta_{k-1})
\end{equation}
for $1< k\le n$. The scaling by $n^2$ is chosen so that $\sigma_k(t)$
converges to a function $\sigma(t,s)$ as
$n\to \infty$. The numbers $\sigma$ physically represent the
tensions in each link. We set $\sigma_0=0$ so the same equation is
valid when $k=1$.

The constraint equations determine the $\sigma$. Differentiating \eqref{constrainerizer} twice
with respect to time and using \eqref{pendulaevolution}, we get
\begin{multline}\label{sigmadiscreteraw}
-\lvert \dot{\eta}_{k+1}-\dot{\eta}_{k}\rvert^2=
n^2\sigma_{k+1}\langle \eta_{k+2} - \eta_{k+1},
\eta_{k+1}-\eta_{k}\rangle - 2 \sigma_{k}\\ + n^2 \sigma_{k-1}
\langle \eta_{k}-\eta_{k-1}, \eta_{k+1}- \eta_{k}\rangle
\end{multline}
for $1\le k< n$ (again using $\sigma_0=0$), while for $k=n$ we get (using $\eta_{n+1}={0}$)
\begin{equation}\label{sigmadiscreteendpoint}
-\lvert \dot{\eta}_n\rvert^2 = -\sigma_n - n^2 \sigma_{n-1}\langle
\eta_n, \eta_n-\eta_{n-1}\rangle.
\end{equation}

We note that if 
\begin{equation}\label{discretizationprime} 
\eta_k(t) = -\frac{1}{n} \sum_{j=k}^n \etadiffs\big(t,\tfrac{j}{n}\big) \quad \text{and} \quad \sigma_k(t) = \sigma(t, \tfrac{k}{n}),
\end{equation} where $\eta\colon\mathbb{R}\times[0,1]\to\mathbb{R}^d$ and $\sigma\colon \mathbb{R}\times [0,1]\to\mathbb{R}$ are $C^{\infty}$, then as $n\to\infty$, the formal limit of \eqref{pendulaevolution} is \eqref{generalstringevol} and the formal limit of \eqref{sigmadiscreteraw} is \eqref{generalstringconst}. (Note that this discretization of $\eta$ ensures that $\lvert (\fordiff\eta)_k\rvert = 1$ for all $k$, since $\lvert \etadiffs(t, \tfrac{j}{n})\rvert = 1$ for all $j$. We will refine this in Section \ref{interpolation}.) If $\eta(t,1)={0}$ and $\sigma(t,0)=0$, then this choice also gives $\eta_{n+1}(t)={0}$ and $\sigma_0(t)=0$, as desired. Hence the chain equations \eqref{pendulaevolution} and \eqref{sigmadiscreteraw} form a discretization of the whip equations \eqref{generalstringevol} and \eqref{generalstringconst} which conserves energy as well as preserving the geometry. 

The analysis of the chain equations becomes much simpler if we can avoid using separate 
equations for the boundary terms. An easy way to do this is to extend
$\eta_k$ and $\sigma_k$ beyond $k=n$ by demanding that $\eta$ be 
odd through $k=n+1$ and that $\sigma$ be even, which is exactly what 
we had to do for the whip in Section \ref{whipsection}. So for $k\ge n+1$ we set 
\begin{equation}\label{discreteodd}
\begin{split}
\eta_k &= -\eta_{2n+2-k},\\
\sigma_k &= \sigma_{2n+1-k}.
\end{split}
\end{equation}
Then it is easy to see that the evolution equation \eqref{pendulaevolution} still holds for the fixed 
point at $k=n+1$ and that \eqref{sigmadiscreteraw} for $k=n$ yields the tension boundary
condition \eqref{sigmadiscreteendpoint}. 

A further simplification comes from using difference operators. (See for example \cite{LL}.)
First recall that for a sequence $f$
defined on some subset of $\mathbb{Z}$, the (forward) shift operator $E$ is given by 
$(Ef)_k = f_{k+1}$. The backward shift is denoted by $E^{-1}$, so that $(E^{-1}f)_k = f_{k-1}$, and powers of $E$ signify composition.
We define the (forward) difference operator $\fordiff$ by
\begin{equation}\label{differenceoperator}
(\fordiff f)_k = n[f_{k+1}-f_k],
\end{equation} 
so that if $I$ denotes the identity operator, then $\fordiff = n(E-I)$.
It is also sometimes convenient to work with the backward difference operator
$\backdiff$, defined by $(\backdiff f)_k = n[f_k-f_{k-1}]$, so that $\backdiff = E^{-1}\fordiff = n(I-E^{-1})$.
In this notation\footnote{The more usual finite-difference notation is $\Delta$ for the forward difference and $\nabla$ for the backward difference; we use $\fordiff$ and $\backdiff$ instead to avoid confusion with the Laplacian on smooth functions, and since our rescaled version is not standard. We prefer the rescaling since if the sequence $f_k$ converges to a smooth function $f(s)$ as $n\to\infty$, then $(\fordiff f)_k$ converges to $f'(s)$.} equations \eqref{pendulaevolution} and \eqref{sigmadiscreteraw} become
\begin{align}
\ddot{\eta} &= \backdiff(\sigma \fordiff \eta), \label{fancydiscreteevolution} \\
\langle \fordiff\eta, \backdiff\fordiff (\sigma \fordiff \eta)\rangle &= -\lvert \fordiff \dot{\eta}\rvert^2\label{fancydiscretetension},
\end{align}
where the equations are valid when any subscript $1\le k\le n$ is placed on all the terms simultaneously. 
We can thus write all the discrete equations without specific reference to subscripts, which 
simplifies the notation.

The following formulas will be useful when working with difference operators and sums: both follow
from the simplest product formula $\fordiff (fg) = g \fordiff f + Ef \fordiff g$.
\begin{align}
\fordiff^{\ell} (fg) &= \sum_{j=0}^{\ell} {\ell \choose j} (E^j\fordiff^{\ell-j}f) (\fordiff^jg) \label{doubleproduct} \\
\frac{1}{n} \sum_{k=0}^{n-1} g_k \fordiff f_k &= -\frac{1}{n} \sum_{k=1}^{n} f_k \backdiff g_k +
 f_ng_n - f_0 g_0  \label{summationbyparts}
\end{align}

We can rewrite \eqref{fancydiscretetension} in a more useful form, solving for the second difference $\backdiff\fordiff\sigma$
in terms of everything else, using $\lvert \fordiff\eta\rvert^2\equiv 1$ to simplify the terms. We obtain
\begin{equation}\label{usefultension}
\backdiff\fordiff\sigma = \frac{E\sigma}{2} \lvert \fordiff^2\eta\rvert^2 + \frac{E^{-1}\sigma}{2} \lvert \backdiff\fordiff\eta\rvert^2 - \lvert \fordiff\dot{\eta}\rvert^2,
\end{equation}
and the resemblance to the continuous version \eqref{generalstringconst} is obvious.

\section{The Green function for the tension}\label{greenfunctionsection}

At each fixed time, equation \eqref{generalstringconst} is a linear nonhomogeneous ordinary differential equation for the tension $\sigma$. 
Hence there is a Green function $G(t,s,x)$ depending on $\lvert \etadiffss\rvert$, such that 
$$ \sigma(t,s) = \int_0^1 G(t,s,x) \lvert \etadifftx(t,x)\rvert^2 \, dx.$$ Similarly, equation \eqref{usefultension} can be thought of as a linear nonhomogeneous matrix equation for $\sigma$, for which the solution takes the analogous form 
$$ \sigma_k(t) = \frac{1}{n} \sum_{j=1}^n G_{kj}(t) \lvert \fordiff\dot{\eta}_j(t)\rvert^2$$
for some ``discrete Green function'' $G_{kj}$. Naturally we expect that if $\frac{k_n}{n} \to s$ and $\frac{j_n}{n}\to x$, then $G_{k_nj_n}(t) \to G(t,s,x)$ as $n\to\infty$; this can be proved as a consequence of our general convergence result for $\eta$. Our goal in this section is to establish properties of these Green functions. In particular we establish that the Green function is always nonnegative for a whip in a sufficiently smooth configuration, while the Green function is nonnegative for a chain as long as all the angles between links are obtuse. Furthermore we want to establish upper and lower bounds for the ratios $\frac{G(t,s,x)}{s}$ and $\frac{G_{kj}(t)}{s_k}$, where $s_k=\frac{k}{n}$, in order to be able to compare the norms weighted by powers of $\sigma(t,s)$ to the norms weighted by powers of $s$.

\subsection{Basic properties of the Green functions}

First we discuss the solution operator of the whip tension.
To keep the notation relatively simple, we will suppress
the time dependence.
\begin{proposition}\label{sigmasolprop}
For any fixed time $t$, the solution $\sigma(s)$ of \eqref{generalstringconst} is given by
\begin{equation}\label{sigmagreen}
\sigma(s) = \int_0^1 G(s,x) \lvert \etadifftx(x)\rvert^2 \, dx,
\end{equation}
where $G$ is the Green function given by
\begin{equation}\label{greencontinuous}
\Gdiffss(s,x)
- \lvert \etadiffss(s)\rvert^2
G(s,x) = -\delta(s-x), \quad G(0,x) = 0, \quad \Gdiffs(1,x)=0.
\end{equation}

The Green function is symmetric, i.e., $G(s,x)=G(x,s)$. It satisfies $G(s,x)>0$ whenever
$x>0$ and $0<s\le 1$. In addition if $0<x<1$, we have $\Gdiffs(s,x)>0$ for $0<s<x$ and $\Gdiffs(s,x)\le 0$ for
$x<s<1$.
\end{proposition}

\begin{proof}
The existence of the Green function and the symmetry property
$G(s,x)=G(x,s)$ is a well-known result of the general theory for
second-order equations with homogeneous boundary conditions. See for
example Courant-Hilbert~\cite{ch}.

To prove the other statements, we first show that $G(x,x)> 0$ for any $x\in (0,1)$.
For any fixed $x$, multiplying \eqref{greencontinuous} by $G(s,x)$, integrating from $s=0$ to $s=1$,
and using integration by parts with the homogeneous boundary conditions 
shows that 
$$ G(x,x) = \int_0^1 \Gdiffs(s,x)^2 \, ds + \int_0^1 \lvert \etadiffss(s)\rvert^2 G(s,x)^2 \, ds,$$
which forces $G(x,x)\ge 0$. Because of the jump condition 
$$\displaystyle \lim_{s\to x^-} \Gdiffs(s,x)-\lim_{s\to x^+} \Gdiffs(s,x)=1,$$ we cannot have $\Gdiffs(s,x)$ identically
zero if $0<x<1$. So $G(x,x)>0$ if $0<x<1$.
It is then easy to prove the other statements in the intervals $(0,x)$ and $(x,1)$ using the boundary conditions.
\end{proof}

Now let us do the same for the tension operator for the chain.
The equation \eqref{sigmadiscreteraw}, or the more elegant version \eqref{usefultension}, makes clear that the 
vector $(\sigma_1, \ldots, \sigma_n)$ of tensions comes from inverting a tridiagonal matrix. Since this is one of the easiest 
matrix types to invert, we get a relatively explicit formula for the solution, which will be useful in constructing estimates on the 
maximum and minimum tension.

\begin{proposition}\label{solverprop}
The solution of the constraint equations \eqref{fancydiscretetension}
is
\begin{equation}\label{finitetensionsolution} 
\sigma_k = \frac{1}{n} \sum_{j=1}^n G_{kj} \lvert \fordiff\dot{\eta}_j\rvert^2,
\end{equation}
where the discrete Green
function $G_{kj}$ is constructed by
\begin{equation}\label{discretegreen}
G_{kj} = \frac{1}{n} \sum_{i=1}^{\min{\{j,k\}}}
\frac{p_{ij}p_{ik}}{\beta_i}, \text{ where }
p_{ij} = \displaystyle\prod_{m=i}^{j-1} \frac{\alpha_m}{\beta_{m+1}}, \quad
\alpha_i = \langle \fordiff\eta_{i+1}, \fordiff\eta_i\rangle\end{equation}
and $\beta$ satisfies the recursion
\begin{equation}\label{betarecursion}
\beta_n = 1, \quad \beta_i = 2 - \frac{\alpha_i^2}{\beta_{i+1}} \text{ for $1\le i\le
n-1$}.\end{equation}
In \eqref{discretegreen} we use the convention that the empty product when $j=i$ is $1$.

The tensions $\sigma_k$ are positive for every nontrivial choice of
$\fordiff\dot{\eta}$ if and only if $\alpha_i> 0$ for every $i$.
\end{proposition}

\begin{proof}
The system \eqref{sigmadiscreteraw} and \eqref{sigmadiscreteendpoint} is of the form $A \sigma = w$, where 
$A$ is a symmetric nonnegative diagonally-dominant tridiagonal matrix and $w$ is the vector
of angular velocities $w_i = \lvert \fordiff\dot{\eta}_i\rvert^2$. There are several standard
algorithms for inverting such a matrix; the formula \eqref{finitetensionsolution} 
is given in the review paper of Meurant~\cite{meurant}.

Clearly $\alpha_i^2 \le 1$ for all $i$, so that inductively we have $1\le \beta_i\le 2$ for all $i$. 
Hence if all $\alpha_i$ are positive, then all $p_{ij}$ are positive and hence all terms $G_{kj}$ are positive
for $1\le j,k\le n$.
Thus if any $\fordiff\dot{\eta}_j$ is nonzero\footnote{Of course, the only way every $\fordiff\dot{\eta}_j$ is zero is if the chain is 
stationary, since $\dot{\eta}_{n+1}=0$ always.}, then \eqref{finitetensionsolution} says that all $\sigma_k$ are 
positive for $1\le k\le n$.
It is easy to see that if $\alpha_i\le 0$ for some $i$, then there is some choice of $\fordiff\dot{\eta}$ so that some $\sigma$ is nonpositive.
\end{proof}

\begin{figure}[ht]
\begin{minipage}[b]{0.25\linewidth}
\includegraphics[scale=0.15]{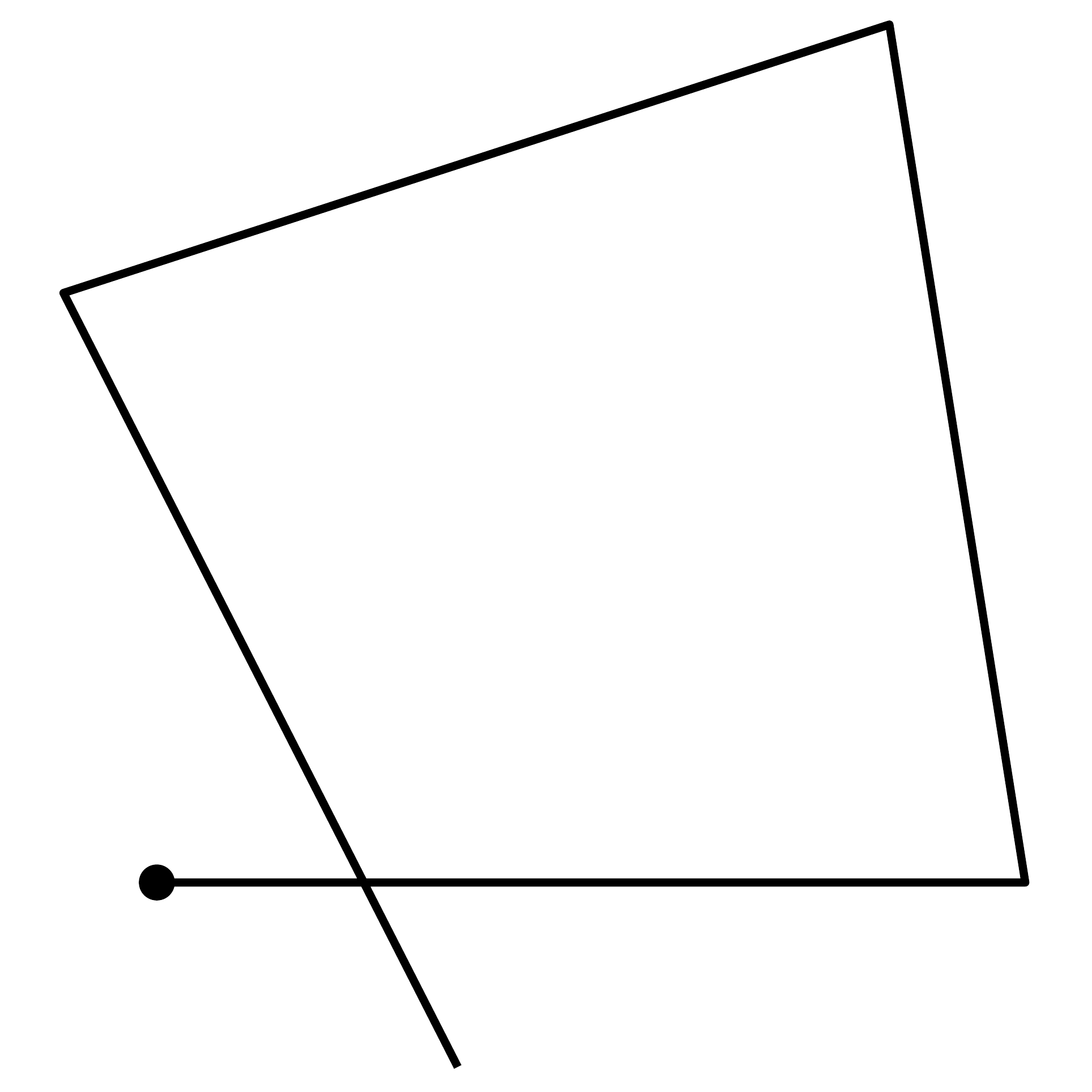}
\end{minipage}
\begin{minipage}[b]{0.25\linewidth}
\includegraphics[scale=0.15]{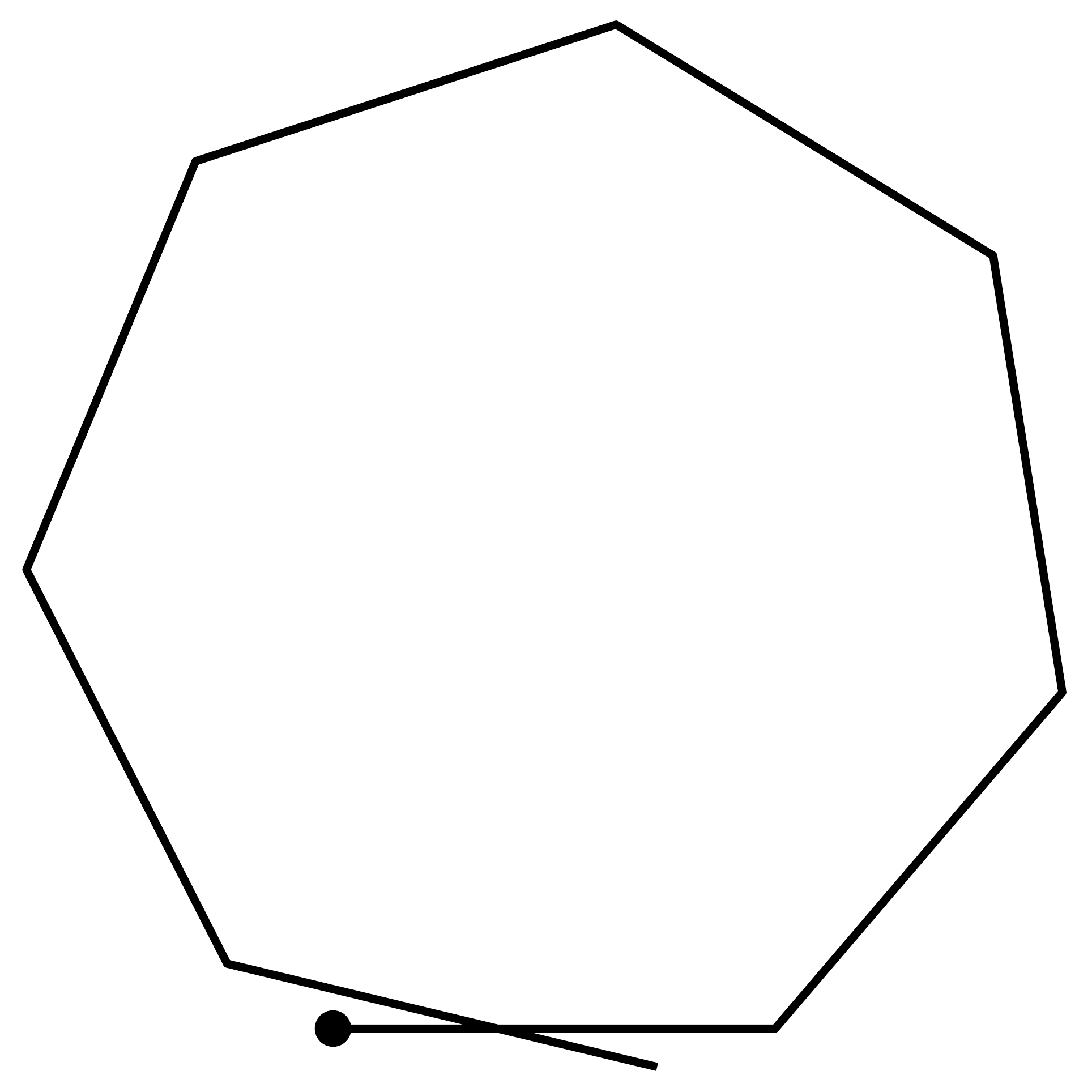}
\end{minipage}
\begin{minipage}[b]{0.25\linewidth}
\includegraphics[scale=0.15]{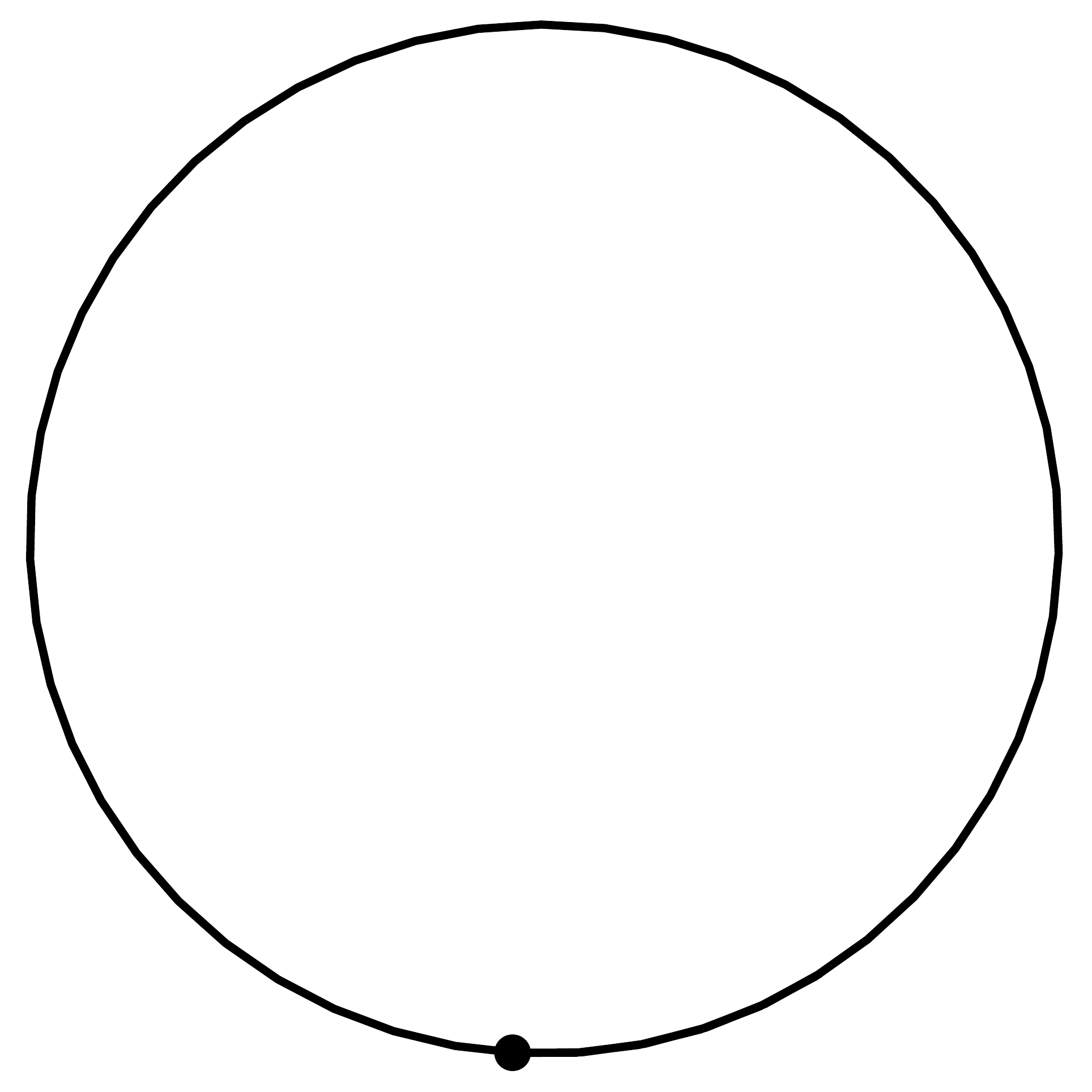}
\end{minipage}
\\
\begin{minipage}[b]{0.25\linewidth}
\includegraphics[scale=0.15]{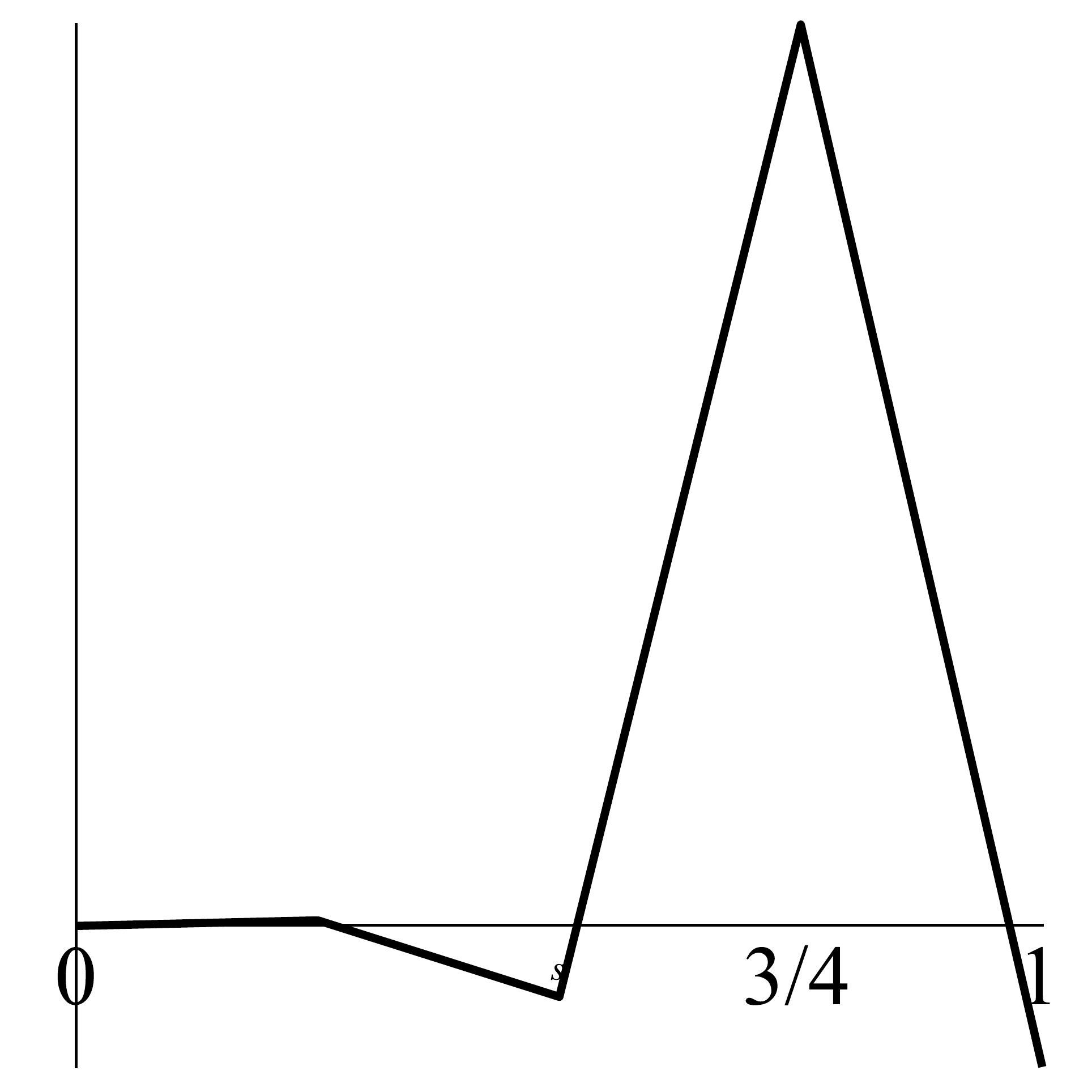}
\end{minipage}
\begin{minipage}[b]{0.25\linewidth}
\includegraphics[scale=0.15]{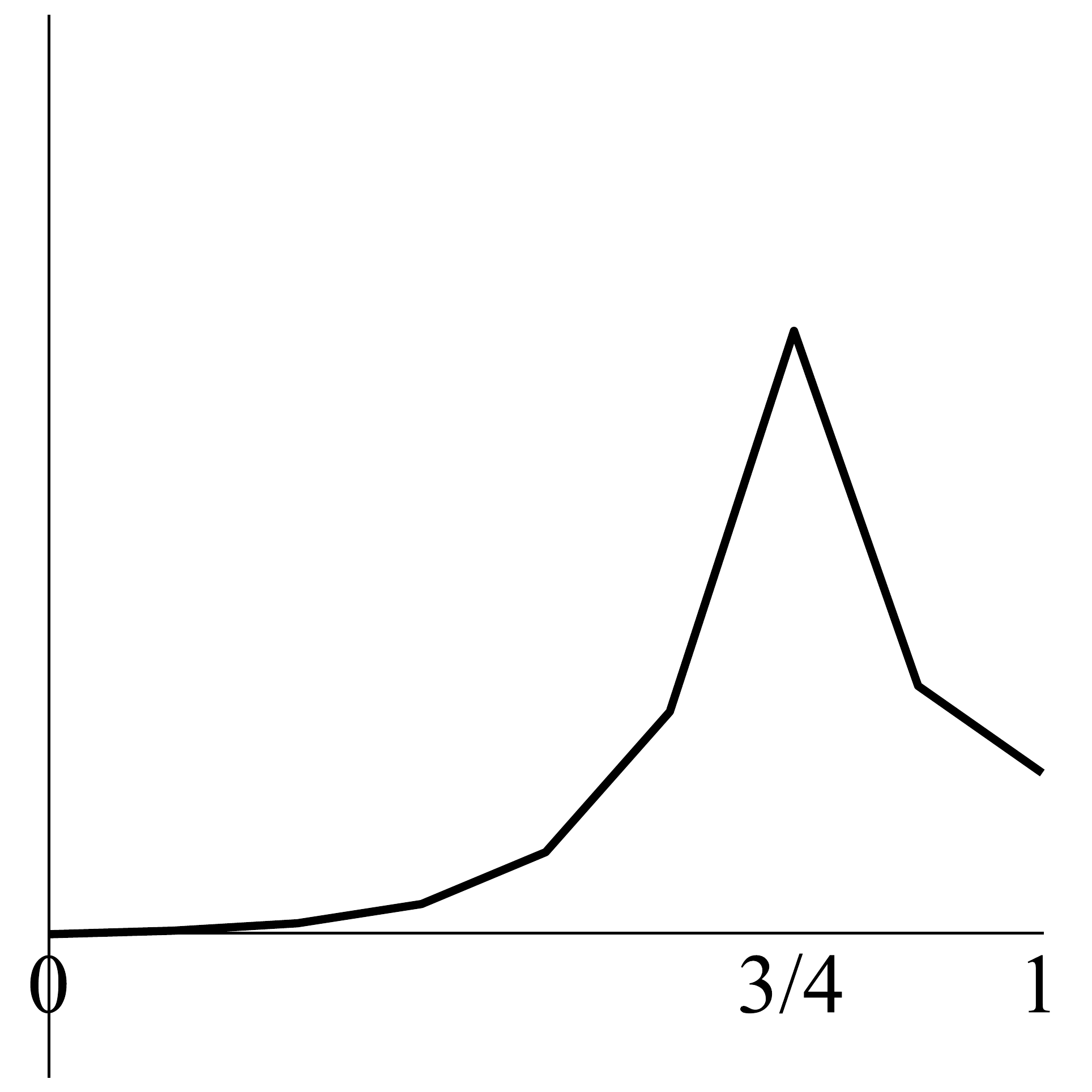}
\end{minipage}
\begin{minipage}[b]{0.25\linewidth}
\includegraphics[scale=0.15]{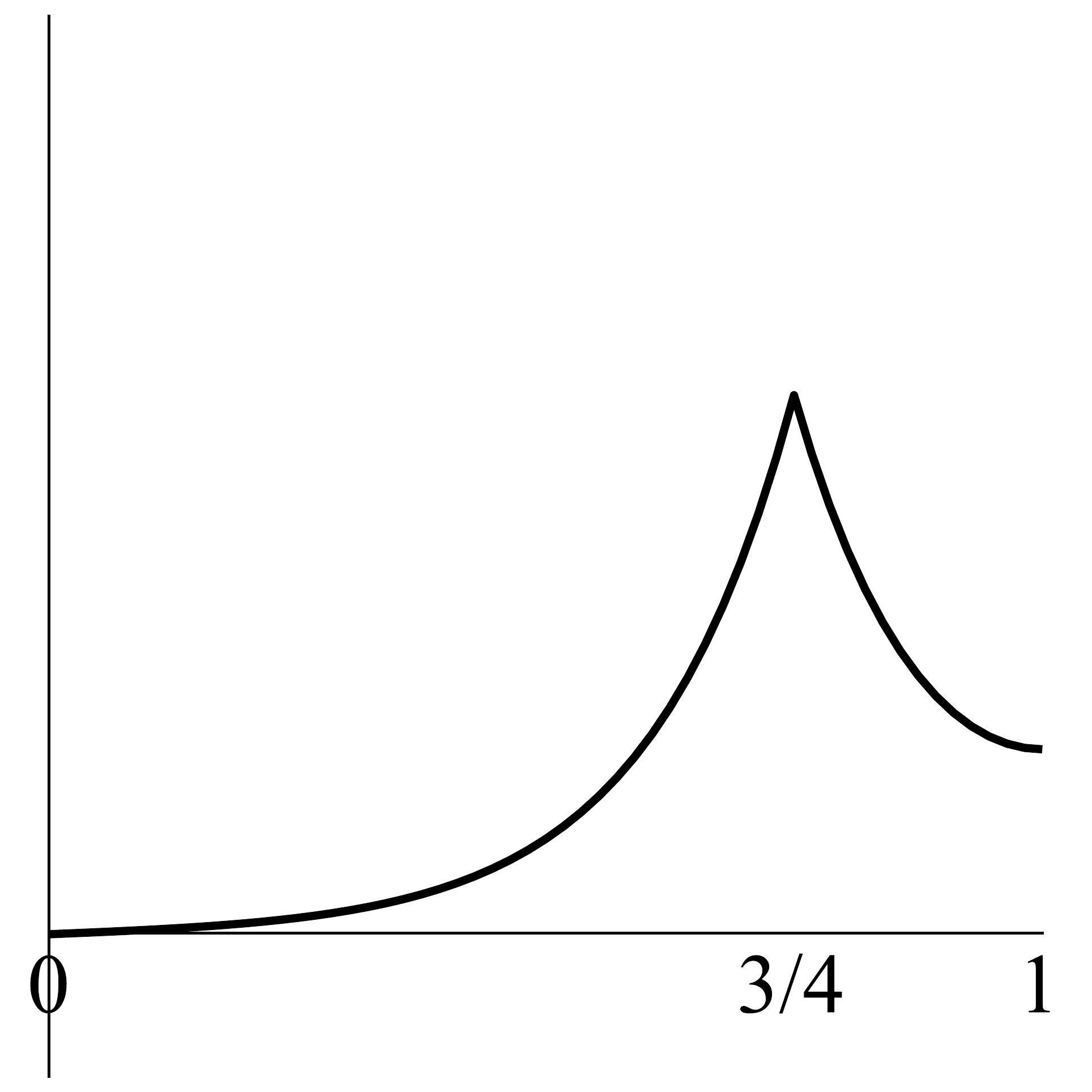}
\end{minipage}
\caption{In the top row we plot chains with constant angles 
between consecutive segments, with $n=4$, $n=8$, and $n=56$ links.
In the bottom row we plot the corresponding discrete Green function $G_{kj}$
as a function of $\frac{k}{n}$, evaluated at $j=\frac{3n}{4}$, for each configuration.
Notice that when $n=4$ the angles in the chain are acute, which is what allows the tension to
become negative in that case. Also notice that as $n\to\infty$, the discrete Green function
approaches the Green function for the differential equation.}
\end{figure}

\subsection{Upper and lower bounds for the Green functions}

Proposition \ref{sigmasolprop} implies that if $0<x\le 1$, then 
$G(s,x)/s$ is a positive function of $s$ on $[0,1]$, since $\lim_{s\to 0} G(s,x)/s = \Gdiffs(0,x)>0$. 
We now want to know exactly how large or small this positive function can be;
ultimately our interest will be in the quantities
$\sup_{s\in [0,1]} \sigma(s)/s$ and $\inf_{s\in [0,1]} \sigma(s)/s$, 
which are completely determined by the bounds on the Green
function. We are especially interested in the discrete analogues, 
$\max_{1\le k\le n} nG_{kj}/k$ and $\min_{1\le k\le n} nG_{kj}/k$. 
We end up with the same upper bound in both cases, which is relatively easy to prove,
while the lower bound is much more complicated and necessarily weaker in the discrete case.

First we establish the upper bound. 

\begin{proposition}\label{greenupperbounds}
If $\etadiffss$ is smooth, then the Green function $G(s,x)$ defined by Proposition
\ref{sigmasolprop} satisfies the following bounds.
\begin{equation}\label{upperboundsmooth}
\sup_{0\le s,x\le 1} \lvert \Gdiffs(s,x)\rvert \le 1, \text{ and }
\sup_{0\le s,x\le 1} \frac{G(s,x)}{s} \le 1.
\end{equation}

Furthermore, suppose $G_{kj}$, $\eta_k$, $\alpha_k$, and $\beta_k$ are as defined
in Proposition \ref{solverprop}, and that 
$\alpha_k = \langle \fordiff\eta_{k+1}, \fordiff\eta_k\rangle \ge 0$ for all
$k$, so that $G_{kj}\ge 0$ for all $j,k$. 
Then if $(\backdiffone G)_{kj}$ denotes the partial difference $(\backdiffone G)_{kj} \equiv
n(G_{kj}-G_{k-1,j})$, using the convention $G_{0j}=0$, then
\begin{equation}\label{upperbounddisc}
\lvert (\backdiffone G)_{kj}\rvert \le 1 \text{ and } \frac{nG_{kj}}{k} \le 1 \text{ for all $1\le j,k\le n$.}
\end{equation}
\end{proposition}

\begin{proof} 
The proof of \eqref{upperboundsmooth} is easy: by Proposition \ref{sigmasolprop}, 
the partial derivative $\Gdiffs(s,x)$ is positive for $s<x$ and nonpositive
for $s>x$, and jumps by $-1$ at $s=x$. Since $\Gdiffss\ge 0$ whenever $s\ne x$, we know $\Gdiffs$ is increasing on each interval. We therefore must have $0<\Gdiffs(s,x)\le 1$ for $s<x$
and $-1< \Gdiffs(s,x)\le 0$ for $x<s\le 1$; either way, $\lvert \Gdiffs(s,x)\rvert \le 1$.
Then using the fact that $G(0,x)=0$, we have 
$0 \le G(s,x) = \int_0^s \partial_rG(r,x)\,dr \le \int_0^s \, dr = s,$
which yields \eqref{upperboundsmooth}.

The proof of \eqref{upperbounddisc} is more complicated, but uses the same basic ideas.
First, from Proposition \ref{solverprop} we know that $G_{kj}\ge 0$ for all
$k$ and $j$ since every $\alpha_i\ge 0$. 

Assume first that $j\ne n$. Then rewriting \eqref{sigmadiscreteraw}--\eqref{sigmadiscreteendpoint},
we see that the discrete Green function satisfies the equation
$$ \alpha_{k+1} G_{k+1,j} - 2 G_{kj} + \alpha_{k-1} G_{k-1,j} = -\frac{1}{n} \delta_{kj}$$
for $1\le k<n-1$, while
$$ -G_{nj} + \alpha_{n-1} G_{n-1,j} = 0 $$
Since $\alpha_k = \langle \fordiff\eta_k, \fordiff\eta_{k+1}\rangle$ with $\lvert
\fordiff\eta_k\rvert = 1$, we have $\alpha_k\le 1$. Thus for $k\ne j$ we can easily see the second
partial difference satisfies
$(\fordiffone\backdiffone G)_{kj} 
\ge 0.$
Since the second partial differences are nonnegative except at the diagonal, the first partial differences are increasing except at the diagonal, i.e.,  
\begin{equation}\label{Gdiscreteconvex}
(\backdiffone G)_{k+1,j}-(\backdiffone G)_{kj}\ge 0
\text{ for all $k\ne j$.}
\end{equation}
When $k=j$ we can check that
\begin{equation}\label{Gdiscretejump}
(\backdiffone G)_{j+1,j}-(\backdiffone G)_{jj} 
\ge -1.
\end{equation}

Now look at the endpoint terms: at the left endpoint, we know $(\backdiffone G)_{1j} = nG_{1j}\ge 0$.
At the right endpoint, if $j\ne n$ then we have $-G_{nj}+G_{n-1,j} =
(1-\alpha_{n-1})G_{n-1,j}\ge 0$, so that $(\backdiffone G)_{nj} \le 0$.
Thus combining \eqref{Gdiscreteconvex} and \eqref{Gdiscretejump}, we conclude that if $j\ne n$ then
$$ 0\le (\backdiffone G)_{1j} \le \cdots \le (\backdiffone G)_{jj} \le 1 + (\backdiffone G)_{j+1,j} \le \cdots \le 1+(\backdiffone G)_{nj} \le 1.$$
Hence we must have $\lvert (\backdiffone G)_{kj}\rvert \le 1$ for all $k$, as
long as $j\ne n$.

If $j=n$, the situation is slightly different; in that case we 
get
$$ 0 \le (\backdiffone G)_{1n} \le (\backdiffone G)_{2n} \le \cdots \le (\backdiffone G)_{n-1,n}\le (\backdiffone G)_{nn}\le 1,$$
so that $\lvert (\backdiffone G)_{kj}\rvert \le 1$ even if $j=n$. This completes the proof of \eqref{upperbounddisc}.
\end{proof}

\begin{rem}
Unfortunately we cannot bound $\frac{G(s,x)}{sx}$ from above. If we denote by $G_0(s,x)$ the Green function when $\lvert \etadiffss\rvert \equiv 0$, then we easily compute that $G_0(s,x) = \min{\{s,x\}}$, so that $\frac{G_0(s,x)}{sx} = \min{\{\frac{1}{s}, \frac{1}{x}\}}$ is unbounded on $[0,1]\times [0,1]$. Note that by the Sturm comparison theorem, we have that $G(s,x)\le G_0(s,x)$ for any Green function satisfying \eqref{greencontinuous}. 
However it is easy to see that for any $0\le p\le 1$, we have $\frac{G(s,x)}{s^p x^{1-p}} \le \frac{G_0(s,x)}{s^px^{1-p}} \le 1$. This will be useful in the proof of Theorem \ref{uniqueness}.

It is easy to check that the discrete Green function satisfies the same inequality, $\lvert G_{kj}\rvert \le \frac{1}{n}\min{\{j,k\}}$, using formula \eqref{discretegreen} and the fact that $\lvert \alpha_i\rvert \le 1$ and $\beta_i \ge 1$ for all $i$. In fact this bound is valid even if not all $G_{kj}$ are positive. 
\end{rem}

Now we establish the lower bound. This is the only time in the paper where we get 
a weaker result for the chain than for the whip; the reason is that we need to
make strong assumptions in order to prevent sharp kinks in the chain, 
to ensure nonnegative tension. Smoothness of the whip, on the other hand, 
ensures that the tension in the whip is nonnegative automatically. 

\begin{proposition}\label{sigmalowerbound}
Suppose $G_{kj}$, $\eta_k$, $\alpha_k$, and $\beta_k$ are defined
as in Proposition \ref{solverprop}. Assume the $\eta_k$ are such
that, for some $\upsilon\in (0, \frac{2\sqrt{n}}{5}]$, we have 
\begin{equation}\label{threehalves}
\frac{k^{3/2}}{n^{3/2}} \lvert \fordiff^2\eta_k\rvert^2 \le \upsilon \quad \text{for all $1\le k\le n-1$.}
\end{equation}
Then for all $1\le j,k\le n$ we have
\begin{equation}\label{lowerbounddisc}
\frac{n^2G_{kj}}{jk} \ge e^{-2\upsilon}. 
\end{equation}

If $G$ solves \eqref{greencontinuous} and $\etadiffss$ is a smooth function, then we have
\begin{equation}\label{lowerboundcont}
\inf_{0\le s,x\le 1} \frac{G(s,x)}{sx} \ge \frac{e^{-\varrho}}{1+\varrho}
\text{ where } \varrho = \int_0^1 s\lvert \etadiffss\rvert^2 \, ds.
\end{equation}
\end{proposition}

\begin{proof}
The two estimates are proved in slightly different ways, but the main point for both estimates 
is to show that the minimum is attained at the off-diagonal corner, then estimate this value
either using the direct formula \eqref{discretegreen} (for the chain) or through a substitution (for the whip). 
The full proof is in Appendix \ref{A1}.
\end{proof}

The assumption \eqref{threehalves} for the discrete case is much stronger than the assumption 
$\int_0^1 s\lvert \etadiffss\rvert^2 \, ds<\infty$ for the continuous case, but such a pointwise bound is necessary to ensure
every $\alpha_k>0$ in order to get all tensions positive (by Proposition \ref{solverprop}), even when $k=1$. The exponent $\frac{3}{2}$ 
is important: the exponent $1$ would work to prove the estimate, but we cannot prove that such an estimate 
actually holds for all values of $t$; the exponent $2$ is not enough to get a lower bound for the tension.

\begin{rem}\label{cantweaken}
Note that we could easily get a stronger estimate than \eqref{lowerboundcont} 
if we simply assumed an upper bound on $\lvert
\etadiffss \rvert$, using the Sturm-Liouville comparison theorem. 
However, we prefer the weaker assumption that $\int_0^1 s
\lvert \etadiffss \rvert^2 \, ds < \infty$, since it allows for the
possibility of the curvature at the free end of the whip approaching
infinity (a possibility not precluded by the equations due to the
degeneracy there). Even if the weighted energy $E_3$ is finite---the 
condition under which we will prove local existence---we will not necessarily
have $\lvert \etadiffss\rvert$ bounded on $[0,1]$; an example is when 
$\theta(s)=s^q$ for some $q\in (\frac{1}{2}, 1)$, using the spherical representation
\eqref{spherical}. See Example \ref{ABCnecessary} for details.

\end{rem}

\section{Weighted Sobolev norms}\label{weightedsobolevsection}

\subsection{Motivation}

In order to demonstrate existence and uniqueness, we want to apply the usual technique of energy estimates in Sobolev spaces. By showing that we have sequences of solutions of the chain equations \eqref{fancydiscreteevolution}--\eqref{fancydiscretetension} for which the energy is uniformly bounded, we can extract a convergent subsequence to establish existence; the same sort of energy estimates can also be used to establish uniqueness. Several issues arise to complicate this strategy.

Ordinarily for a wave equation like \eqref{generalstringevol}, one would try to bound an energy like 
$$ \tilde{E}_0 = \int_0^1 \lvert \etadifft\rvert^2 + \sigma \lvert \etadiffs\rvert^2 \, ds$$
by computing its derivative and using Gronwall's lemma. It is easy to compute that 
$$
\frac{d\tilde{E}_0}{dt} = (\sigma \langle \etadiffs, \etadifft\rangle)\big|_{s=0}^{s=1} + \int_0^1 \sigmadifft \lvert \etadiffs\rvert^2 \, ds
\le \left(\sup_{0\le s\le 1} \frac{\sigmadifft(t,s)}{\sigma(t,s)}\right) \, \tilde{E}_0(t),
$$
using the boundary conditions $\sigma(t,0)=0$ and $\etadifft(t,1)=0$. 
Unfortunately we cannot bound $\sigmadifft$ or even $\sigma$ in terms only of $\tilde{E}_0$. Indeed, it is hard to even 
make sense of equation \eqref{generalstringconst} unless both $\etadiffst$ and $\etadiffss$ are in $L^2$, which means we have to consider higher energies. 
 
Here a complication arises. The usual approach would be to consider an energy like 
$$ \tilde{F}_1 = \tilde{E}_0 + \int_0^1 \lvert \etadiffst\rvert^2 + \sigma \lvert \etadiffss\rvert^2 \, ds.$$
Its derivative is, using \eqref{generalstringevol},
$$ \frac{d\tilde{F}_1}{dt} \le \left( \sup_{0\le s\le 1}\frac{\sigmadifft(t,s)}{\sigma(t,s)} \right) \tilde{F}_1(t) + 2 \int_0^1 \sigmadiffs \langle \etadiffst, \etadiffss\rangle \, ds.$$
Here the boundary term vanishes since $\sigma(t,0)=0$ and $\etadiffss(t,1)=0$ (recall we assume $\eta$ extends to an odd function through $s=1$). Furthermore since $\lvert \etadiffs\rvert^2 \equiv 1$, we have $\langle \etadiffs, \etadiffst\rangle \equiv 0$.
The problem is that if we want to get the right side in terms of $\tilde{F}_1$ alone, we need to use the Cauchy-Schwarz inequality to get
$$ 
\left\lvert \int_0^1 \sigmadiffs\langle \etadiffst, \etadiffss\rangle \, ds\right\rvert 
\le \left( \sup_{0\le s\le 1} \frac{\lvert \sigmadiffs(t,s)\rvert}{\sqrt{\sigma(t,s)}} \right) \tilde{F}_1(t),
$$
but the right side is \emph{not} bounded. We always have $\sigma(t,0)=0$, while we will generally not have $\sigmadiffs(t,0)=0$. 

Instead we want an energy for which the integration by parts cancels out this highest-order remainder. The only such quantity of the form $\int_0^1 A \lvert \etadiffst\rvert^2 + B \lvert \sigmadiffss\rvert^2 \, ds$ for which this works is 
$$ \tilde{E}_1 = \tilde{E}_0 + \int_0^1 \sigma \lvert \etadiffst\rvert^2 + \sigma^2 \lvert \etadiffss\rvert^2 \, ds.$$
With such a choice we get 
$$ \frac{d\tilde{E}_1}{dt} \le 2 \left( \sup_{0\le s\le 1} \frac{\sigmadifft(t,s)}{\sigma(t,s)}\right) \tilde{E}_1(t),$$
which we can manage once we understand how $\sigma$ and $\sigmadifft$ behave. The same phenomenon continues for the higher energies as well, which motivates us to define 
\begin{equation}\label{correctenergy}
\tilde{E}_m = \sum_{\ell=0}^m \int_0^1 \sigma^{\ell} \lvert \partial_s^{\ell} \etadifft\rvert^2 + \sigma^{\ell+1} \lvert \partial_s^{\ell+1}\eta\rvert^2 \, ds.
\end{equation}
With this definition, we have 
\begin{equation}\label{nineandahalf}
\begin{split}
\frac{d\tilde{E}_m}{dt} &= \sum_{\ell=0}^m \bigg[ \ell \int_0^1 \sigma^{\ell-1} \sigmadifft \lvert \partial_s^{\ell}\etadifft\rvert^2 \, ds + (\ell+1) \int_0^1 \sigma^{\ell} \sigmadifft \lvert \partial_s^{\ell+1}\eta\rvert^2 \, ds \\
&\qquad\qquad +  \sum_{i=0}^{\ell-1} 2\textstyle{\ell+1\choose i} \int_0^1 \sigma^{\ell} \partial_s^{\ell+1-i} \sigma \langle \partial_s^{\ell}\etadifft, \partial_s^{i+1}\eta\rangle \, ds\bigg],
\end{split}
\end{equation}
with the three remaining terms integrating to give $\sigma^{\ell+1} \langle \partial_s^{\ell} \etadifft, \partial_s^{\ell+1}\eta\rangle \vert_{s=0}^{s=1} = 0$ due to oddness of $\eta$ through $s=1$. 

Our primary goal will be to bound \eqref{nineandahalf} in terms of the energies $\tilde{E}_m$. More specifically, using the fact that $\sigma(t,s)$ degenerates like $s$ near $s=0$, we want to get bounds in terms of the simpler weighted energies 
$E_m$ defined by 
\begin{equation}\label{correctweightenergy}
E_m = \sum_{\ell=0}^m \int_0^1 \Big( s^{\ell} \lvert \partial_s^{\ell} \etadifft\rvert^2 + s^{\ell+1} \lvert \partial_s^{\ell+1}\eta\rvert^2 \Big) \, ds.
\end{equation}

To do this, we will need several estimates. So our first goal is establishing basic Sobolev-type and Wirtinger-type inequalities for such weighted norms. In addition we need to show the energies \eqref{correctweightenergy} are equivalent to the tension-dependent energies \eqref{correctenergy}, which means we have to bound $\sup_s \sigma(t,s)/s$ and $\inf_s \sigma(t,s)/s$ away from zero. (The constants in these bounds will also turn out to depend on the energies \eqref{correctweightenergy}.) Most of the work for this was done in Section \ref{greenfunctionsection}.

\subsection{Definitions and properties of weighted seminorms}

First let us define the weighted Sobolev and supremum seminorms we need.
\begin{defn}\label{weightednorms}
The weighted Sobolev seminorm of a function $f\colon [0,1] \to \mathbb{R}^d$ is defined by 
\begin{equation}\label{weightedsobolevdef}
\lVert f\rVert^2_{r, m} = \int_0^1 s^r \big\lvert f^{(m)}(s) \big\rvert^2 \, ds.
\end{equation}
The weighted supremum seminorm of $f$ is 
\begin{equation}\label{weightedsupremumdef}
\supnorm{f}^2_{r, m} = \sup_{0\le s\le 1} s^r \big\lvert f^{(m)}(s)\big\rvert^2.
\end{equation}
\end{defn}

We want to define a discrete analogue of each of these, for a sequence $\{f_1, \cdots, f_n\}$ with values in $\mathbb{R}^d$. For this purpose,
it is convenient to set 
\begin{equation}\label{risingfactorial}
s_k^{(r)} = \frac{\Gamma(k+r)}{n^r \Gamma(k)} \text{ for $k\in \{1,\ldots, n\}$ and for any real $r>-1$},
\end{equation}
where $\Gamma$ is the usual gamma function satisfying $\Gamma(x+1)=x\Gamma(x)$ for $x>0$ and $\Gamma(k) = (k-1)!$ for $k$ a natural number.
These are rising factorials, which are more convenient for our purposes than the 
falling factorials typically used in difference equations; either is much more convenient
in studying difference equations than simply using the powers $(\frac{k}{n})^r$; see \cite{LL}. Clearly
if $k_n$ is a sequence such that $\lim_{n\to\infty} \frac{k_n}{n} = s$, then we have $\lim_{n\to\infty} s_{k_n}^{(r)} = s^r$.

Recalling the definition \eqref{differenceoperator} of the difference operator $\fordiff$, we define 
our discrete analogues of \eqref{weightedsobolevdef}--\eqref{weightedsupremumdef} by
\begin{equation}\label{weightedsobolevdiscretedef}
\lVert f\rVert^2_{r,m} = \frac{1}{n} \sum_{k=1}^{n-m} s_k^{(r)} \lvert \fordiff^mf_k\rvert^2
\end{equation}
and 
\begin{equation}\label{weightedsupremumdiscretedef}
\supnorm{f}^2_{r,m} = \max_{1\le k\le n-m} s_k^{(r)} \lvert \fordiff^mf_k\rvert^2.
\end{equation}
We use the same notation for both norms to emphasize the analogy; every estimate we prove for 
the discrete seminorms \eqref{weightedsobolevdiscretedef}--\eqref{weightedsupremumdiscretedef} will have constants independent of $n$, so that we get corresponding estimates for the smooth seminorms \eqref{weightedsobolevdef}--\eqref{weightedsupremumdef}. Clearly if $f$ is smooth and we define $f_k = f(\frac{k}{n})$ for each $n$, then $\lim_{n\to\infty} \lVert f_n\rVert_{r,m} = \lVert f\rVert_{r,m}$ and $\lim_{n\to\infty} \supnorm{f_n} = \supnorm{f}$.

Now let us describe the main estimates. For two norms $\lVert \cdot\rVert_1$ and $\lVert \cdot \rVert_2$ on functions, we use the notation $ \lVert f\rVert_1 \lesssim \lVert f\rVert_2$ to mean $\lVert f\rVert \le C \lVert f\rVert$ for some constant $C$ independent of $f$. If $f$ is instead a sequence, then this notation will imply that $C$ is also independent of $n$.

For unweighted norms of smooth functions, we have the Wirtinger inequality 
$$ \int_0^1 \lvert f(s)\rvert^2 \, ds \lesssim \Big\lvert \int_0^1 f(s)\, ds\Big\rvert^2 + \int_0^1 \lvert f'(s)\rvert^2 \, ds. $$
We also have the Sobolev inequality 
\begin{equation}\label{usualsobolev}
\sup_{0\le s\le 1} \lvert f(s)\rvert^2 \lesssim \int_0^1 \lvert f(s)\rvert^2\,ds + \int_0^1 \lvert f'(s)\rvert^2\,ds.
\end{equation}
Our weighted versions of each are as follows.

\begin{theorem}\label{basicinequalitiesthm}
Let $f\colon [0,1]\to \mathbb{R}^d$ be $C^{\infty}$. Then for any
$r>0$ the norms \eqref{weightedsobolevdef} and \eqref{weightedsupremumdef} 
satisfy the weighted inequalities
\begin{align}
\lVert f\rVert^2_{r-1, m} &\lesssim \lVert f\rVert^2_{r, m} + \lVert f\rVert^2_{r+1, m+1} \text{ and} \label{weightedpoincare} \\
\supnorm{f}^2_{r, m} &\lesssim \lVert f\rVert^2_{r, m} + \lVert f\rVert^2_{r+1, m+1}. \label{weightedsobolev}
\end{align}
If in addition we have $f^{(m)}(1)=0$, then these inequalities can be simplified to 
\begin{align}
\lVert f\rVert^2_{r-1,m} &\lesssim \lVert f\rVert^2_{r+1,m+1} \text{ and}\label{tieddownpoincare} \\
\supnorm{f}^2_{r,m} &\lesssim \lVert f\rVert^2_{r+1,m+1}. \label{tieddownsobolev}
\end{align}

If $f$ is instead a sequence $\{f_1,\cdots, f_n\}$ with values in $\mathbb{R}^d$, then the 
inequalities \eqref{weightedpoincare}--\eqref{weightedsobolev} also hold if the norms are 
interpreted as \eqref{weightedsobolevdiscretedef} and \eqref{weightedsupremumdiscretedef},
while the inequalities \eqref{tieddownpoincare}--\eqref{tieddownsobolev} hold if $f_{n-m}^{(m)}=0$. 
\end{theorem}

\begin{proof}
The continuous version of this inequality appears in Adams-Fournier~\cite{adams}. We prove the discrete
version in Appendix \ref{A2}, from which the continuous version follows in the limit.
\end{proof}

\begin{rem}\label{nosup}
The example $f(s)=\arcsinh{(\ln{s})}$ demonstrates that the inequalities
\eqref{tieddownpoincare} and \eqref{tieddownsobolev} cannot be extended to
$r=0$: in that case we have $f(1)=0$,
$\int_0^1 \lvert f(s)\rvert^2 \,ds< \infty$, and  $\int_0^1 s\lvert
f'(s)\rvert^2 \, ds < \infty$, while $\int_0^1
\frac{1}{s} \lvert f(s)\rvert^2 \, ds$ and $\sup_{x\in [0,1]} \lvert
f(x)\rvert^2$ are both infinite.
In particular there cannot be constants for the discrete versions that are independent of $n$ when $r=0$.
\end{rem}

The important thing about \eqref{weightedsobolev} is that by Remark \ref{nosup}, the estimate only works when $r>0$. Hence in any computation where a supremum norm is required, we will want a positive power of $s$ attached to be able to use this result. This will show up 
when we need to estimate weighted Sobolev norms of products of three functions: we want to pull out a supremum of one and use Cauchy-Schwarz on the rest,
and we will need a little extra weighting in some cases. Of course, we could use the usual Sobolev inequality \eqref{usualsobolev} to get 
$$ \supnorm{f}_{0,0} \lesssim \lVert f\rVert_{0,0} + \lVert f\rVert_{0,1} \lesssim \lVert f\rVert_{0,0} + \lVert f\rVert_{1,1} + \lVert f\rVert_{2,2},$$
but requiring \emph{two} extra derivatives rather than one is usually not worthwhile (except once in the proof of Theorem \ref{uniqueness}).

Frequently our weighting in discrete norms will be slightly off (for example, we may want to replace 
$s_k^{(p+q)}/s_k^{(q)}$ with $s_k^{(p)}$, or we may want to replace $s_k^{(p)}$ with 
$s_{k+j}^{(p)}$ for some $j$). In the continuous case these formulas are trivial, but in the discrete case, 
bounds such as these come from properties of the gamma function (in particular the fact that
the gamma function is log-convex by the Bohr-Mollerup theorem).
The constants will never be important; what will matter is that they are independent of $k$ and $n$. 
The following estimates are easy to prove.

\begin{proposition}\label{weightsprop}
Let $n\in\mathbb{N}$ and let $k\in\{1,\ldots, n\}$. Let $p$ and $q$ be
positive real numbers. 

Then the weight function $s_k^{(p)} = \frac{\Gamma(k+p)}{n^p\Gamma(k)}$ satisfies
the following inequalities:
\begin{equation}\label{w}
s_k^{(p)} \le \frac{s_k^{(p+q)}}{s_k^{(q)}} \le \frac{\Gamma(p+q+1)}{\Gamma(p+1)\Gamma(q+1)} \, s_k^{(p)}.
\end{equation}
We also have $s_k^{(p)} \le s_{k+j}^{(p)} \le \frac{\Gamma(j+p+1)}{\Gamma(j+1)\Gamma(p+1)} \, s_k^{(p)}$ for any nonnegative integer $j$.
\end{proposition}

Proposition \ref{weightsprop} also gives the following corollary, which is the most useful tool we have for estimating norms of products. To get the higher-difference norms of products, we will use the product rule \eqref{doubleproduct} for differences together with these formulas. The proof is trivial.

\begin{corollary}\label{supproductboundcorollary}
Suppose $(f_1, \ldots, f_n)$ and $(g_1, \ldots, g_n)$ are sequences of real numbers, and let $p$ and $q$ be nonnegative real numbers. Then 
\begin{equation}\label{supregproductbound}
\lVert fg\rVert^2_{p+q, 0} \lesssim \supnorm{f}^2_{p, 0} \lVert g\rVert^2_{q, 0}
\end{equation}
and 
\begin{equation}\label{supsupproductbound}
\supnorm{fg}^2_{p+q, 0} \lesssim \supnorm{f}^2_{p,0} \supnorm{g}^2_{q,0}.
\end{equation}
The formulas are also valid if one of the sequences is in $\mathbb{R}$ and the other in $\mathbb{R}^d$, or if both are in $\mathbb{R}^d$ and we use $\langle f, g\rangle$ or $\lvert f\rvert \lvert g\rvert$ instead.

Estimates \eqref{supregproductbound} and \eqref{supsupproductbound} are also valid if $f$ and $g$ are smooth functions with the norms interpreted as \eqref{weightedsobolevdef} and \eqref{weightedsupremumdef}.
\end{corollary}

\begin{rem}
Typically we will extend the seminorms \eqref{weightedsobolevdiscretedef} and \eqref{weightedsupremumdiscretedef} when used for $\eta$ and $\sigma$; since we have $\sigma_0=0$ and $\eta_{n+1}=0$, it is more convenient to modify the definitions to 
\begin{align*}
\lVert \eta\rVert^2_{r,m} &= \frac{1}{n} \sum_{k=1}^{n-m+1} s_k^{(r)} \lvert (\fordiff^m\eta)_k\rvert^2 \\
\lVert \sigma\rVert^2_{r,m} &= \frac{1}{n} \sum_{k=0}^{n-m} s_k^{(r)} \lvert (\fordiff^m\sigma)_k\rvert^2.
\end{align*}
This does not affect any of the estimates, but it allows us to incorporate the endpoint information. This is convenient 
for example to interpret \eqref{finitetensionsolution} in terms of $\lVert \dot{\eta}\rVert^2_{1,1}$ in Lemma \ref{acboundslemma}. 
\end{rem}

\subsection{Weighted energy norms}

We have already defined the weighted energy \eqref{correctenergy} and \eqref{correctweightenergy} for a whip. Now we want to define corresponding discrete energy  for the chain. The definitions are made much easier
if we use the odd extension \eqref{discreteodd} of $\eta$ to define the differences $\fordiff^m\eta_k$ beyond 
$k=n-m+1$. Furthermore, by analogy with \eqref{risingfactorial}, we define
\begin{equation}\label{risingfactorialsigma}
\sigma_k^{(r)} = \prod_{j=k}^{k+r-1} \sigma_j \text{ for any integer $r\ge 0$.}
\end{equation}

Our time-independent energy will be
\begin{equation}\label{discreteenergydef}
e_m = \frac{1}{n} \sum_{\ell=0}^m \sum_{k=1}^{n-\lfloor\ell/2\rfloor} \Big(s_k^{(\ell)} \lvert \fordiff^{\ell}\dot{\eta}_k\rvert^2 + s_k^{(\ell+1)} \lvert \fordiff^{\ell+1}\eta_k\rvert^2\Big),
\end{equation}
while the time-dependent energy is 
\begin{equation}\label{discretetimeenergydef}
\widetilde{e}_m = \frac{1}{n} \sum_{\ell=0}^m \sum_{k=1}^{n-\lfloor \ell/2\rfloor} \Big( \sigma_k^{(\ell)} \lvert \fordiff^{\ell}\dot{\eta}_k\rvert^2
+ \sigma_k^{(\ell+1)} \lvert \fordiff^{\ell+1}\eta_k\rvert^2\Big).
\end{equation}
Recall that we need to use the time-dependent $\sigma$-weighted quantities to compute the time-derivative of energy
in order to get some cancellation, while only time-independent energies are useful for constructing topologies and 
relating distinct norms.

Clearly if we have sequences $\eta_n$ and $\sigma_n$ defined for each $n\in\mathbb{N}$ as in \eqref{discretizationprime}, then 
\begin{equation*}
E_m[\eta] = \lim_{n\to\infty} e_m[\eta_n]
\quad \text{and}\quad 
\widetilde{E}_m[\eta] = \lim_{n\to\infty} \widetilde{e}_m[\eta_n].
\end{equation*}
So any estimate we obtain on the chain energies $e_m$ and $\widetilde{e}_m$ will become an \emph{a priori} estimate 
on the corresponding whip energies $E_m$ and $\widetilde{E}_m$.

Note that we have 
\begin{equation}\label{norminequality}
e_m \ge \sum_{\ell=0}^m \lVert \dot{\eta}\rVert^2_{\ell, \ell} + \sum_{\ell=1}^{m+1} \lVert \eta\rVert^2_{\ell,\ell},
\end{equation}
in terms of the discrete weighted seminorms \eqref{weightedsobolevdiscretedef}.\footnote{We would have equality if the sums over $k$ went from $k=1$ to $k=n-\ell$ rather than $k=n-\lfloor\ell/2\rfloor$. The reason the sums in \eqref{discreteenergydef} and \eqref{discretetimeenergydef} contain a few extra terms in the sums is in order to make the derivative estimate of Theorem \ref{discreteenergyestimate} simpler: with this definition the endpoint terms of the discrete energy derivative always vanish, as they did in \eqref{nineandahalf}.}

It is also convenient to observe that the lowest-level energy $e_0$ is constant in time. Since it will be useful later in Lemma \ref{bboundslemma}, we separate the terms and define
\begin{equation}\label{u0def}
u_0 = \frac{1}{n} \sum_{k=1}^n \lvert \dot{\eta}_k\rvert^2 = \lVert \dot{\eta}\rVert^2_{0,0}, \qquad v_0 = \frac{1}{n} \sum_{k=1}^n s_k \lvert \fordiff\eta_k\rvert^2.
\end{equation}
\begin{lemma}\label{u0v0constant}
If $\eta$ satisfies \eqref{fancydiscreteevolution} with $\lvert \fordiff\eta\rvert\equiv 1$, then both $u_0$ and $v_0$ are constant in time.
\end{lemma}

\begin{proof}
The fact that $\lvert \fordiff\eta\rvert \equiv 1$ implies that $v_0 = \frac{1}{2} + \frac{1}{2n}$.
As a corollary, every energy $e_m$ given by \eqref{discreteenergydef} satisfies $e_m \ge \frac{1}{2}$. 

For $u_0$, we just compute 
$$
\frac{du_0}{dt} = \frac{2}{n} \sum_{k=1}^n \langle \dot{\eta}_k, \ddot{\eta}_k\rangle 
= \frac{2}{n} \sum_{k=1}^n \langle \dot{\eta}_k, \backdiff(\sigma \fordiff\eta)_k \rangle.
$$
Using the summation by parts formula \eqref{summationbyparts} along with $\langle \backdiff \eta, \backdiff\dot{\eta}\rangle \equiv 0$ and the endpoint conditions $\sigma_0=0$ and $\eta_{n+1}=0$, it is easy to show this sum vanishes.
\end{proof}

For a smooth solution of \eqref{generalstringevol}--\eqref{generalstringconst}, we clearly have that the analogous quantities $U_0 = \int_0^1 \lvert \etadifft(t,s)\rvert^2 \, ds$ and $V_0 = \int_0^1 s \lvert \etadiffs(t,s)\rvert^2 \, ds$ satisfy $V_0=\frac{1}{2}$ and $U_0$ is constant in time. Thus $E_0$ is also constant in time.

Our primary use of Theorem \ref{basicinequalitiesthm} will be the following formulas, which follow easily from 
\eqref{norminequality}.
\begin{lemma}\label{normtoenergy}
For any $i\ge 0$ and $0\le j<i$, we have
\begin{alignat}{2}
\lVert \eta\rVert^2_{i-j, i} &\lesssim e_{i+j-1}, 
& \lVert \dot{\eta}\rVert^2_{i-j, i} &\lesssim e_{i+j}, \label{squaredotenergybound}\\
\supnorm{\eta}^2_{i-j,i} &\lesssim e_{i+j}, 
& \supnorm{\dot{\eta}}^2_{i-j,i} &\lesssim e_{i+j+1} \label{supenergybound} \\
&\supnorm{\dot{\eta}}^2_{1/2, i} \lesssim e_{2i+1}. \label{supdotenergybound} &&
\end{alignat}
\end{lemma}

\begin{rem}
We need the extra power of $1/2$ in estimate \eqref{supdotenergybound}, since \eqref{weightedsobolev} is not valid when $r=0$. This is important once at the end of the proof of Lemma \ref{sigmasobolevlemma} and once at the end of the proof of Theorem \ref{discreteenergyestimate}.
\end{rem}

\section{Bounds for the tension in terms of the energy}\label{tensionboundssection}

Before bounding the energy itself, we first want bounds for the tension $\sigma$ given by either \eqref{generalstringconst} or \eqref{fancydiscretetension}. To compare the energies $\tilde{E}_m$ and $E_m$, we want upper and lower bounds for $\sigma/s$. In addition, to compute the time derivative of $\tilde{E}_m$, we need to know a bound for $\sigmadifft/s$, by formula \eqref{nineandahalf}.

For a smooth solution $(\eta, \sigma)$ of \eqref{generalstringevol}--\eqref{generalstringconst}, 
we define quantities $A$, $B$, and $C$ by the formulas
\begin{equation}\label{ABCdef}
A(t) = \sup_{0\le s\le 1} \lvert \sigmadiffs(t,s)\rvert, \quad 
B(t) = \sup_{0\le s\le 1} \frac{s}{\sigma(t,s)}, \quad C(t) = \sup_{0\le s\le 1} \lvert \sigmadiffst(t,s)\rvert.
\end{equation}
Observe that since $\sigma(t,0)=0$, we have 
$$ \lvert \sigma(t,s)\rvert = \left\lvert \int_0^s \sigma_x(t,x) \, dx \right\rvert \le s\sup_{0\le x\le 1} \lvert \sigma_x(t,x)\rvert,$$
so that we have 
\begin{equation}\label{supfracderiv}
\sup_{0\le s\le 1} \frac{\sigma(t,s)}{s} \le A(t), \quad \text{and similarly} \quad \sup_{0\le s\le 1} \frac{\sigmadifft(t,s)}{s} \le C(t).
\end{equation}
Generally the bounds \eqref{supfracderiv} will be much more useful to us, although occasionally we will need the actual definition \eqref{ABCdef}.

Similarly, we define the discrete analogues of the quantities \eqref{ABCdef}; as with the energy, we use upper-case and lower-case for norms of the whip or chain respectively. Recall that $s_k=\frac{k}{n}$, while our convention is that $\sigma_0(t)=0$, and recall the definition $(\backdiff\sigma)_k = n(\sigma_k-\sigma_{k-1})$. We therefore set
\begin{equation}\label{abcdef}
a(t) = \max_{1\le k\le n} \lvert (\backdiff\sigma)_k(t)\rvert, \quad b(t) = \max_{1\le k\le n} \frac{s_k}{\sigma_k(t)}, \quad 
c(t) = \max_{1\le k\le n} \lvert (\backdiff\dot{\sigma})_k(t)\rvert.
\end{equation}
As above, the fact that $\sigma_0(t)=0$ means we can write 
$$ \lvert \sigma_k(t)\rvert = \left\lvert \frac{1}{n} \sum_{j=1}^k (\backdiff \sigma)_j(t) \right\rvert \le \frac{k}{n} \left( \max_{1\le j\le n} \lvert (\backdiff\sigma)_j(t)\rvert\right)$$ to obtain
\begin{equation}\label{supfracnabla}
\max_{1\le k\le n} \frac{\sigma_k}{s_k} \le a \qquad \text{and}\qquad \max_{1\le k\le n} \frac{\lvert \dot{\sigma}_k\rvert}{s_k} \le c,
\end{equation}

\begin{lemma}\label{acboundslemma}
If $\eta$ and $\sigma$ form a smooth solution of \eqref{generalstringevol}--\eqref{generalstringconst}, then the quantities defined by \eqref{ABCdef} satisfy the bounds
\begin{align}
A(t) &\lesssim E_2(t), \label{Abound} \\
C(t) &\lesssim E_2(t)^{3/2} E_3(t)^{1/2}. \label{Cbound}
\end{align}

Similarly, suppose $(\eta_1(t), \ldots, \eta_n(t))$ and $(\sigma_1(t), \ldots, \sigma_n(t))$ form a solution of \eqref{fancydiscreteevolution} and \eqref{fancydiscretetension}
with the odd extensions \eqref{discreteodd}. Suppose also that we have $\alpha_i = \langle \fordiff\eta_i, \fordiff\eta_{i+1}\rangle \ge 0$ for $1\le i\le n-1$.                                                                                                                                                                                                                                                                                                                                                                                                                                                                                                                                                                                                                             
Then $a$ and $c$ given by \eqref{abcdef}
satisfy the bounds 
\begin{align}
a(t) &\lesssim e_2(t) \label{abound} \\
c(t) &\lesssim e_2(t)^{3/2} e_3(t)^{1/2}.\label{cbound}
\end{align}
\end{lemma}

\begin{proof}
We will just prove the discrete bounds for $a$ and $c$ in detail; the bounds for $A$ and $C$
can be proved using the same techniques, or we can view them as a limiting case of the bounds for $a$ and $c$. 

The estimate \eqref{abound} for $a$ is easy: by Proposition \ref{solverprop}, we have 
$$ (\backdiff\sigma)_k = \frac{1}{n} \sum_{j=1}^n (\backdiffone G)_{kj} \lvert \fordiff\dot{\eta}_j\rvert^2.$$
By Proposition \ref{greenupperbounds}, we have $\lvert (\backdiffone G)_{kj}\rvert \le 1$, and thus 
$$
(\backdiff\sigma)_k \le \frac{1}{n} \sum_{j=1}^n \lvert \fordiff\dot{\eta}_j\rvert^2
= \lVert \dot{\eta}\rVert^2_{0,1} 
\lesssim e_2$$
by \eqref{squaredotenergybound}. We then have \eqref{abound}. The inequality \eqref{Abound} is proved identically.

The bound \eqref{cbound} for $c$ is more complicated. The first step is to differentiate the equation
\eqref{fancydiscretetension} in time to get 
\begin{align*}
-\langle \fordiff\eta, \backdiff\fordiff (\dot{\sigma} \fordiff \eta)\rangle &= 
\langle \fordiff\dot{\eta}, \backdiff\fordiff (\sigma \fordiff\eta)\rangle + 
\langle \fordiff\eta, \backdiff\fordiff (\sigma \fordiff\dot{\eta})\rangle \\
&\qquad\qquad
+ 2 \langle \fordiff \dot{\eta}, \fordiff\ddot{\eta} \rangle \\
&= 3 \langle \fordiff\dot{\eta}, \backdiff\fordiff (\sigma\fordiff\eta)\rangle + \langle \fordiff\eta, \backdiff\fordiff(\sigma\fordiff\dot{\eta})\rangle.\end{align*}
Thus $\dot{\sigma}$ satisfies the same kind of equation as $\sigma$ with the same endpoint conditions, so we can use Proposition \ref{solverprop} to write
$$(\backdiff\dot{\sigma})_k = \frac{1}{n} \sum_{j=1}^n (\backdiffone G)_{kj} \Big( 3 \langle \fordiff\dot{\eta}_j, \backdiff\fordiff (\sigma\fordiff\eta)_j\rangle + \langle \fordiff\eta_j, \backdiff\fordiff(\sigma\fordiff\dot{\eta})_j\rangle\Big).$$
As above, the fact that $\lvert (\backdiffone G)_{kj}\rvert\le 1$ implies that $c\le \lvert \Lambda\rvert$, where
\begin{equation}\label{cbound1}
\Lambda = \frac{3}{n} \sum_{j=1}^n \Big(\langle \fordiff\dot{\eta}_j, \backdiff\fordiff (\sigma\fordiff\eta)_j\rangle + \langle \fordiff\eta_j, \backdiff\fordiff(\sigma\fordiff\dot{\eta})_j\rangle\Big).
\end{equation}

Applying the summation by parts formula \eqref{summationbyparts} to this, using the endpoint
conditions $\eta_{n+1}=\dot{\eta}_{n+1}=\sigma_0$, and performing some manipulations with the formula
$\langle \fordiff\eta, \fordiff\dot{\eta}\rangle\equiv 0$, 
we get
$$ \Lambda = \frac{1}{n} \sum_{j=1}^n \Big[ \fordiff\sigma_j \Big(3\langle \fordiff\dot{\eta}_j, \fordiff^2\eta_j\rangle - \langle \fordiff^2\eta_j, \fordiff\dot{\eta}_{j+1}\rangle\Big)
- 4 \sigma_j \langle \fordiff^2\eta_j, \fordiff^2\dot{\eta}_j\rangle\Big].$$

Using the bounds $\lvert \fordiff\sigma\rvert \lesssim e_2$ and $\frac{\sigma_k}{s_k} \lesssim e_2$ from above, we obtain

$$c \le \lvert \Lambda\rvert 
\le e_2 \lVert \dot{\eta}\rVert_{0,1} \lVert \eta\rVert_{0,2} + e_2 \lVert \eta\rVert_{1,2} \lVert \dot{\eta}\rVert_{1,2}
\lesssim e_2 \sqrt{e_2} \sqrt{e_3},
$$
using Lemma \ref{normtoenergy}, which gives \eqref{cbound}. The proof of \eqref{Cbound} is almost identical.
\end{proof}

\begin{lemma}\label{bboundslemma}
Suppose $\eta$ and $\sigma$ solve \eqref{generalstringevol}--\eqref{generalstringconst}. Then $B$ defined by \eqref{ABCdef} satisfies the estimate 
\begin{equation}\label{Bbound}
B(t) \lesssim \frac{1}{U_0}(1+\lVert \eta\rVert^2_{1,2} ) e^{\lVert \eta\rVert^2_{1,2}},
\end{equation}
where $U_0 = \int_0^1 \lvert \etadifft\rvert^2 \, ds$ is a constant depending on the initial condition.

Similarly, suppose $(\eta_1(t), \ldots, \eta_n(t))$ and $(\sigma_1(t), \ldots, \sigma_n(t))$ form a solution of \eqref{fancydiscreteevolution} and \eqref{fancydiscretetension}
with the odd extensions \eqref{discreteodd}. 
Suppose also that $\supnorm{\eta}_{3/2,2} \le \frac{2\sqrt{n}}{5}$. 
Then the quantity $b$ in \eqref{abcdef} satisfies the estimate 
\begin{equation}\label{bbound}
b(t) \lesssim \frac{1}{u_0} e^{2\supnorm{\eta}^2_{3/2,2}},
\end{equation}
where $u_0$ is defined by \eqref{u0def} (and is constant due to Lemma \ref{u0v0constant}).
\end{lemma}

\begin{proof}
The bound \eqref{Bbound} comes directly from \eqref{lowerboundcont}: we have by Proposition \ref{sigmasolprop} that
\begin{align*}
\inf_{0\le s\le 1} \frac{\sigma(s)}{s} \ge \int_0^1 \inf_{0\le s,x\le 1} \frac{G(s,x)}{sx}\cdot x \lvert \etadifftx(t,x)\rvert^2 \, dx \ge \frac{e^{-\lVert \eta\rVert^2_{1,2}}}{1+\lVert \eta\rVert^2_{1,2}} \, \lVert \etadifft\rVert^2_{1,1}.
\end{align*}
The fact that $U_0 = \lVert \etadifft\rVert^2_{0,0} \lesssim \lVert \etadifft\rVert^2_{1,1}$ follows from
\eqref{tieddownpoincare}, since $\etadifft(t,1)=0$. 

The proof of \eqref{bbound} is similar, using \eqref{lowerbounddisc} and \eqref{finitetensionsolution} for the discrete Green function. 
The bound $u_0\lesssim \lVert \dot{\eta}\rVert^2_{1,1}$ similarly follows from \eqref{tieddownpoincare} since $\dot{\eta}_{n+1}\equiv 0$. 
\end{proof}

\begin{exmp}\label{ABCnecessary}
In terms of the weighted energy \eqref{correctweightenergy}, Lemmas \ref{acboundslemma} and \ref{bboundslemma} give upper bounds for $A$ and $B$ if the energy $E_2$ is finite, while we only get an upper bound for $C$ if  $E_3$ is finite. Simple examples show that these conditions are necessary: we can have $E_1$ bounded while $A$ and $B$ are unbounded, and we can have $E_2$ bounded while $C$ is unbounded. 

To obtain the examples for $A$ and $B$, we consider (at time $t=0$) the whip position
\begin{equation}\label{singularcurvature}
\eta(0,s) = \left( \frac{3(1-s\cos{ (\tfrac{2}{3}\ln{s})  } )}{\sqrt{13}}, -\frac{3s\sin{(\tfrac{2}{3}\ln{s})}}{\sqrt{13}}\right),
\end{equation}
which satisfies $\eta(0,1)=0$, $\lvert \etadiffs\rvert \equiv 1$, and $\lvert \etadiffss\rvert = \frac{2}{3s}$. 
This corresponds to a whip where the free end sits at $(\frac{3}{\sqrt{13}}, 0)$ despite making infinitely many rotations around it (as $s\to 0$). See Figure \ref{logarithmicspiral}.

\begin{figure}[ht]
\begin{center}
\includegraphics[scale=0.5]{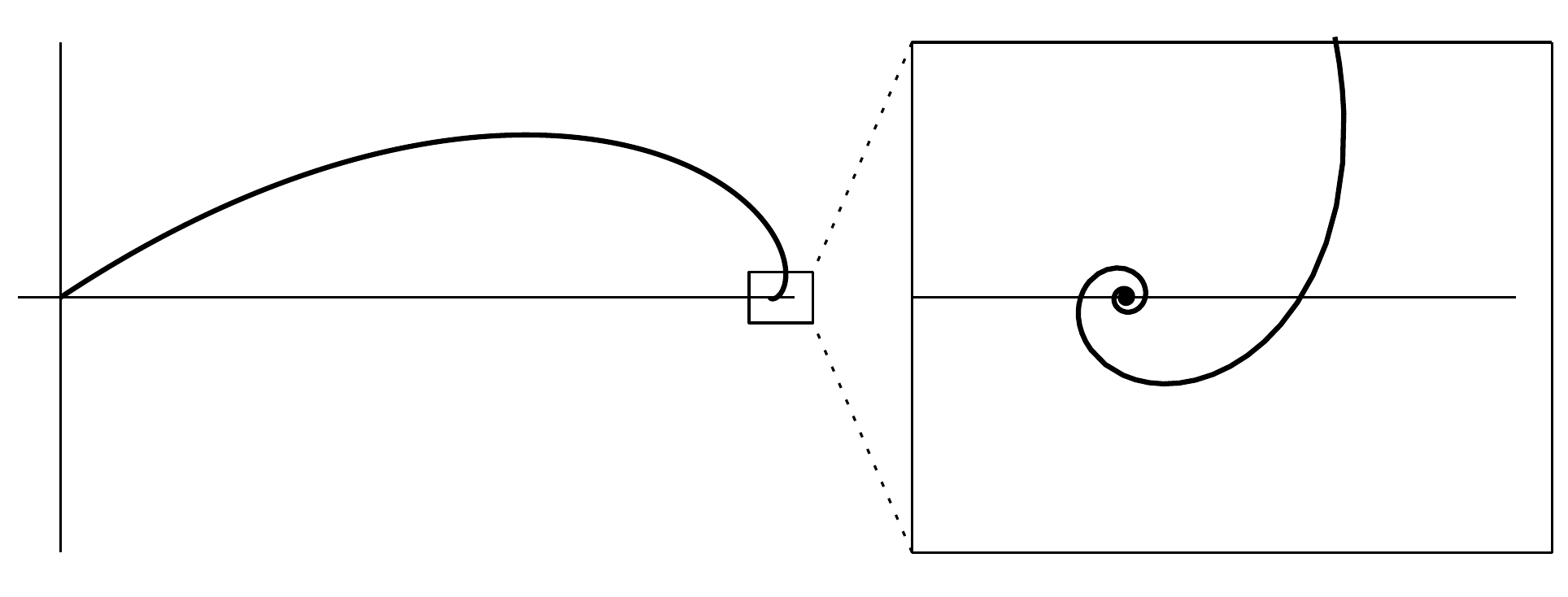}
\caption{The curve defined by \eqref{singularcurvature}, for which the curvature approaches infinity at the free end. Although the curve has length one, its free end wraps around the limiting point $\frac{3}{\sqrt{13}}$ infinitely many times. We have plotted this heuristically in the inset, although the actual curve wraps itself up too tightly for these loops to be visible.}\label{logarithmicspiral}
\end{center}
\end{figure}

The Green function \eqref{greencontinuous} can be computed explicitly to obtain $\sigma$ from $\lvert \etadiffst\rvert$. If $\eta(0,s)$ satisfies \eqref{singularcurvature} and $\lvert \etadiffst(0,s)\rvert = s^{-3/4}$, one computes that $E_1$ is finite while $E_2$ and $A$ are both infinite. If on the other hand $\lvert \etadiffst(0,s)\rvert \equiv 1$, we easily see that $E_1$ is still finite while $E_2$ and $B$ are both infinite.

The example for $C$ is a bit more involved. Suppose $\eta$ is given by \eqref{spherical}, where
$\thetadiffs(0, s) = s^{-3/4}$ and $\thetadifft(0, s) = s^{-1/4}$. It is easy to verify that $E_2$ is finite at this instant, while $E_3$ is infinite. We have $\sigma(s)=s$ at this instant by \eqref{generalstringconst}, so that 
differentiating \eqref{generalstringconst} with respect to time and using \eqref{generalstringevol} gives
$ \sigmadifftss(s) - s^{-3/2} \sigmadifft(s) = -3/s$ with boundary conditions $\sigmadifft(0)=0$ and $\sigmadiffst(1)=0$. In this case we can verify that $C$ is infinite.
\end{exmp}

The fact that we cannot bound $C$ unless $E_3$ is bounded is one of the main reasons why the energy estimates only close up at $E_3$. Taking a time derivative of \eqref{correctenergy} as in \eqref{nineandahalf} gives a number of terms of the form $\sigmadifft$ which can only be bounded in terms of $C$, and thus in terms of $E_3$.

We are now ready for an a priori estimate on the tension $\sigma$. Although the norms of $\sigma$ and $\sigmadiffs$ are easier to measure using the supremum, it is convenient to use weighted Sobolev norms for the higher derivatives of $\sigma$. Thus we define the squared norm for a whip:
\begin{equation}\label{sigmasobolevnormsmooth}
D_m = \sum_{\ell=0}^{m-1} \lVert \sigma\rVert^2_{\ell+3/2, \ell+2} \quad \text{for $m\ge 1$}.
\end{equation}
The discrete version is defined by the same formula:
\begin{equation}\label{sigmasobolevnorm}
d_m = \sum_{\ell=0}^{m-1} \lVert \sigma\rVert^2_{\ell+3/2, \ell+2} \quad \text{for $m\ge 1$}.
\end{equation}

\begin{lemma}\label{sigmasobolevlemma}
If $\sigma$ is a smooth solution of \eqref{generalstringconst}, then the norms \eqref{sigmasobolevnormsmooth} can 
be bounded by the energy \eqref{correctweightenergy} via
\begin{equation}\label{lowdsmooth}
D_1\lesssim E_3^4, \qquad D_2\lesssim E_3^4, \qquad \text{and} \quad D_3\lesssim E_3^6,
\end{equation} 
while for $m>3$ we have 
\begin{equation}\label{sigmaboundsmooth}
D_m \le P_m(E_{m-1}) E_m,
\end{equation}
where $P_m$ depends only on $E_{m-1}$. 

Similarly if $\sigma$ satisfies \eqref{usefultension} with the condition $\sigma_0=0$, then the norms \eqref{sigmasobolevnorm} can be bounded by the energy \eqref{discreteenergydef} via
\begin{equation}\label{lowd}
d_1\lesssim e_3^4, \qquad d_2\lesssim e_3^4, \qquad \text{and} \quad d_3\lesssim e_3^6,
\end{equation} 
while for $m>3$ we have 
\begin{equation}\label{sigmabound}
d_m \le P_m(e_{m-1}) e_m,
\end{equation}
where $P_m$ depends only on $e_{m-1}$. 
\end{lemma}

\begin{proof}
We will just prove the discrete estimates \eqref{lowd}--\eqref{sigmabound}; the estimates \eqref{lowdsmooth}--\eqref{sigmaboundsmooth} are proved using the exact same technique. The full proof is in Appendix \ref{A3}; the basic idea is just to take iterated differences of \eqref{usefultension} and estimate using Corollary \ref{supproductboundcorollary} and Lemma \ref{normtoenergy}.
\end{proof}

\section{The main energy estimate}\label{aprioriestimates}\label{energyestimatesection}

In order to construct the solution of the partial differential equations
\eqref{basicPDE}--\eqref{tensionODE}, we want to find bounds on all the discrete energies \eqref{discreteenergydef}--\eqref{discretetimeenergydef}
which are independent of the initial conditions and of the number $n$ of links.
Then in Section \ref{convergence} we will find a subsequence that converges
to a solution.
As a consequence, we can show that the motion of a chain converges to
the motion of a whip as $n$ approaches infinity, in the sense that
position, velocity, and acceleration all converge.

We now want to estimate the time evolution of the energy $\tilde{e}_m$. Our strategy will be to bound $d\tilde{e}_m/dt$ in terms of the energies $e_m$; we will then use the fact that $e_m$ and $\tilde{e}_m$ are equivalent (since $\sigma/s$ is bounded above and below by Lemmas \ref{acboundslemma} and \ref{bboundslemma}) to get an inequality for $d\tilde{e}_m/dt$ in terms of $\tilde{e}_m$.
In proving it we will use 
Lemmas \ref{acboundslemma}--\ref{bboundslemma} and \ref{sigmasobolevlemma} in an essential way.

\begin{theorem}\label{discreteenergyestimate}
Let $n\in \mathbb{N}$, and suppose $(\eta_1(t), \ldots, \eta_n(t))$ and $(\sigma_1(t), \ldots, \sigma_n(t))$ form a solution of \eqref{fancydiscreteevolution} and \eqref{fancydiscretetension}
with $\sigma_0(t)\equiv 0$ and $\eta_{n+1}(t)\equiv 0$, and that $\eta$ and $\sigma$ extend to sequences  satisfying the oddness condition \eqref{discreteodd}. 

Then the energies \eqref{discreteenergydef} and \eqref{discretetimeenergydef} satisfy the estimates
\begin{equation}\label{energyderivative}
\frac{d\widetilde{e}_3}{dt} \le M_3 e_3^7
\end{equation}
for some $M_3$ independent of the initial data and of $n$.
In addition the higher energies satisfy 
\begin{equation}\label{generalenergyderivative}
\frac{d\widetilde{e}_m}{dt} \lesssim M_m(e_{m-1}) e_m
\end{equation}
for every $m>3$, where $M_m$ depends only on $e_{m-1}$. 

Analogously, if $\eta$ and $\sigma$ form a smooth solution of \eqref{generalstringevol} and \eqref{generalstringconst}, then the energies \eqref{correctenergy} and \eqref{correctweightenergy} satisfy the estimates
\begin{equation}\label{energyderivativesmooth}
\frac{d\widetilde{E}_3}{dt} \le M_3 E_3^7
\end{equation}
for some $M_3$ independent of the initial data.
In addition the higher energies satisfy 
\begin{equation}\label{generalenergyderivativesmooth}
\frac{d\widetilde{E}_m}{dt} \lesssim M_m(E_{m-1}) E_m
\end{equation}
for every $m>3$, where $M_m$ depends only on $E_{m-1}$. 
\end{theorem}

\begin{proof}
The proof is in Appendix \ref{A4}.
\end{proof}

The fact that the energy estimates only close up at $m=3$ is perhaps
explained by the following observation, which is easier to understand in terms
of the spherical representation \eqref{spherical}.

\begin{proposition}\label{r4}
Let $D^4$ denote the unit ball in $\mathbb{R}^4$. The pair $(\theta, \sigma)$ is a smooth solution of \eqref{thetaevol}--\eqref{thetaconst}
if and only if the functions $ \varphi\colon D^4 \to
S^1 $ and $\alpha \colon D^4 \to \mathbb{R}^+$ defined by
$$ 
\varphi(x) = \theta(\lvert x\rvert^2) \quad \text{and}\quad
\alpha(x) = \frac{\sigma(\lvert x\rvert^2)}{4\lvert x\rvert^2}
$$
are spherically symmetric solutions of the equations
\begin{equation}\label{r4equations}
\begin{split}
\phidifftt &= \diver{(\alpha \grad \varphi)} \\
\Laplacian \alpha - \lvert \grad \varphi\rvert^2 \alpha &= -\lvert
\phidifft \rvert^2.
\end{split}
\end{equation}
with Neumann boundary condition $\partial_{\nu} \varphi = 0$ 
for $\varphi$ and Robin boundary condition $
\partial_{\nu}\alpha
+ 2 \alpha = 0$ for $\alpha$ on $\partial D^4
= S^3$. Furthermore any smooth solution has $\alpha>0$ everywhere,
so that the hyperbolic equation for $\varphi$ is nondegenerate.
\end{proposition}

\begin{proof}
Setting $\sigma(s) = 4s\alpha(s)$ and changing variables by $s=r^2$, 
we easily see that \eqref{thetaevol}--\eqref{thetaconst} become
\begin{align*}
\phidifftt &= \alpha\left( \phidiffrr + \tfrac{3}{r}
\phidiffr \right) + 2 \alphadiffr \phidiffr \\
-\lvert \phidifft \rvert^2 &= \alphadiffrr + \tfrac{3}{r} \alphadiffr -
\lvert \phidiffr \rvert^2 \alpha.
\end{align*}
Now the operator $\partial_r^2 + \frac{3}{r} \partial_r$ is familiar as the spherically
symmetric Laplacian on $\mathbb{R}^4$, and hence we recognize both
terms above as coming from the Laplacian on $\mathbb{R}^4$ under the
assumption that $\alpha$ and $\theta$ are both spherically
symmetric. The boundary conditions are easy to check.
\end{proof}

The fact that the degeneracy can be removed if we work in a
higher-dimensional space, and thus in some sense the equations
naturally ``live'' there, is essentially the reason why we need
higher than usual Sobolev order for the estimates to close.

\section{Local existence and uniqueness of the solution}\label{convergence}\label{existenceuniquenesssection}

Now we can finally prove the local existence theorem for the system
\eqref{generalstringevol}--\eqref{generalstringconst} of partial differential equations. The fact that Theorem
\ref{discreteenergyestimate} gives us estimates for $\tilde{e}_m$ in terms of $e_m$ 
that are independent of $n$ allows us to construct the solution as a
limit of a subsequence of discrete solutions as $n\to \infty$,
following the technique of Ladyzhenskaya~\cite{lady} and references
therein.

\subsection{The discrete interpolation}\label{interpolation}

We first need to establish the interior approximation of the space of whips by the space of chains,
which allows us to go from estimates on $e_m$ given by \eqref{discreteenergydef} to estimates on 
$E_m$ given by \eqref{correctweightenergy} and back.

Consider any function $\eta\colon [0,1]\to \mathbb{R}^d$ such that $\lvert \etadiffs\rvert\equiv 1$,
with $\eta$ extending to an odd function through $s=1$, such that the seminorms $\lVert \eta\rVert_{\ell, \ell}$
for $2\le \ell \le m$ are all finite. For each $n\in\mathbb{N}$ we want to approximate $\eta$ by a sequence 
$\eta_k \in \mathbb{R}^d$ for $1\le k\le n$, extend it for $k>n$ by $\eta_{k} = -\eta_{2n+2-k}$, have it satisfy $\lvert \fordiff\eta_k\rvert \equiv 1$, and have uniform bounds on the discrete Sobolev seminorms $\lVert \eta\rVert_{\ell,\ell}$ in terms of the smooth seminorms 
that are independent of $n$.

The complication arises from handling the constraint $\lvert \etadiffs\rvert \equiv 1$. Although it is relatively easy to approximate
functions by sequences in the norms we need, the typical discrete approximation will not satisfy the condition $\lvert \fordiff\eta\rvert\equiv 1$,
which means it does not actually represent a chain. We deal with this by using the spherical representation
$\etadiffs(s) = \big( \cos{\theta(s)}, \sin{\theta(s)}\big)$ as in \eqref{spherical}. (Although this formula works only when $d=2$, we can use a similar procedure in higher
dimensions, using generalized spherical coordinates.)
Using $\eta(1) = 0$, we can easily reconstruct 
$\eta$ if $\theta$ is known. We can then approximate the function $\theta$ by a sequence $\theta_k$ and rebuild $\eta_k$ 
using the formula $\eta_{n+1}=0$ and $\fordiff\eta_k = (\cos{\theta_k}, \sin{\theta_k})$ when $d=2$, with a similar formula in higher
dimensions.

Fortunately, the Sobolev norms of $\eta$ and $\theta$ are closely related. 

\begin{proposition}\label{etatothetasobolev}
If $\eta\colon [0,1]\to\mathbb{R}^2$ is related to $\theta\colon [0,1]\to \mathbb{R}$ by the formula \eqref{spherical}, 
with $\eta(1)=\eta''(1)$ and $\theta'(1)=0$,
then boundedness of the squared norm
\begin{equation}\label{thetasobolev} 
A = \int_0^1 \Big( s^2 \lvert \theta'(s)\rvert^2 + s^3 \lvert \theta''(s)\rvert^2 + s^4 \lvert \theta'''(s)\rvert^2\Big) \, ds
\end{equation}
is equivalent to boundedness of the squared norm
\begin{equation}\label{etathirdsobolev}
B = \int_0^1 \Big( s^2 \lvert \eta''(s)\rvert^2 + s^3 \lvert \eta'''(x)\rvert^2 + s^4 \lvert \eta^{(4)}(x)\rvert^2 \Big) \, dx.
\end{equation}
\end{proposition}

\begin{proof}
We easily compute that \begin{equation}\label{etatothetaderiv}
\begin{split}
&\lvert \eta''(s)\rvert^2 = \theta'(s)^2, \quad \lvert \eta'''(s)\rvert^2 = \theta''(s)^2 + \theta'(s)^4, \\
&\text{and}\quad \lvert \eta^{(4)}(s)\rvert^2 = \big(\theta'''(s)-\theta'(s)^3\big)^2 + 9\theta'(s)^2\theta''(s)^2.
\end{split}
\end{equation}
An integration by parts using $\theta'(1)=0$ shows that \eqref{thetasobolev} and \eqref{etathirdsobolev} are related by
\begin{equation}\label{BtoA}
B = A + \int_0^1 \Big(  s^3 \theta'(s)^4 - 8s^3 \theta'(s)^4 + 15 s^4 \theta'(s)^2 \theta''(s)^2 + s^4 \theta'(s)^6\Big) \, ds.
\end{equation}
Repeated use of the basic weighted Sobolev inequalities of Theorem \ref{basicinequalitiesthm} allows us to express every term 
on the right side in terms of $A$, so we get an inequality of the form $B\lesssim A+A^2+A^3$. 

In the other direction, \eqref{BtoA} gives 
$$ A \le 
B + 8\int_0^1 s^3 \lvert \eta''(s)\rvert^4 \, ds,$$ 
and again using Theorem \ref{basicinequalitiesthm} gives $A\lesssim B+B^2$.
\end{proof}

We can derive the same sort of result for difference quotients. If we write $\fordiff\eta_k = (\cos{\theta_k}, \sin{\theta_k})$ for 
$1\le k\le n$, then the analogues of \eqref{etatothetaderiv} are as follows:
\begin{align*}
\lvert \fordiff^2\eta\rvert^2 &= 4n^2 \sin^2{\left(\tfrac{\fordiff\theta}{2n}\right)}, \\
\lvert \fordiff^3\eta\rvert^2 &= 16n^4 \sin^2{\left(\tfrac{\fordiff\theta}{2n} + \tfrac{\fordiff^2\theta}{2n^2}\right)} \sin^2{\left(\tfrac{\fordiff\theta}{2n}\right)}
+ 4n^4 \sin^2{\left( \tfrac{\fordiff^2\theta}{2n^2}\right)} \\
\lvert \fordiff^4\eta\rvert^2 &= 4n^6 \left\lvert \sin{\left( \tfrac{\fordiff^3\theta + n\fordiff^2\theta + n^2\fordiff\theta}{2n^3}\right)} 
- 3 \sin{\left( \tfrac{\fordiff^2\theta+n\fordiff\theta}{2n^2}\right)} \right\rvert^2 \\
&\qquad\qquad + 48 \sin^2{\left( \tfrac{\fordiff^3\theta + 2n\fordiff^2\theta}{4n^3}\right)}
\sin{\left( \tfrac{\fordiff^2\theta + n\fordiff\theta}{2n^2}\right)} \sin{\left( \tfrac{\fordiff^3\theta + n\fordiff^2\theta + n^2\fordiff\theta}{2n^3}\right)}.
\end{align*}
For sufficiently large $n$, we can proceed as in Proposition \ref{etatothetasobolev} to show that the discrete squared norms 
$$ \frac{1}{n} \sum_{k=1}^{n-1} s_k^{(2)} \lvert \fordiff^2\eta_k\rvert^2 + \frac{1}{n} \sum_{k=1}^{n-2} s_k^{(3)} \lvert \fordiff^3\eta_k\rvert^2 
+ \frac{1}{n} \sum_{k=1}^{n-3} s_k^{(4)} \lvert \fordiff^4\eta_k\rvert^2 $$
and 
$$ \frac{1}{n} \sum_{k=1}^{n-1} s_k^{(2)} \lvert \fordiff\theta_k\rvert^2 + \frac{1}{n} \sum_{k=1}^{n-2} s_k^{(3)} \lvert \fordiff^2\theta_k\rvert^2 
+ \frac{1}{n} \sum_{k=1}^{n-3} s_k^{(4)} \lvert \fordiff^3\theta_k\rvert^2$$
can each be bounded in terms of the other.

Thus for either whips or chains in two dimensions, we can work directly in terms of Sobolev norms of $\theta$. The most convenient way 
to map from Sobolev spaces of continuous maps to Sobolev spaces of discrete sequences is to use orthogonal polynomials. (A direct
approach, using values of the function on a discrete grid, does not work for our purposes since bounds on the differences require more
smoothness of the function than we have.)

The only complication is the oddness requirement on $\eta$ (and the discrete oddness criterion \eqref{discreteodd}). In terms of the spherical variable $\theta$, oddness of $\eta$ through $s=1$ translates into evenness of $\theta$, i.e., there is an extension of $\theta$ to $[0,2]$ such that $\theta(2-s)=\theta(s)$. Similarly the discrete oddness condition \eqref{discreteodd} translates into the discrete evenness condition 
$\theta_{2n+1-k}=\theta_k$. These conditions are easy to handle if we extend the interval to $[0,2]$ (or extend the sequence to $\{1, 2, \ldots, 2n\}$) and use Sobolev seminorms with symmetric weights 
\begin{equation}\label{symmetrizedsobolevsmooth}
\llangle \theta, \theta\rrangle_{\rho, j} = \int_0^1 \rho(s)^{j+1} \lvert \theta^{(j)}(s)\rvert^2 \, ds \quad\text{where}\quad \rho(s) = s(2-s)
\end{equation}
and 
\begin{equation}\label{symmetrizedsobolevdiscrete}
\llangle \theta, \theta\rrangle_{\rho, j} = \frac{1}{n} \sum_{k=1}^{n-\lfloor j/2\rfloor} \rho_k^{(j+1)} \lvert \fordiff^j\theta_k\rvert^2 \quad \text{where}\quad 
\rho_k = \frac{k(2n+1-k)}{n^2}.
\end{equation}
These norms are clearly topologically equivalent to the weighted norms we have been using, where the weights are $s$ and $\frac{k}{n}$. 

\begin{theorem}\label{orthogonalexpansion}
There are polynomials $Q_m(s)$ on $[0,2]$ satisfying $Q_m(2-s)=Q_m(s)$,
and such that
\begin{equation}\label{capitalorthogonal}
\llangle Q_{\ell}, Q_m\rrangle_{\rho, j} = \delta_{\ell m} r_{mj}, \quad \text{where} \quad r_{mj} = \frac{(2m+j)!}{(2m-j-2)!(2m)(2m-1)}.
\end{equation}
in the weighted Sobolev seminorms \eqref{symmetrizedsobolevsmooth} for all $j\ge 0$. Thus we can expand $\theta(s) = \sum_{m=1}^{\infty} A_m Q_m(s)$ and obtain 
\begin{equation}\label{thetasobolevcontnorm}
\llangle \theta, \theta\rrangle_{\rho, j} = \sum_{m=1}^{\infty} r_{mj} A_m^2
\end{equation}
for all $j\ge 0$.

There are also, for each $n\in\mathbb{N}$, discrete polynomials $q_m(\frac{k}{n})$ defined for $1\le k\le 2n$ and $1\le m\le n$, satisfying $q_m(\frac{2n+1-k}{n})=q_m(\frac{k}{n})$ and such that 
\begin{equation}\label{lowercaseorthogonal}
\llangle q_{\ell}, q_m\rrangle_{\rho, j} = \delta_{\ell m} r_{mj},
\end{equation}
where $r_{mj}$ is as in \eqref{capitalorthogonal}.
Hence if $\theta_k = \sum_{m=1}^n a_m q_m(\tfrac{k}{n})$ for $1\le k\le n$, then 
\begin{equation}\label{thetasobolevdiscretenorm}
\llangle \theta, \theta\rrangle_{\rho, j} = \sum_{m=1}^n r_{mj} a_m^2.
\end{equation}
\end{theorem}

\begin{proof}
The desired polynomials come from a slight variation on the classical Legendre polynomials and the Chebyshev polynomials of a discrete variable. (These are special cases of the Jacobi and Hahn polynomials respectively, with parameters $\alpha=\beta=0$.) The desired formulas follow from general properties of continuous and discrete orthogonal polynomials; see Nikiforov et al.~\cite{NSU} for a good reference.

To obtain \eqref{capitalorthogonal}, we set $Q_m(s) = K_m \, P_{2m-1}'(1-s)$, where $P_r(x) = \frac{1}{2^r r!} \frac{d^r}{dx^r} (x^2-1)^r$ is the usual Legendre polynomial given by the Rodrigues formula and $K_m$ is a constant chosen to make $Q_m$ orthonormal when $j=0$. To obtain \eqref{lowercaseorthogonal}, we set 
$$q_m(\tfrac{k}{n}) = \frac{k_{mn}}{n^{2m-2}} \big[h_{2m-1}^{(0,0)}(k, 2n+1)-h_{2m-1}^{(0,0)}(k-1,2n+1)\big],$$ where $$h_r^{(0,0)}(x,N) = \frac{(-1)^r}{r!} (1-E^{-1})^r[(x+1)\cdots (x+r)(N-1-x)\cdots (N-r-x)]$$ is the Hahn polynomial given in terms of a discrete Rodrigues formula, with $E$ denoting the integer shift operator, and again $k_{mn}$ is a constant chosen to give orthonormality when $j=0$. 

Checking all the conditions is routine using the formulas in \cite{NSU}.
\end{proof}

For each $n\in\mathbb{N}$, we can define the map $F_n$ which takes a continuous angular function $\theta(s)$ to a discrete approximation $\theta_k$, and the map $G_n$ which takes a discrete angular sequence $\theta_k$ to a continuous angular function $\theta(s)$, by the formulas
\begin{equation}\label{discretization}
\begin{split}
\theta(s) &= \sum_{m=1}^{\infty} A_m Q_m(s) \mapsto \sum_{m=1}^n A_m q_m(\tfrac{k}{n}) = \theta_k \\
\theta_k &= \sum_{m=1}^n a_m q_m(\tfrac{k}{n}) \mapsto \sum_{m=1}^n a_m Q_m(s) = \theta(s),
\end{split}
\end{equation}
where the coefficients are obtained using orthonormality by $$A_m = \int_0^1 \theta(s) Q_m(s)\,ds \quad \text{and}\quad a_m = \frac{1}{n} \sum_{k=1}^n \theta_k q_m(\tfrac{k}{n}).$$
By the formulas \eqref{thetasobolevcontnorm} and \eqref{thetasobolevdiscretenorm}, we can bound the continuous and discrete Sobolev norms of any order $j$ in terms of each other using this map. Furthermore $G_n$ is an isometry, $F_n\circ G_n$ is the identity, and $G_n\circ F_n$ converges strongly to the identity as $n\to\infty$ in any weighted $(\rho, j)$-norm.

Thus given an initial condition $\eta(0,s) = \gamma(s)$, we can write the discrete initial condition $\gamma_n$ as
$$(\gamma_n)_k = -\frac{1}{n} \sum_{j=k}^n (\fordiff \gamma_n)_j(t) = -\frac{1}{n} \sum_{j=k}^n \big( \cos{\theta_k}, \sin{\theta_k}\big)$$
where $\theta_k$ is the discretization obtained from \eqref{discretization}. And conversely, if we solve the discrete chain equations to obtain $\eta_k(t)$, we can construct an approximate whip solution by finding, for each $t$, the angles $\theta_k(t)$ and using \eqref{discretization} to obtain the function $\theta(t,s)$, then reconstructing $\eta(t,s) = -\int_s^1 \big(\cos{\theta(t,x)}, \sin{\theta(t,x)}\big) \, dx$.

We clearly have a similar construction for the velocity $\etadifft(t,s)$ in terms of the angular velocity $\thetadifft(t,s)$, which works based on the formulas 
\begin{align*}
\etadiffst(t,s) &= (-\sin{\theta(t,s)}, \cos{\theta(t,s)}) \thetadifft(t,s) \\
\fordiff\dot{\eta}_k(t) &= (-\sin{\theta_k(t)}, \cos{\theta_k(t)}) \dot{\theta}_k(t).
\end{align*}

These constructions ensure that we can go back and forth between whips and chains while preserving the Sobolev norms as well as the constraint equation.

\subsection{Uniform energy bounds}

Now suppose that the initial whip conditions $\eta(0,s)=\gamma(s)$ 
and $\etadifft(0,s)=w(s)$ have bounded energy $E_3(0)$ given by \eqref{correctweightenergy}, as well as satisfying the constraints $\lvert \gamma'(s)\rvert^2 \equiv 1$ and $\langle \gamma'(s), w'(s)\rangle \equiv 0$, and have odd extensions through $s=1$. Using the procedure of the preceding section, we know that for each $n\in\mathbb{N}$ there are discrete initial conditions $\gamma_n$ and $w_n$ such that the discrete energy $e_3$ given by \eqref{discreteenergydef} is bounded uniformly, independently of $n$. These approximate conditions converge strongly in $N_4[0,1]$ and $N_3[0,1]$ respectively to the actual initial conditions.

\begin{lemma}\label{uniformestimateslemma}
Suppose $\gamma$ and $w$ are initial conditions as in Theorem \ref{theorem1}, 
and suppose discretizations $\gamma_n$ and $w_n$ 
are defined as in Section \ref{interpolation}.
Let $(\eta_n)_k(t)$ and $(\sigma_n)_k(t)$, for $1\le k\le n$, be the
solution of equations \eqref{fancydiscreteevolution} and \eqref{fancydiscretetension} with
$\eta_n(0)=\gamma_n$ and $\dot{\eta}_n(0)=w_n$.

Then there is a $T>0$ such that the discrete energy
$ e_3(t)$ defined by \eqref{discreteenergydef}
is bounded uniformly on $[0,T]$ and uniformly in $n$.
\end{lemma}

\begin{proof}
Since the discrete energy $e_3(0)$ is bounded
uniformly for all $n$, we conclude by Lemma \ref{normtoenergy} that $\supnorm{\eta_n(0)}_{3/2,2}$ is uniformly bounded for all $n$. In particular we know the hypotheses of Lemma \ref{bboundslemma} are satisfied for $n$ sufficiently large.
Fix such an $n$.

By Lemma \ref{bboundslemma} we have 
\begin{equation}\label{bestimate}
\frac{s_k}{\sigma_k(t)} \le b(t) \lesssim \frac{1}{u_0}e^{2\supnorm{\eta(t)}^2_{3/2, 2}}.
\end{equation}
We want an estimate for the evolution of $\supnorm{\eta}^2_{3/2,2}$. For any $k\in \{1,\ldots, n-1\}$, we have by 
Proposition \ref{weightsprop} and Lemma \ref{normtoenergy} that
\begin{align*} \frac{d}{dt}\Big( s_k^{(3/2)} \lvert \fordiff^2\eta_k(t)\rvert^2\Big) &= 2s_k^{(3/2)} \langle \fordiff^2\eta_k(t), \fordiff^2\dot{\eta}_k(t)\rangle \\
&\le \supnorm{\eta(t)}_{1,2} \supnorm{\dot{\eta}(t)}_{2,2} \\
&\lesssim e_3(t).
\end{align*}
Since this is true for any $k$, we conclude 
\begin{equation}\label{threehalvestime}
\supnorm{\eta(t)}^2_{3/2, 2} \le \supnorm{\eta(0)}^2_{3/2, 2} + \frac{L}{4} \int_0^t e_3(\tau)\,d\tau
\end{equation}
for some constant $L$, independent of $t$ and $n$. This bound also ensures that the hypotheses of 
Lemma \ref{bboundslemma} are satisfied for sufficiently large $n$ as long as $e_3(t)$ is bounded.  

By the definitions \eqref{discreteenergydef} and \eqref{discretetimeenergydef}, we clearly have 
$$ e_3(t) \le \max{\{1,b(t)\}}^4 \tilde{e}_3(t),$$
and we conclude by combining \eqref{bestimate} and \eqref{threehalvestime} that 
\begin{equation}\label{tilderelation}
e_3(t) \le K \exp{\left( L \int_0^t e_3(\tau)\,d\tau\right)} \tilde{e}_3(t),
\end{equation}
for some constant $K$ which is also independent of $t$ and $n$.

Let $y(t) = \int_0^t e_3(\tau)\,d\tau$ and let $z(t) = \int_0^t \tilde{e}_3(\tau)\,d\tau$. Then \eqref{tilderelation} 
can be written as 
$$ e^{-Ly(t)} \frac{dy}{dt} \le K\frac{dz}{dt},$$
and integrating both sides yields 
\begin{equation}\label{yz}
e^{-Ly(t)} \ge 1 - KLz(t).
\end{equation}

Now we use Theorem \ref{discreteenergyestimate} to get $\frac{d\tilde{e}_3}{dt} \le M_3 e_3(t)^7$ for some 
$M_3$ independent of $t$ and $n$. Using \eqref{tilderelation} and \eqref{yz}, 
we have $$e_3^7 \le K^7 e^{7Ly(t)} \tilde{e}_3^7 \le \frac{K^7}{(1-KLz(t))^7} \, \left(\frac{dz}{dt}\right)^7,$$
so we obtain 
\begin{equation}\label{zprimeprime}
\frac{d^2z}{dt^2} \le \frac{M_3K^7}{(1-KLz(t))^7} \, \left( \frac{dz}{dt}\right)^7.
\end{equation}
Dividing by $(\frac{dz}{dt})^6$ and integrating, we obtain 
\begin{equation}\label{zprime}
z'(t) \le z'(0) \left[ 1 + J\left(1-\frac{1}{(1-KLz(t))^6}\right)\right]^{-1/5},
\end{equation}
where $J=\frac{5M_3K^6z'(0)^5}{6L}$. Another integration gives a bound for 
$z(t)$ on some time interval $[0,T]$, which depends only on $e_3(0)$, and 
\eqref{zprime} gives a uniform bound on $\tilde{e}_3(t)$. Combining this with \eqref{tilderelation} 
and \eqref{yz}, we get a uniform bound on $e_3(t)$ as well on the same time interval.
\end{proof}

Now having obtained a sequence of chain solutions $\eta_n(t)$, bounded uniformly
in the discrete weighted Sobolev norms uniformly on an interval $[0,T]$, we use the 
technique of Section \ref{interpolation} to interpolate. For each $n$ we obtain an
approximate whip solution $\overline{\eta}_n \colon [0,T]\times [0,1]\to \mathbb{R}^d$ for which the energy $E_3(t)$ is bounded on $[0,T]$ independently of $n$. We can then extract a subsequence which 
converges in the weak-* topology on $L^{\infty}([0,T], N_4[0,1])$.

Before doing this, we prove one final lemma, a compactness result analogous to the usual 
Rellich theorem.

\begin{lemma}\label{compactnesslemma}
Let $\overline{N}_m[0,2]$ denote the space of functions $\eta\colon [0,2]\to\mathbb{R}^d$ such that 
the norm
\begin{equation}\label{overlinenorm}
\lVert \eta\rVert^2_{\overline{N}_m} = \sum_{\ell=0}^m \int_0^2 s^{\ell}(2-s)^{\ell} \left\lvert \frac{d^{\ell}\eta}{ds^{\ell}}\right\rvert^2 \, ds
\end{equation}
is finite.

Then $\overline{N}_{m+1}[0,2]$ is compact in $\overline{N}_m[0,2]$ for each $m\ge 0$.
\end{lemma}

\begin{proof}
Expand $\eta(s) = \sum_{j=0}^{\infty} w_j P_j(1-s)$, where $P_j$ are the standard Legendre polynomials. Then as discussed in Section \ref{interpolation}, we have 
$$ \lVert \eta\rVert^2_{\overline{N}_m} = \sum_{\ell=0}^m \sum_{j=\ell}^{\infty} \frac{2}{2j+1} \frac{(j+\ell)!}{(j-\ell)!} w_j^2.$$
Hence the embedding $\iota\colon \overline{N}_{m+1}\to \overline{N}_m$ is a norm limit of operators with finite-dimensional range, so it is compact.
\end{proof}

As noted in Section \ref{interpolation}, for functions on $[0,1]$ that are restrictions of odd functions on $[0,2]$, the norm on $\overline{N}_m[0,2]$ given by \eqref{overlinenorm} is equivalent to the norm on $N_m[0,1]$ given by \eqref{Nnorm}, and thus we get compactness of $N_{m+1}[0,1]$ in $N_m[0,1]$ for functions with an odd extension through $s=1$.

We now establish the existence part of Theorem \ref{theorem1}.

\begin{theorem}\label{convergencetheorem}
Given initial conditions $\gamma$ and $w$ as in Theorem \ref{theorem1}, 
there is a $T>0$ such that there is a solution $\eta$ of
the system \eqref{fullsystem} in $L^{\infty}([0,T], N_4[0,1])\cap
W^{1,\infty}([0,T], N_3[0,1])$.
\end{theorem}

\begin{proof}
For each fixed $t$ and each $n\in \mathbb{N}$, construct a continuous approximation
of the chain $\eta_n(t)$ as in Section \ref{interpolation},
and call it $\overline{\eta}_n(t)$.
Then by Lemma \ref{uniformestimateslemma} we get a uniform bound on
on $E_3(t)$ in some short time interval $[0,T]$; in other words, the
family $\overline{\eta}_n$ is bounded in $L^{\infty}([0,T], N_4[0,1])\cap
W^{1,\infty}([0,T], N_3[0,1])$. By the Alaoglu theorem, there is a subsequence $\eta_{n_k}$
that converges in the weak-* topology 
to $\eta \in L^{\infty}([0,T], N_4[0,1])\cap W^{1,\infty}([0,T], N_3[0,1])$. 

By the compactness Lemma \ref{compactnesslemma}, there is a
sub-subsequence $\tilde{\eta}_{n_{k_j}}$ which converges strongly to $\eta$
in $L^{\infty}([0,T], N_3[0,1])\cap W^{1,\infty}([0,T], N_2[0,1])$.
For any $\epsilon>0$ the convergence is strong in $H^3[\epsilon,1]$,
and thus by the usual Sobolev embedding theorem also in
$C^2[\epsilon,1]$. So we can take the limit of the system
\eqref{fancydiscreteevolution} and \eqref{fancydiscretetension} 
pointwise to see that we have a solution
of \eqref{fullsystem}.
\end{proof}

The fact that all the estimates close up at the level of $E_3$, with all other energies satisfying linear differential inequalities, implies that the only way a solution which is initially $C^{\infty}$ can fail to be $C^{\infty}$ for all time is if $E_3$ becomes infinite in finite time. This gives a crude blowup criterion.

\begin{corollary}\label{blowupcriterion}
Suppose $\eta$, $\sigma$ is a solution of the system
\begin{multline}\label{fullsystem}
\etadifftt =
\partial_s(\sigma \etadiffs), \qquad  \sigmadiffss - \lvert
\etadiffss\rvert^2 = -\lvert \etadiffst\rvert^2, \qquad \lvert \etadiffs\rvert^2\equiv 1, \\
\eta(t,1)=0, \;\sigmadiffs(1)=0,\; \sigma(0)=0, \quad
\eta(0,s)=\gamma(s),\; \etadifft(0,s) = w(s),
\end{multline}
where we assume that $\gamma$ and $w$ satisfy the conditions of Theorem \ref{theorem1}.

Assume that in some time interval $[0,T]$, the
energy $E_3(t)$ is bounded uniformly. Assume further
that $E_m(0)$ is bounded for all $m>3$. Then $E_m(t)$ 
is also bounded in $[0,T]$ for all $m>3$.
\end{corollary}

\begin{proof}
By equation \eqref{generalenergyderivativesmooth}, we have for $k\ge 4$ that
$$ \frac{d\tilde{E}_k}{dt} \lesssim M_k(E_{k-1}) E_k.$$ Furthermore since $E_2(t)$
is bounded, so is $B(t) = \sup_s \frac{s}{\sigma(t,s)}$ by Lemma \ref{bboundslemma}, 
and thus
$$ \frac{d\tilde{E}_k}{dt} \le \tilde{M}_k(E_{k-1}) \tilde{E}_k$$
for some function $\tilde{M}_k$. 
So by Gronwall's inequality, $\tilde{E}_k(t)$ 
is bounded on $[0,T]$ in terms of $\tilde{E}_k(0)$. Thus finally $E_k(t)$ 
is also bounded on $[0,T]$.
\end{proof}

We now complete the proof of Theorem \ref{theorem1} by proving uniqueness.

\begin{theorem}\label{uniqueness}
Suppose $\gamma$ and $w$ are functions on $[0,1]$ as in Theorem \ref{theorem1}.
If $(\eta_1, \sigma_1)$ and $(\eta_2, \sigma_2)$ are two solutions of \eqref{fullsystem}, both in $$L^{\infty}([0,T], N_4[0,1])\cap
W^{1,\infty}([0,T], N_3[0,1]),$$
with the same initial conditions  $$\eta_1(0,s)=\eta_2(0,s)=\gamma(s) \quad \text{and} \quad \partial_t\eta_1(0,s)=\partial_t\eta_2(0,s) = w(s),$$  then $\eta_1(t,s)=\eta_2(t,s)$ and $\sigma_1(t,s) = \sigma_2(t,s)$ for all $t\in [0,T]$ and all $s\in [0,1]$. 
\end{theorem}

\begin{proof}
The proof relies on an energy estimate for the differences at the level of the first energy,
\begin{multline*} E_1[\eta_1-\eta_2] = \int_0^1 \Big( \lvert \partial_t\eta_1-\partial_t\eta_2\rvert^2
+ s \lvert \partial_s\eta_1-\partial_s\eta_2\rvert^2\\ + s \lvert \partial_t\partial_s\eta_1-\partial_t\partial_s\eta_2\rvert^2 
+ s^2 \lvert \partial_s^2\eta_1-\partial_s^2\eta_2\rvert^2 \Big) \, ds.
\end{multline*}
We estimate this energy using a Gronwall inequality, as in Theorem \ref{discreteenergyestimate} and Lemma \ref{uniformestimateslemma}.
The reason this works is that since $\eta_1-\eta_2$ satisfies a linear PDE whose coefficients involve the known quantities $\eta_1$, $\eta_2$, $\sigma_1$, and $\sigma_2$, we can use Corollary \ref{supproductboundcorollary} in a more effective way to put all the weights on the known terms. 
The full proof appears in Appendix \ref{A5}.
\end{proof}

Finally we discuss some refinements of these results. First, given a solution $\eta$ of \eqref{basicPDE} with 
$E_3$ finite, we can check using the differential equation that $\lVert \partial_t^2\eta\rVert_{2,2}$, $\lVert \partial_t^3\eta\rVert_{1,1}$, and $\lVert \partial_t^4\eta\rVert_{0,0}$ can all be bounded in terms of $E_3$. Hence the solution is also in
$ 
W^{2,\infty}([0,T],N_2[0,1]) \cap W^{3,\infty}([0,T], N_1[0,1]) \cap W^{4,\infty}([0,T], N_0[0,1]).$
Now a well-known general technique (see e.g., Lemma 11.9 of \cite{RR}) shows that $\eta$ is continuous
as a curve in $N_4[0,1]$, $\eta$ is $C^1$ as a curve in $N_3[0,1]$, etc.

\section{General remarks and future research}\label{futureresearchsection}

In this paper we considered the whip with one fixed and one free end
as boundary conditions. The other possibilities are to have two free
ends, to have two fixed ends, and to have periodicity. All of
the estimates in this paper have analogues in those cases. When
there are two free ends, the tension must satisfy $\sigma(0)=0$ and
$\sigma(2)=0$, so the appropriate weighted norms look like the square root of $\int_0^2
s^k (2-s)^k \lvert f^{(k)}(s)\rvert^2 \,ds$. Since we have essentially
solved the problem with one fixed end by constructing an odd extension
in order to turn the problem into a string with two free ends on $[0,2]$,
we expect that the same estimates prove existence for an inextensible string
with two free ends. 
When there are two
fixed ends, or when the whip is periodic, the problem becomes
simpler since we can use ordinary Sobolev spaces for the estimates. In this 
case we expect the energy estimates to close up at the level of $e_2$ rather than 
$e_3$.

The addition of gravity brings some complications. One is that the
boundary conditions change, and oddness through the fixed point is
no longer enough to satisfy the conditions automatically. (This is already an 
issue even for the wave equation with constant coefficients, if an external force
is imposed which does not respect the boundary conditions.) The other complication
is that if the whip is above the fixed point, the
tension may become negative: the effect of gravity is to change the boundary
condition in \eqref{generalstringconst} to $\sigmadiffs(t,1) = \langle g, \etadiffs(t,1)\rangle$,
where $g$ is the gravitational acceleration vector, and if $\sigmadiffs(t,1)<0$ then it is 
possible to have $\sigma(t,s)<0$ for some $t$ and $s$. In that case the evolution equation
becomes elliptic, so the discussion becomes much more complicated.

The blowup criterion Corollary \ref{blowupcriterion}, that a smooth solution remains smooth up to time $T$ 
iff $\sup_{0\le t\le T}E_3(t) < \infty$, can certainly be
improved. Once we know a solution exists, we can use alternative
methods to get better a priori bounds on it. Thess et al. have
speculated that blowup for the periodic loop might be controlled by
the $L^{\infty}$ norms of $\lvert \etadiffss\rvert$ and $\lvert
\etadiffst\rvert$, analogous to the way blowup for the ideal Euler
equations is controlled by the $L^{\infty}$ norm of vorticity. This
is an interesting problem to study, since we have a much greater
handle on all aspects of this one-dimensional problem. We will
explore this in a future paper.

In addition, the geometry of the space of inextensible curves is
interesting in its own right. Although the geometric objects are not
smooth in the Sobolev topology (unlike on the group of
volumorphisms), the curvature formulas still make sense, and one can
compute formally that all sectional curvatures are nonnegative. We can thus
try to study stability of the motion from the geometric point of view (as
in \cite{ak}), as well as the geometry of blowup. See
\cite{whipcurvature} for details on this.

A similar problem in higher dimensions is given by the motion of a flag 
attached to a pole in $3$-space. Here our configuration space would be the 
space of maps of a rectangle into $\mathbb{R}^3$ which are isometric immersions
with one side of the rectangle held fixed. We expect to see a similar nonlocal
coupled degenerate system, the only obvious difference being that the ordinary
differential equation \eqref{tensionODE} becomes an elliptic equation in 
the spatial variables. 

The whip-chain equations are interesting partly in and of
themselves, but especially as a ``toy model'' of inviscid, incompressible fluids. There are
some structural similarities between the equations \eqref{basicPDE}
and \eqref{tensionODE} and the Euler equation for an ideal fluid, 
given in Lagrangian form by
\begin{equation*}
\etadifftt(t,x) = -\grad p\big(t, \eta(t,x)\big)
\quad \text{and}\quad
\Laplacian p =
-\text{Tr}\big([D\etadifft(t,x)\circ\eta^{-1}(t,x)]^2\big),
\end{equation*}
with some boundary condition to determine $\grad p$ uniquely. Both systems 
involve a hyperbolic evolution equation for a constrained function, 
where the right side is given in terms of a function determined by a purely spatial differential equation.
The technique of approximating a continuous system with a
discrete system preserving the geometry may be interesting to apply
to fluids directly. For example, in two dimensions we could consider a rectangular grid on a torus, the
vertices of which are free to move as long as all quadrilateral
areas are preserved. Although such a model may not have global
existence (as edges of a quadrilateral may collapse to give a
triangle without changing the area), we might still get some useful
insight out of it.

\appendix

\section{Longer proofs}\label{appendix}

\subsection{Proof of Proposition \ref{sigmalowerbound}}\label{A1}

\begin{prop}
Suppose $G_{kj}$, $\eta_k$, $\alpha_k$, and $\beta_k$ are defined
as in Proposition \ref{solverprop}. Assume the $\eta_k$ are such
that, for some $\upsilon\in (0, \frac{2\sqrt{n}}{5}]$, we have  
\begin{equation}\label{threehalvesapp}
\frac{k^{3/2}}{n^{3/2}} \lvert \fordiff^2\eta_k\rvert^2 \le \upsilon \quad \text{for $1\le k\le n-1$.}
\end{equation}
Then for every $1\le j,k\le n$, we have
\begin{equation}\label{lowerbounddiscapp}
\frac{n^2G_{kj}}{jk} \ge e^{-2\upsilon}.
\end{equation}

If $G$ solves \eqref{greencontinuous} and $\etadiffss$ is a smooth function, then we have
\begin{equation}\label{lowerboundcontapp}
\inf_{0\le s,x\le 1} \frac{G(s,x)}{sx} \ge \frac{e^{-\varrho}}{1+\varrho}
\text{ where } \varrho = \int_0^1 s\lvert \etadiffss\rvert^2 \, ds = \lVert \eta\rVert_{1,2}.
\end{equation}
\end{prop}

\begin{proof}
Our strategy for proving \eqref{lowerboundcontapp} 
will be to first show that the minimum of the ratios  $\frac{n^2G_{kj}}{kj}$ and $\frac{G(s,x)}{sx}$ is attained at the off-diagonal corners; that is, $\min_{1\le j,k\le n} \frac{n^2G_{kj}}{kj} = n G_{1n}$ and $\inf_{0< s,x,\le 1} \frac{G(s,x)}{sx} = \lim_{s\to 0} \frac{G(s,1)}{s} = \Gdiffs(0,1)$. The proofs are nearly identical in both cases, so we will just give the discrete proof. Then we estimate the size of this value; here the proofs are different, and we can get a sharper estimate for the continuous case.

We first define
a matrix $F$ by
$F_{kj} = \frac{n^2}{kj} G_{kj}$. Clearly $F$ is symmetric since $G$ is.
We want to prove that $F_{kj} \ge F_{1n}$.
Note that for $1<k\le n$ we
have \begin{equation}\label{Fdifference}
 F_{kj}-F_{k-1,j} 
=\frac{n^2}{jk(k-1)} \left[ (k-1)(G_{kj}-G_{k-1,j}) - G_{k-1,j}
\right].
\end{equation}

First we show that we can decrease $F_{kj}$ by increasing the larger index.
If $k>j$ we know from Proposition \ref{greenupperbounds} that 
$(\backdiffone G)_{kj} \le 0$, so that $G_{kj}\le G_{k-1,j}$, and thus
$F_{kj}-F_{k-1,j}\le 0$. Thus 
\begin{equation}\label{increaselarger}
F_{kj}\ge F_{nj} \text{ if $k\ge j$}.
\end{equation}

Next we show that we can decrease $F_{kj}$ by decreasing the smaller index, which
is a bit more involved.
Inspired by \eqref{Fdifference}, we define for $1\le k\le n$ the auxiliary quantity
$H_{kj} = (k-1)(G_{kj}-G_{k-1,j})-G_{k-1,j}$; then it is easy to
compute that $(\fordiffone H)_{kj} = k(\backdiffone\fordiffone G)_{kj}$,
and we conclude using \eqref{Gdiscreteconvex} that if $k<j$ then
$H_{k+1,j}-H_{kj}\ge 0$. Since $H_{1j} = 0$, this shows that
$H_{kj}\ge 0$ as long as $k\le j$. Then since $F_{kj}-F_{k-1,j} =
\frac{n^2}{jk(k-1)} H_{kj}$ for $k>1$, we have $F_{kj} \ge F_{k-1,j}$
as long as $1<k\le j$, and hence 
\begin{equation}\label{decreasesmaller}
F_{kj}\ge F_{1j} \text{ for $1\le k\le j$}.
\end{equation} 

Combining \eqref{increaselarger} and \eqref{decreasesmaller}, and using the 
fact that $F_{kj}=F_{jk}$, we obtain
\begin{equation}\label{Flowerbound}
\min_{1\le j,k\le n} F_{kj} = F_{1n}.
\end{equation}
We finally want to bound $F_{1n}$ from below. Using the formula \eqref{discretegreen} we have that
\begin{equation}\label{discreteproductbound}
F_{1n} = \frac{p_{11}p_{1n}}{\beta_1} = \frac{1}{\beta_1} \prod_{m=1}^{n-1}
\frac{\alpha_m}{\beta_{m+1}}
\end{equation}
It is easier to estimate sums than products, so we rewrite \eqref{discreteproductbound} as \begin{equation}\label{Flog}
\ln{F_{1n}} = \sum_{k=1}^{n-1} \ln{\alpha_k} - \sum_{k=1}^{n-1} \ln{\beta_k},
\end{equation}
recalling that $\beta_n=1$.

First we get an upper estimate for $\sum_{k=1}^{n-1} \ln{\beta_k}$. 
Rearranging \eqref{betarecursion} and using $\alpha_k = 1-\frac{1}{2n^2} \lvert \fordiff^2\eta_k\rvert^2$, we have 
$$ \beta_k -\beta_{k+1} = -(\beta_k-1)(\beta_{k+1}-1) + \frac{\lvert \fordiff^2\eta_k\rvert^2}{n^2} - \frac{\lvert \fordiff^2\eta_k\rvert^4}{4n^4}.$$
Recalling that $1\le \beta_k\le 2$ for each $k$, we conclude
$$ \beta_k - \beta_{k+1} \le \frac{1}{n^2} \lvert \fordiff^2\eta_k\rvert^2,$$
and since $\beta_n=1$, we find
\begin{equation}\label{betastep}
\beta_j = \beta_n + \sum_{k=j}^{n-1} (\beta_k - \beta_{k+1}) \le 1 + \frac{1}{n^2} \sum_{k=j}^{n-1} \lvert \fordiff^2\eta_k\rvert^2.
\end{equation}

Now since $\beta_j\ge 1$, we have $\ln{\beta_j} \le \beta_j-1$, so that (incorporating the assumption \eqref{threehalvesapp})
\begin{equation}\label{betauppersum}
\begin{split} 
\sum_{j=1}^{n-1} \ln{\beta_k} &\le \frac{1}{n^2} \sum_{j=1}^{n-1} \sum_{k=j}^{n-1} \lvert \fordiff^2\eta_k\rvert^2
= \frac{1}{n^2} \sum_{k=1}^{n-1} k\lvert \fordiff^2\eta_k\rvert^2 \\
&\le \frac{\upsilon}{n} \sum_{k=1}^{n-1} \sqrt{\frac{n}{k}} \le \upsilon\left(2 + \frac{\zeta(1/2)}{\sqrt{n}}\right),
\end{split}
\end{equation}
where $\zeta(1/2) \approx -1.46$ is a value of the Riemann zeta 
function.\footnote{Having a precise estimate of this remainder is useful to make part of
\eqref{betauppersum} cancel out \eqref{alphalowersum}, in order to make the estimate 
\eqref{lowerbounddiscapp} independent of $n$ and thus a bit more elegant.}

Next we get a lower estimate for $\sum_{k=1}^{n-1} \ln{\alpha_k}$. Since $\alpha_k = 1-\frac{\lvert \fordiff^2\eta_k\rvert^2}{2n^2}$, the 
assumption \eqref{threehalvesapp} yields $\alpha_k \ge 1 - \frac{\upsilon}{2n^{1/2}k^{3/2}}$.
Since we have assumed $\upsilon \le \frac{2\sqrt{n}}{5}$, we 
have $\alpha_k\ge 0.8>0$ for all $1\le k\le n-1$.
Now we want to get a lower bound for $\ln{\alpha_k}$; this is a bit more delicate than an upper bound for the logarithm of $\ln{\beta_k}$.
Define $c=-2\zeta(\frac{1}{2})/\zeta(\frac{3}{2}) \approx 1.118$ in terms of the Riemann zeta function. 
It is not difficult to verify that when $0.8\le \alpha_k\le 1$, then $\ln{\alpha_k} \ge -c(1-\alpha_k)$. Thus we have 
$\ln{\alpha_k} \ge -(c\upsilon)/(2n^{1/2}k^{3/2})$ for every $k$, from which we conclude
\begin{equation}\label{alphalowersum}
\sum_{k=1}^{n-1} \ln{\alpha_k} 
\ge -\frac{c\upsilon}{2n^{1/2}} \sum_{k=1}^{n-1} \frac{1}{k^{3/2}} 
\ge -\frac{c\upsilon\zeta(\frac{3}{2})}{2n^{1/2}} 
= \frac{\upsilon\zeta(\frac{1}{2})}{\sqrt{n}}.
\end{equation}

Combining \eqref{betauppersum} with \eqref{alphalowersum} and plugging into \eqref{Flog}, we obtain $\ln{F_{1n}} \ge -2\upsilon$. Using 
\eqref{Flowerbound}, we obtain \eqref{lowerbounddiscapp} as desired.

Now we will just sketch the proof of \eqref{lowerboundcontapp}. We similarly establish that the infimum of $\frac{G(s,x)}{sx}$ is attained when $s=0$ and $x=1$, which works the same way as in the discrete case. So we just need to estimate $\lim_{s\to 0} \frac{G(s,1)}{s} = \Gdiffs(0,1)$. 
Letting $J(s)=G(s,1)$, we see that $J$ satisfies
\begin{equation}\label{varphieq}
J''(s) - \lvert \eta''(s)\rvert^2 J(s) = 0, \qquad J(0)=0, \quad J'(1)=1.
\end{equation}
The minimum is then $\Gdiffs(0,1) = J'(0)$. 

Set $\lambda(s) = \ln{[J(s)/s]}$; then we can verify by explicit computation that \eqref{varphieq} can be rewritten in two ways:
\begin{align*} 
\frac{d}{ds} \big[ s^2 \lambda'(s)\big] &= -s^2 \lambda'(s)^2 + s^2 \lvert \eta''(s)\rvert^2, \text{ and}\\
\frac{d}{ds} \big[ s(1-s)\lambda'(s)\big] + \lambda'(s) &= -s(1-s) \lambda'(s)^2 + s(1-s) \lvert \eta''(s)\rvert^2.
\end{align*}
Integrating the first equation from $s=0$ to $s=1$ gives 
\begin{equation}\label{lowsigma1}
\lambda'(1) \le \int_0^1 s^2 \lvert \eta''(s)\rvert^2 \, ds \le \varrho,
\end{equation}
and integrating the second from $s=0$ to $s=1$ gives 
\begin{equation}\label{lowsigma2}
\lambda(1)-\lambda(0) \le \int_0^1 s(1-s) \lvert \eta''(s)\rvert^2 \, ds \le \varrho.
\end{equation}
Since $\lambda'(1) = \frac{1}{J(1)} - 1$ and $\lambda(1)-\lambda(0) = \ln{[J(1)/J(0)]}$, 
estimates \eqref{lowsigma1} and \eqref{lowsigma2} combine to give \eqref{lowerboundcontapp}.
\end{proof}

\subsection{Proof of Theorem \ref{basicinequalitiesthm}}\label{A2}

\begin{theo}
Let $f\colon [0,1]\to \mathbb{R}^d$ be $C^{\infty}$. Then for any
$r>0$ the norms \eqref{weightedsobolevdef} and \eqref{weightedsupremumdef} 
satisfy the weighted inequalities
\begin{align}
\lVert f\rVert^2_{r-1, m} &\lesssim \lVert f\rVert^2_{r, m} + \lVert f\rVert^2_{r+1, m+1} \text{ and} \label{discretepoincareapp} \\
\supnorm{f}^2_{r, m} &\lesssim \lVert f\rVert^2_{r, m} + \lVert f\rVert^2_{r+1, m+1}. \label{discretesobolevapp}
\end{align}
If in addition we have $f^{(m)}(1)=0$, then these inequalities can be simplified to 
\begin{align}
\lVert f\rVert^2_{r-1,m} &\lesssim \lVert f\rVert^2_{r+1,m+1} \text{ and}\label{tieddownpoincareapp} \\
\supnorm{f}^2_{r,m} &\lesssim \lVert f\rVert^2_{r+1,m+1}. \label{tieddownsobolevapp}
\end{align}

If $f$ is instead a sequence $\{f_1,\cdots, f_n\}$ with values in $\mathbb{R}^d$, then the 
inequalities \eqref{discretepoincareapp}--\eqref{discretesobolevapp} also hold if the norms are 
interpreted as \eqref{weightedsobolevdiscretedef} and \eqref{weightedsupremumdiscretedef},
while the inequalities \eqref{tieddownpoincareapp}--\eqref{tieddownsobolevapp} hold if $f_{n-m}^{(m)}=0$. 
\end{theo}

\begin{proof}
It is clearly sufficient to prove these inequalities when $m=0$. 
To derive the discrete versions, we do the following.
Let $p$ be any real number. Then for any $k\in\{1,\ldots, n-1\}$, we have
\begin{equation*}
\lvert (p+k)f_{k+1}-kf_k\rvert^2 = p^2\lvert f_{k+1}\rvert^2 +
\tfrac{k(k+p)}{n^2} \lvert \fordiff f_k\rvert^2 + pk\big(\lvert f_{k+1}\rvert^2 -
\lvert f_k\rvert^2\big),\end{equation*} so that
\begin{equation}\label{lovestinks}
0 \le p^2\lvert f_{k+1}\rvert^2 +
\tfrac{k(k+p)}{n^2} \lvert \fordiff f_k\rvert^2 + pk\big(\lvert f_{k+1}\rvert^2 -
\lvert f_k\rvert^2\big).
\end{equation}

Note that if we define $f_0$ in any way at all, the equation is still satisfied at $k=0$. 
Furthermore we have the easy-to-verify formulas $ks_{k+1}^{(q-1)} = ns_k^{(q)}$
and $\fordiff s_k^{(q)} = qs_{k+1}^{(q-1)}$, which are valid for $k\ge 0$. 

For any real $q>0$, multiply \eqref{lovestinks} through by $s_{k+1}^{(q-1)}$ and simplify to get
$$
0 \le p^2 s_{k+1}^{(q-1)} \lvert f_{k+1}\rvert^2 
+ \tfrac{k+p}{n} s_k^{(q)} \lvert \fordiff f_k\rvert^2
 + np s_k^{(q)} \big(\lvert f_{k+1}\rvert^2 -
\lvert f_k\rvert^2\big).$$
Now notice that the last term simplifies to 
\begin{align*}
np s_k^{(q)} \big( \lvert f_{k+1}\rvert^2 - \lvert f_k\rvert^2\big) &= 
p\Big[\fordiff \big( s_k^{(q)} \lvert f_k\rvert^2\big) + n(s_k^{(q)} - s_{k+1}^{(q)}) \lvert f_{k+1}\rvert^2\Big] \\
&= p\fordiff \big( s_k^{(q)} \lvert f_k\rvert^2 \big) - pq s_{k+1}^{(q-1)} \lvert f_{k+1}\rvert^2,
\end{align*}
using $\fordiff s_k^{(q)} = q s_{k+1}^{(q-1)}$. Thus we have 
$$ 0\le p(p-q) 
s_{k+1}^{(q-1)} \lvert f_{k+1}\rvert^2
+ \frac{k+p}{n} s_k^{(q)} \lvert \fordiff f_k\rvert^2 + 
p \fordiff\big( s_k^{(q)} \lvert f_k\rvert^2\big).$$

Now let $i$ and $j$ be any integers with $0\le i<j\le n$. Summing
all the terms from $k=i$ to $k=j-1$ and using the telescope formula
$\frac{1}{n} \sum_{k=i}^{j-1} \fordiff b_k = b_j-b_i$ for any sequence $\{b_k\}$, we obtain
\begin{multline}\label{discreteproductrule}
0\le \frac{p(p-q)}{n} \sum_{k=i+1}^j 
s_k^{(q-1)} \lvert f_k\rvert^2
+ \frac{1}{n} \sum_{k=i}^{j-1} \frac{k+p}{n} s_k^{(q)} \lvert \fordiff f_k\rvert^2 \\
+ 
p\big( s_j^{(q)} \lvert f_j\rvert^2 - s_i^{(q)} \lvert f_i\rvert^2\big),
\end{multline}
after reindexing the first sum on the right side. This is the basic building block for
all the other inequalities in this proof. Now we consider some special cases which will
together prove \eqref{discretepoincareapp}--\eqref{discretesobolevapp}. Take any $r>0$.

\begin{itemize}
\item 
For any integer $i$ with $1\le i\le n$, 
if we set $p=q=r$ and $j=n$ in \eqref{discreteproductrule}, we get 
$$ s_i^{(r)} \lvert f_i\rvert^2 \le \frac{1}{rn} \sum_{k=i}^{n-1} \frac{k+r}{n} s_k^{(r)} \lvert \fordiff f_k\rvert^2 \\
+ s_n^{(r)} \lvert f_n\rvert^2,$$
and since $\frac{k+r}{n} s_k^{(r)} = s_k^{(r+1)}$, we have 
\begin{equation}\label{basiceq1}
s_i^{(r)} \lvert f_m\rvert^2 \le s_n^{(r)} \lvert f_n\rvert^2 + \frac{1}{r} \lVert f\rVert_{r+1, 1}^2.
\end{equation}
We use this to obtain \eqref{discretesobolevapp}; if $f_n=0$ we obtain \eqref{tieddownsobolevapp}.

\item Next, if we set $p=\frac{r}{2}$, $q=r$, $i=0$, and $j=n$, then we get
$$
0\le -\frac{r^2}{4n} \sum_{k=2}^n s_k^{(r-1)} \lvert f_k\rvert^2 + \frac{1}{n} \sum_{k=1}^{n-1} \frac{k+r/2}{n} s_k^{(r)} \lvert \fordiff f_k\rvert^2+ \frac{r}{2} s_n^{(r)} \lvert f_n\rvert^2.
$$
Noting that 
$k+\frac{r}{2} < k+r$ for any $r>0$ and any $k$, we obtain after solving for $\lVert f\rVert^2_{r-1,0}$ that
\begin{equation}\label{basiceq2}
\lVert f\rVert^2_{r-1,0} \le \frac{4}{r^2} \lVert f\rVert^2_{r+1, 1} + \frac{2}{r} s_n^{(r)} \lvert f_n\rvert^2,
\end{equation}
which is used to bound \eqref{discretepoincareapp}. If $f_n=0$ we obtain \eqref{tieddownpoincareapp}.

\item Finally we get an upper bound for $\lvert f_n\rvert$. Choose $q=r+1$ and $p=-(r+1)$ 
with $i=0$ and $j=n$. Then we have
\begin{align*}
0&\le \frac{2(r+1)^2}{n} \sum_{k=1}^n 
s_k^{(r)} \lvert f_k\rvert^2
+ \frac{1}{n} \sum_{k=0}^{n-1} \frac{k-r-1}{n} s_k^{(r+1)} \lvert \fordiff f_k\rvert^2 \\
&\qquad\qquad\qquad- (r+1) s_n^{(r+1)} \lvert f_n\rvert^2.
\end{align*}
from which we conclude
\begin{equation}\label{basiceq3}
s_n^{(r+1)} \lvert f_n\rvert^2 \le  2(r+1) \lVert f\rVert^2_{r,0} + \frac{1}{r+1}\lVert f\rVert^2_{r+2,1}.
\end{equation}
\end{itemize}

Now we obviously have 
$$ \lVert f\rVert^2_{r+2,1} 
\le \frac{n+r}{n} \lvert f\rVert^2_{r+1,1},$$
and plugging into \eqref{basiceq3} gives
\begin{equation}\label{boundaryestimate}
s_n^{(r)} \lvert f_n\rvert^2 \le \frac{2r^2+4r+1}{r(r+1)} \lVert f\rVert^2_{r+1,1} + 4(r+1) \lVert f\rVert^2_{r,0},
\end{equation}

Combining \eqref{basiceq2} with \eqref{boundaryestimate}, we get
\eqref{discretepoincareapp}. Combining \eqref{basiceq1} with
\eqref{boundaryestimate} and taking the maximum, we get \eqref{discretesobolevapp}.
\end{proof}

\subsection{Proof of Lemma \ref{sigmasobolevlemma}}\label{A3}

\begin{lemm}
If $\sigma$ satisfies \eqref{usefultension} with the condition $\sigma_0=0$, then the norms \eqref{sigmasobolevnorm} can be bounded by the energy \eqref{discreteenergydef} via
\begin{equation}\label{lowdapp}
d_1\lesssim e_3^4, \qquad d_2\lesssim e_3^4, \qquad \text{and} \quad d_3\lesssim e_3^6,
\end{equation} 
while for $m>3$ we have 
\begin{equation}\label{sigmaboundapp}
d_m \le P_m(e_{m-1}) e_m,
\end{equation}
where $P_m$ depends only on $e_{m-1}$. 
\end{lemm}

\begin{proof}
Applying the shift operator $E$ to \eqref{usefultension}, we obtain
$$ \fordiff^2 \sigma = \frac{E^2\sigma}{2} \lvert E\fordiff^2\eta\rvert^2 + \frac{\sigma}{2} \lvert \fordiff^2\eta \rvert^2 - \lvert E\fordiff\dot{\eta}\rvert^2.$$
Thus as a first step, we have 
\begin{align*}
\lVert \sigma\rVert_{\ell+3/2, \ell+2} &= \lVert \fordiff^2\sigma\rVert_{\ell+3/2, \ell} \\
&= \Big\lVert \tfrac{1}{2} E^2\sigma \lvert E\fordiff^2\eta\rvert^2 + \tfrac{1}{2} \sigma \lvert \fordiff^2\eta\rvert^2 - \lvert E\fordiff\dot{\eta}\rvert^2 \Big\rVert_{\ell+3/2, \ell} \\
&\lesssim \lVert \sigma \lvert \fordiff^2\eta\rvert^2\rVert_{\ell+3/2, \ell} + \lVert \lvert \fordiff\dot{\eta}\rvert^2\rVert_{\ell+3/2, \ell}
\end{align*}
using the fact that $E$ is a bounded operator in any norm. (The technique we use will make it clear that the norm of $E^2\sigma \lvert E\fordiff^2\eta\rvert^2$ is comparable to that of $\sigma \lvert \fordiff^2\eta\rvert$, so there is no reason to study it separately.) 

To simplify notation a bit, let $f = \lvert \fordiff^2\eta\rvert^2$ and $g=\lvert \fordiff\dot{\eta}\rvert^2$. Then the inequality above is
\begin{equation}\label{initialsigmaestimate}
\lVert \sigma\rVert^2_{\ell+3/2, \ell+2} \lesssim \lVert \sigma f\rVert^2_{\ell+3/2, \ell} + \lVert g\rVert^2_{\ell+3/2, \ell}.
\end{equation}

Our first goal is to bound $\lVert \sigma f \rVert_{\ell+3/2, \ell}$ in terms of the norms of $\sigma$ and the 
norms of $f$. Note that we have bounds on $\sigma$ and $\fordiff \sigma$ in the maximum norm by Lemma \ref{acboundslemma}, while for higher differences of $\sigma$ the bounds are expressed in terms of Euclidean-type norms. So the complication comes from taking this into account.

Using the general product formula \eqref{doubleproduct} for differences, we have 
\begin{align*}
\lVert \sigma f\rVert^2_{\ell+3/2, \ell} &= \frac{1}{n} \sum_{k=1}^n s_k^{(\ell+3/2)} \lvert \fordiff^{\ell}(\sigma f)_k\rvert^2 \\
&= \frac{1}{n} \sum_{k=1}^n s_k^{(\ell+3/2)} \left\lvert \sum_{j=0}^{\ell} {\ell \choose j} \fordiff^{\ell-j} \sigma_{k+j} \fordiff^jf_k\right\rvert^2 \\
&\lesssim \sum_{j=0}^{\ell} \frac{1}{n} \sum_{k=1}^n s_k^{(\ell+3/2)} \lvert \fordiff^{\ell-j}\sigma_{k+j}\rvert^2 \lvert \fordiff^jf_k\rvert^2.
\end{align*}
Now the exceptional cases are when $j=\ell$ or $j=\ell-1$, because there we want to use the maximum norm on $\sigma$ directly. In the other cases, we still
need to use the maximum norm, but we will give it some extra weighting so that we can use \eqref{weightedsobolev}. So we have, using $ \frac{\sigma_k}{s_k} \le a$ and $\lvert \fordiff\sigma\rvert \le a$,
that
\begin{align*}
\lVert \sigma f\rVert^2_{\ell+3/2, \ell} &\lesssim \frac{a^2}{n} \sum_{k=1}^n s_k^{(\ell+3/2)} s_{k+\ell}^2 \lvert \fordiff^{\ell}f_k\rvert^2 + \frac{a^2}{n} \sum_{k=1}^n s_k^{(\ell+3/2)} \lvert \fordiff^{\ell-1}f_k\rvert^2 \\ 
&\qquad\qquad + \sum_{j=0}^{\ell-2} \Big(\max_{1\le k\le n} s_k^{(j+2)} \lvert \fordiff^jf_k\rvert^2 \Big)  \left( \frac{1}{n} \sum_{k=1}^n \frac{s_k^{(\ell+3/2)}}{s_k^{(j+2)}} \lvert \fordiff^{\ell-j}\sigma_{k+j}\rvert^2\right).
\end{align*}
Using the bounds $s_k^{(\ell+3/2)} s_{k+\ell}^2 \lesssim s_k^{(\ell+7/2)}$ and $\frac{s_k^{(\ell+3/2)}}{s_k^{(j+2)}} \lesssim s_{k+j}^{(\ell-j-1/2)}$ from Proposition \ref{weightsprop}, we obtain
\begin{equation}\label{2stepsigmaestimate}
\lVert \sigma f\rVert^2_{\ell+3/2, \ell} \lesssim a^2 \lVert f\rVert^2_{\ell+7/2, \ell} + a^2 \lVert f\rVert^2_{\ell+3/2, \ell-1} + \sum_{j=0}^{\ell-2} 
\supnorm{f}^2_{j+2, j} \lVert \sigma\rVert^2_{\ell-j-1/2, \ell-j}.
\end{equation}

A straightforward application of \eqref{doubleproduct} and the basic estimates of Theorem \ref{basicinequalitiesthm} proves the inequalities
\begin{align}
\lVert f\rVert^2_{\ell+7/2, \ell} &\lesssim e_3 e_{\ell+1} + e_{\ell}^2 \text{ for $\ell \ge 0$;} \label{sigmaA}\\
\lVert f\rVert^2_{\ell+3/2, \ell-1} &\lesssim e_3 e_{\ell+1} + e_{\ell}^2 \text{ for $\ell \ge 1$;} \label{sigmaB} \\
\supnorm{f}^2_{j+2, j} &\lesssim e_3 e_{j+3} + e_{j+2}^2 \text{ for $j\ge 0$;} \label{sigmaC}.
\end{align}

To conclude, we need to estimate the other term in \eqref{initialsigmaestimate}. We will show
\begin{equation}\label{sigmaD}
\lVert g\rVert^2_{\ell+3/2, \ell} \lesssim e_3 e_{\ell+1} + e_{\ell}^2 \text{ for $\ell\ge 0$.}
\end{equation}
We will see that the extra half-power in the weighting is only necessary for the norm of $g$; otherwise we could have worked with $\lVert f\rVert^2_{\ell+3, \ell}$, etc. instead.
We have
\begin{align*}
\lVert g\rVert^2_{\ell+3/2, \ell} &\lesssim \textstyle \sum_{j=0}^{\ell} \lVert \langle \fordiff^{j+1}\dot{\eta}, E^j\fordiff^{\ell+1-j}\dot{\eta}\rangle\rVert^2_{\ell+3/2, 0} \\
&\lesssim \supnorm{\dot{\eta}}^2_{1/2, 1} \lVert \dot{\eta}\rVert^2_{\ell+1, \ell+1} + \textstyle \sum_{j=1}^{\ell-1} \supnorm{\dot{\eta}}^2_{j+1, j+1} \lVert \dot{\eta}\rVert^2_{\ell-j, \ell+1-j} \\
&\lesssim e_3 e_{\ell+1} + \textstyle \sum_{j=1}^{\ell-1} e_{j+2} e_{\ell+2-j} \lesssim e_3 e_{\ell+1} + e_{\ell}^2.
\end{align*}

Plugging \eqref{sigmaA}--\eqref{sigmaD} into \eqref{2stepsigmaestimate} and \eqref{initialsigmaestimate}, we obtain 
\begin{equation}\label{finalsigmaestimate}
\lVert \sigma\rVert^2_{\ell+3/2, \ell+2} \lesssim (1+a^2) \big(e_3 e_{\ell+1} + e_{\ell}^2\big) + \textstyle\sum_{j=0}^{\ell-2} (e_3 e_{j+3} + e_{j+2}^2) \lVert \sigma\rVert^2_{\ell-j-1/2, \ell-j}.
\end{equation}

Recall from Lemma \ref{acboundslemma} that $a\lesssim e_2$, so that $1+a^2 \lesssim e_2^2$. 
So when $\ell=0$ or $\ell=1$ we get 
$\lVert \sigma\rVert^2_{3/2, 2} \lesssim e_2^2 e_3 e_1 \le e_3^4$ and 
$\lVert \sigma\rVert^2_{5/2, 3} \lesssim e_2^2 e_3 e_2 \le e_3^4$. From this we have 
$d_1\lesssim e_3^4$ and $d_2\lesssim e_3^4$. 
Now we can derive from \eqref{finalsigmaestimate} for $m\ge 3$ the recursive inequality
\begin{equation*}
d_m 
\lesssim e_2^2 e_3 e_m + e_2^2 e_{m-1}^2 + \textstyle\sum_{j=0}^{m-3} (e_3 e_{j+3} + e_{j+2}^2) d_{m-2-j}.
\end{equation*}
Plugging in $m=3$ gives the base case \eqref{lowdapp}, and induction on $m$ gives \eqref{sigmaboundapp}.
\end{proof}

\subsection{Proof of Theorem \ref{discreteenergyestimate}}\label{A4}

\begin{theo}
Let $n\in \mathbb{N}$, and suppose $(\eta_1(t), \ldots, \eta_n(t))$ and $(\sigma_1(t), \ldots, \sigma_n(t))$ form a solution of \eqref{fancydiscreteevolution} and \eqref{fancydiscretetension}
with $\sigma_0(t)\equiv 0$ and $\eta_{n+1}(t)\equiv 0$, along with the odd extensions \eqref{discreteodd}. 

Then the energies \eqref{discreteenergydef} and \eqref{discretetimeenergydef} satisfy the estimates
\begin{equation}\label{energyderivativeapp}
\frac{d\widetilde{e}_3}{dt} \le M_3 e_3^7
\end{equation}
for some $M_3$ independent of the initial data and of $n$.
In addition the higher energies satisfy 
\begin{equation}\label{generalenergyderivativeapp}
\frac{d\widetilde{e}_m}{dt} \lesssim M_m(e_{m-1}) e_m
\end{equation}
for every $m>3$, where $M_m$ depends only on $e_{m-1}$. 
\end{theo}

\begin{proof}
As with the proof of Lemma \ref{sigmasobolevlemma}, the estimates for the whip and chain are proved in the same way, so we will just focus on the harder case of the chain (where nontrivial technical issues such as Lemma \ref{notnaivelemma} arise). The essential step is the discrete analogue of the computation \eqref{nineandahalf}, together with the integration by parts employed to cancel out the highest-order term. Then we simply estimate the remainder terms using Corollary \ref{supproductboundcorollary} and Lemma \ref{normtoenergy}.

The first step is just to differentiate. We deal with the terms in \eqref{discretetimeenergydef} one at a time. So fix 
an integer $\ell\ge 0$. Then 
\begin{equation}\label{firstenergyderivative}
\frac{d}{dt} \frac{1}{n} \sum_{k=1}^{n-\lfloor \ell/2\rfloor} \Big( \sigma_k^{(\ell)} \lvert \fordiff^{\ell}\dot{\eta}_k\rvert^2
+ \sigma_k^{(\ell+1)} \lvert \fordiff^{\ell+1}\eta_k\rvert^2\Big) = \frac{1}{n} \sum_{k=1}^{n-\lfloor\ell/2\rfloor} (\termi + 2 \termii)
\end{equation}
where 
\begin{equation}\label{sigmaderivative}
\termi = \left(\frac{d}{dt}\sigma_k^{(\ell)}\right) \lvert \fordiff^{\ell}\dot{\eta}_k\rvert^2 + \left(\frac{d}{dt} \sigma_k^{(\ell+1)}\right) \lvert \fordiff^{\ell+1}\eta_k\rvert^2
\end{equation}
and
\begin{equation}\label{messyterm}
\termii = \sigma_k^{(\ell)} \langle \fordiff^{\ell}\dot{\eta}_k, \fordiff^{\ell}\ddot{\eta}_k\rangle + \sigma_k^{(\ell+1)} \langle \fordiff^{\ell+1}\eta_k, \fordiff^{\ell+1}\dot{\eta}_k\rangle.
\end{equation}

For \eqref{sigmaderivative}, if $\ell=1$ then we have $\frac{d\sigma_k}{dt} \le cs_k$, while if $\ell\ge 2$ the derivative $\frac{d\sigma_k^{(\ell)}}{dt}$ is a derivative of a product of $\ell$ terms. Using $\sigma_k \le as_k$ and $\frac{d\sigma_k}{dt} \le cs_k$, we clearly have 
$$ \left\lvert \frac{d}{dt} \sigma_k^{(\ell)}\right\rvert \le \ell a^{\ell-1} c s_k^{(\ell)},$$
which obviously holds in addition for $\ell=0$ and $\ell=1$. 
Thus the sum of \eqref{sigmaderivative} is bounded by 
\begin{equation}\label{termifinal}
\left\lvert \frac{1}{n} \sum_{k=1}^{n-\lfloor \ell/2\rfloor} \termi\right\rvert \le a^{\ell-1}c \lVert \dot{\eta}\rVert^2_{\ell, \ell} + a^{\ell} c \lVert \eta\rVert^2_{\ell+1, \ell+1}  \lesssim e_3^{\ell+2} e_{\ell} \text{ for $\ell\ge 0$}
\end{equation}
using Lemma \ref{acboundslemma}.

So our primary concern is \eqref{messyterm}. By analogy with the technique used to derive \eqref{nineandahalf}, we want to pull out the worst terms of \eqref{messyterm} and collect them into a single exact difference which will sum to zero. It is easy to check that
if $\psi = \sigma^{(\ell+1)} \langle E\fordiff^{\ell}\dot{\eta}, \fordiff^{\ell+1}\eta\rangle$, then 
\begin{multline*}
\backdiff \psi =  (E^{-1}\sigma^{(\ell+1)}) \langle \fordiff^{\ell}\dot{\eta}, \fordiff^{\ell+1} \backdiff\eta\rangle \\+ \sigma^{(\ell+1)} \langle \fordiff^{\ell+1}\eta, \fordiff^{\ell+1}\dot{\eta}\rangle + (\backdiff \sigma^{(\ell+1)}) \langle \fordiff^{\ell}\dot{\eta}, \fordiff^{\ell+1}\eta\rangle.
\end{multline*}
The middle term of $\backdiff\psi$ is precisely the second term of \eqref{messyterm}, so we want to show that the difference of the remaining terms is relatively simple.\footnote{The complication in this computation is the fact that $\ddot{\eta} = \backdiff (\sigma\fordiff\eta) = E^{-1}\fordiff(\sigma\fordiff\eta)$: we have nothing but $\fordiff$ in the rest of the formula, so the appearance of one $\backdiff$ operator (or, equivalently, of one backward shift $E^{-1}$) necessitates rederiving the formulas to get rid of it, rather than using a formula like \eqref{doubleproduct} directly. The reason we don't want to see a mix of $\backdiff$ and $\fordiff$ operators is because later in the proof we will need to use the fact that $\langle \fordiff\eta, \fordiff\dot{\eta}\rangle \equiv 0$, while there is no simple formula for $\langle \backdiff\eta, \fordiff\dot{\eta}\rangle$.} We will do this computation in Lemma \ref{notnaivelemma}.

To check that this all works, we observe that the backward difference $\backdiff \psi$ sums to zero: we have 
$$ \frac{1}{n} \sum_{k=1}^{n-\lfloor\ell/2\rfloor} \backdiff \psi = \psi_{n-\lfloor\ell/2\rfloor} - \psi_0,$$
which can be checked to vanish due to $\sigma_0=0$ and the oddness conditions \eqref{discreteodd}.
(Recall that this is precisely the reason that the summands in our energies \eqref{discreteenergydef}--\eqref{discretetimeenergydef} all terminate at $n-\lfloor\ell/2\rfloor$.) Hence we have 
$$ \sum_{k=1}^{n-\lfloor\ell/2\rfloor} (\termii)_k = \sum_{k=1}^{n-\lfloor\ell/2\rfloor} (\termii-\backdiff\psi)_k.$$

Now by Lemma \ref{notnaivelemma}, we have 
\begin{equation}\label{firstexp}
\termii - \backdiff \psi = \sum_{i=1}^{\ell} R_i,
\end{equation}
where the remainder terms are given by
\begin{equation}\label{remainder}
R_i = \sum_{j=0}^{\ell-i} \textstyle{{j+i\choose i}} \sigma^{(\ell)} (E^{j-1} \fordiff^{i+1}\sigma) \langle \fordiff^{\ell}\dot{\eta}, \fordiff^{\ell+1-i}\eta\rangle.
\end{equation}

Therefore if $\ell\ge 1$ we have 
\begin{equation}\label{termiisum}
\left\lvert \frac{1}{n} \sum_{k=1}^{n-\lfloor\ell/2\rfloor} (\termii)_k \right\rvert = \left\lvert \frac{1}{n} \sum_{k=1}^{n-\lfloor\ell/2\rfloor} \sum_{i=1}^{\ell} (R_i)_k \right\rvert
\lesssim a^{\ell} \sum_{i=1}^{\ell} 
S_{i},
\end{equation}
where 
\begin{equation}\label{Sijdef}
S_{i} = \frac{1}{n} \sum_{k=1}^{n-\lfloor\ell/2\rfloor} s_k^{(\ell)} \lvert \fordiff^{i+1}\sigma_k\rvert \big\lvert \langle \fordiff^{\ell}\dot{\eta}_k, \fordiff^{\ell+1-i}\eta_k\rangle\big\rvert
\end{equation}
for $1\le i\le \ell$. Here we use the fact that by Proposition \ref{weightsprop}, the shift operators $E^{j-1}$ that appear in \eqref{remainder} are bounded for any $j$; since having or not having the shift operators attached to $\sigma$ doesn't change the estimates in any way, we may as well ignore them. (Having the shift operators attached to one of the $\eta$ terms \emph{would} cause a problem, which is why we need Lemma \ref{notnaivelemma}.)

Obviously when $\ell=1$, we must have $i=1$ also in the sum \eqref{termiisum}, and in this case every summand in \eqref{Sijdef} involves $ \langle \fordiff\dot{\eta}_k, \fordiff\eta_k\rangle \equiv 0$, so that $S_{1}=0$ if $\ell=1$. Hence we will assume $\ell\ge 2$ to estimate $S_i$.

We will show that for any $1\le i\le \ell$, 
\begin{equation}\label{mainSbound}
S_i \lesssim \begin{cases} 
e_3^3 e_{\ell} & \text{if $\ell=2$ or $\ell=3$,} \\
L_{i,\ell}(e_{\ell-1}) e_{\ell} & \text{if $\ell \ge 4$,}
\end{cases}
\end{equation}
for some function $L_{i,\ell}$ of $e_{\ell-1}$. 
The bounds for $1\le i\le \ell-1$ are all basically the same,
while the bound for the $i=\ell$ term requires another trick.

Using Corollary \ref{supproductboundcorollary}, and Lemma \ref{normtoenergy}, it is straightforward to verify that $S_1 \lesssim \sqrt{d_3}e_{\ell}$ and that $S_i \lesssim \sqrt{e_{\ell+2-i}} \sqrt{d_{i+1}} \sqrt{e_{\ell}}$
for $2\le i\le \ell-1$. Then using Lemma \ref{sigmasobolevlemma}, we obtain \eqref{mainSbound} when $1\le i\le \ell-1$.

The last case in \eqref{mainSbound} is when $i=\ell$. (Recall that we are assuming $\ell\ge 2$ since $S_1=0$ when $\ell=1$.) Since $\lvert \fordiff\eta\rvert^2 \equiv 1$, we have $\langle \fordiff\eta, \fordiff\dot{\eta}\rangle \equiv 0$. Applying the difference operator $\fordiff^{\ell-1}$ to both sides and using the product formula \eqref{doubleproduct}, we get 
\begin{equation}\label{notnotnaive}
\langle \fordiff^{\ell}\dot{\eta}, \fordiff\eta\rangle = -\sum_{p=1}^{\ell-1} {\ell-1\choose p} \langle E^p\fordiff^{\ell-p}\dot{\eta}, \fordiff^{p+1}\eta\rangle.
\end{equation}
Thus 
$$
S_{\ell} 
\lesssim \sum_{p=1}^{\ell-1} \frac{1}{n} \sum_{k=1}^{n-\lfloor\ell/2\rfloor} s_k^{(\ell)} \lvert \fordiff^{\ell+1}\sigma_{k}\rvert \big\lvert \langle E^p\fordiff^{\ell-p}\dot{\eta}_k, \fordiff^{p+1}\eta_k\rangle \big\rvert.
$$

We obviously want to use the Cauchy-Schwarz inequality on this, and to do so we have to assign the weight $s_k^{(\ell)}$ to the individual pieces. Recall that we have an estimate for $\lVert \sigma\rVert_{\ell+1/2, \ell+1} \le \sqrt{d_{\ell}} \le \sqrt{P_{\ell}(e_{\ell-1}) e_{\ell}}$ from Lemma \ref{sigmasobolevlemma}, which means we must pull the power $s^{(\ell+1/2)}$ out with this term, leaving us with $s^{(\ell-1/2)}$ for what remains. The worst term is the one with $p=\ell-1$, for then we have to estimate $\big\lVert \lvert \fordiff\dot{\eta}\rvert \cdot \lvert \fordiff^{\ell}\eta\rvert \big\rVert_{\ell-1/2, 0}$. We want to pull out the supremum norm of $\fordiff\dot{\eta}$, but we need some positive weight on it in order to be able to use \eqref{weightedsobolev}.\footnote{This is the only place in the paper where we actually \emph{need} to split the weight into noninteger powers to make the estimates work. Without doing this, we cannot close the estimates at the level of $e_3$. This is why the Sobolev norms of $\sigma$ from Lemma \ref{sigmasobolevlemma} are defined the way they are.}

Using Corollary \ref{supproductboundcorollary} again, we compute
\begin{align*}
S_{\ell} 
&\lesssim \textstyle \sum_{p=1}^{\ell-1} \lVert E^{-1}\fordiff^{\ell+1}\sigma\rVert_{\ell+1/2, 0} \Big\lVert \lvert E^p\fordiff^{\ell-p}\dot{\eta}\rvert \lvert \fordiff^{p+1}\eta\rvert \Big\rVert_{\ell-1/2, 0} \\
&\lesssim \textstyle \lVert \sigma\rVert_{\ell+1/2, \ell+1} \bigg( \supnorm{E^{\ell-1}\fordiff\dot{\eta}}_{1/2, 0} \lVert \fordiff^{\ell}\eta\rVert_{\ell-1, 0} \\
&\qquad\qquad \textstyle + \sum_{p=1}^{\ell-2}
\supnorm{\fordiff^{p+1}\eta}_{p,0} \lVert E^p\fordiff^{\ell-p}\dot{\eta}\rVert_{\ell-p-1/2, 0}\bigg) \\
&\lesssim \textstyle\sqrt{d_{\ell}} \left( \supnorm{\dot{\eta}}_{1/2, 1} \lVert \eta\rVert_{\ell-1, \ell} + \sum_{p=1}^{\ell-2} \supnorm{\eta}_{p, p+1} \lVert \dot{\eta}\rVert_{\ell-p-1, \ell-p}\right).
\end{align*}

We easily estimate the quantities here, using Lemma \ref{normtoenergy}, to get 
$$ S_{\ell} \lesssim \sqrt{d_{\ell}} \left( \sqrt{e_3} \sqrt{e_{\ell}} + \sum_{p=1}^{\ell-2} \sqrt{e_{p+2}} \sqrt{e_{\ell-p+1}}\right).$$
Putting this together with the bound for $d_{\ell}$ from Lemma \ref{sigmasobolevlemma}, we get 
\eqref{mainSbound} for the cases $i=\ell$.

Now plugging \eqref{mainSbound} into \eqref{termiisum}, we get 
\begin{equation}\label{termiifinal}
\left\lvert \frac{1}{n} \sum_{k=1}^{n-\lfloor\ell/2\rfloor} (\termii)_k \right\rvert \lesssim
a^{\ell} \sum_{i=1}^{\ell} 
S_{i}
\lesssim \begin{cases}
0 & \ell=1, \\
e_3^6 & \ell=2, \\
e_3^7 & \ell=3, \\
K_{\ell}(e_{\ell-1}) e_{\ell} & \ell\ge 4,
\end{cases}
\end{equation}
for some function $K_{\ell}$ of $e_{\ell-1}$. 

Using \eqref{termifinal} and \eqref{termiifinal} in \eqref{firstenergyderivative}, we get 
\begin{align*}
\frac{d}{dt} \frac{1}{n} \sum_{k=1}^{n-\lfloor \ell/2\rfloor} \Big( \sigma_k^{(\ell)} \lvert \fordiff^{\ell}\dot{\eta}_k\rvert^2
+ \sigma_k^{(\ell+1)} \lvert \fordiff^{\ell+1}\eta_k\rvert^2\Big) 
&\lesssim \begin{cases} 
e_3^2 & \ell=0, \\
e_3^4 & \ell=1, \\
e_3^6 & \ell=2, \\
e_3^7 & \ell=3, \\
M_{\ell}(e_{\ell-1}) e_{\ell} &\ell\ge 4.
\end{cases}
\end{align*}
Now summing from $\ell=0$ to $\ell=m$, we get \eqref{energyderivativeapp} and \eqref{generalenergyderivativeapp}.
\end{proof}

To complete the proof, let us establish the formula \eqref{firstexp}.

\begin{lemma}\label{notnaivelemma}
We have the formula 
\begin{multline}\label{notnaive}
\langle E^{-1}\fordiff^{\ell+1}(\sigma\fordiff\eta), \fordiff^{\ell}\dot{\eta}\rangle 
= (E^{-1}\sigma)\langle \fordiff^{\ell+1} \backdiff \eta, \fordiff^{\ell}\dot{\eta}\rangle 
\\ + \big(\backdiff \sigma^{(\ell+1)}\big) \langle \fordiff^{\ell}\dot{\eta}, \fordiff^{\ell+1}\eta\rangle
+ \sum_{i=1}^{\ell} R_i, 
\end{multline}
where $R_i$ is given by \eqref{remainder}.
\end{lemma}

\begin{proof}
The product formula \eqref{doubleproduct} yields
$$
E^{-1}\fordiff^{\ell+1}(\sigma\fordiff\eta) - (E^{-1}\sigma)\fordiff^{\ell+1} \backdiff \eta 
= \sum_{p=0}^{\ell} {\ell+1\choose p+1} (\fordiff^{p+1}E^{-1}\sigma) (E^p\fordiff^{\ell+1-p}\eta),
$$
Then using $E=1+\frac{\fordiff}{n}$ and some binomial expansions and identities, we obtain 
\begin{equation}\label{notnaiveprime}
E^{-1}\fordiff^{\ell+1}(\sigma\fordiff\eta) - (E^{-1}\sigma)\fordiff^{\ell+1} \backdiff \eta
= \sum_{i=0}^{\ell} \sum_{j=0}^{\ell-i} \textstyle{{j+i\choose i}} (E^{j-1} \fordiff^{i+1}\sigma)(\fordiff^{\ell+1-i}\eta).
\end{equation}
We then compute the inner product of all terms with $\fordiff^{\ell}\dot{\eta}$. 

In this last sum of \eqref{notnaiveprime}, notice that when $i=0$ we get 
$$ R_0 = \sigma^{(\ell)} \bigg( \sum_{j=0}^{\ell} E^{j-1} \fordiff\sigma\bigg) \langle \fordiff^{\ell}\dot{\eta}, \fordiff^{\ell+1}\eta\rangle = \big(\backdiff \sigma^{(\ell+1)}\big) \langle \fordiff^{\ell}\dot{\eta}, \fordiff^{\ell+1}\eta\rangle,$$
using the obvious telescoping, which finally yields \eqref{notnaive}.
\end{proof}

\subsection{Proof of Theorem \ref{uniqueness}}\label{A5}

\begin{theo}
Suppose $\gamma$ and $w$ are functions on $[0,1]$ as in Lemma \ref{uniformestimateslemma}.
If $(\eta_1, \sigma_1)$ and $(\eta_2, \sigma_2)$ are two solutions of \eqref{fullsystem} in $L^{\infty}([0,T], N_4[0,1])\cap W^{1,\infty}([0,T], N_3[0,1])$,
with the same initial conditions $$\eta_1(0,s)=\eta_2(0,s)=\gamma(s) \quad \text{and} \quad \partial_t\eta_1(0,s)=\partial_t\eta_2(0,s) = w(s),$$ then $\eta_1(t,s)=\eta_2(t,s)$ and $\sigma_1(t,s) = \sigma_2(t,s)$ for all $t\in [0,T]$ and all $s\in [0,1]$. 
\end{theo}

\begin{proof}
Define the differences by $\varepsilon = \frac{1}{2}(\eta_2-\eta_1)$ and $\delta = \frac{1}{2} (\sigma_2-\sigma_1)$, and the averages by $\overline{\eta}=\frac{1}{2} (\eta_1+\eta_2)$ and $\overline{\sigma}=\frac{1}{2}(\sigma_1+\sigma_2)$. It is easy to compute that these quantities satisfy the equations
\begin{align}
\deltass &= \big( \lvert \etameandiffss \rvert^2 + \lvert \epsilondiffss\rvert^2\big) \delta + 2\langle \epsilondiffss, \etameandiffss \rangle \overline\sigma - 2\langle \epsilondiffst, \etameandiffst \rangle \label{delta} \\
\epsilondifftt  &= \partial_s\big( \overline{\sigma} \epsilondiffs\big) + \partial_s\big( \delta \etameandiffs \big).\label{varepsilon}
\end{align}

The quantities $\overline{\sigma}$ and $\delta$ have the same boundary conditions as those for $\sigma_1$ and $\sigma_2$; similarly $\overline{\eta}$ and $\varepsilon$ must both be odd through $s=1$ since $\eta_1$ and $\eta_2$ are. Furthermore, the fact that $\lvert \partial_s\eta_1\rvert^2 \equiv 1 \equiv \lvert \partial_s\eta_2\rvert^2$ implies that $\langle \epsilondiffs, \etameandiffs \rangle \equiv 0$.

We now estimate the norms of these quantities; the primary goal is to estimate the norm of $\varepsilon$, but we will need the norms of the other terms to do this. For this purpose, we generalize the quantities $A$, $B$, $C$ in 
\ref{ABCdef}, the quantities $D_m$ from \eqref{sigmasobolevnorm}, 
the quantities $E_m$ from \eqref{correctweightenergy}, and the quantities $\tilde{E}_m$ from \eqref{correctenergy}: we will denote 
$A[\overline{\sigma}] = \sup_{0\le s\le 1} \frac{\lvert \overline{\sigma}(s)\rvert}{s}$, $E_m[\varepsilon] = \sum_{\ell=0}^m \lVert \epsilondifft \rVert^2_{\ell,\ell} + \lVert \varepsilon\rVert^2_{\ell+1, \ell+1}$, etc. For the time-dependent energy $\tilde{E}_m$ we use $\overline{\sigma}$ for the weighting.

It is easy to verify that
\begin{alignat*}{2}
A[\overline{\sigma}] &\le \tfrac{1}{2} A[\sigma_1] + \tfrac{1}{2} A[\sigma_2], &\qquad B[\overline{\sigma}] &\le B[\sigma_1]  + B[\sigma_2]\\
C[\overline{\sigma}] &\le \tfrac{1}{2} C[\sigma_1] + \tfrac{1}{2} C[\sigma_2], &
\qquad D_m[\overline{\sigma}] &\le D_m[\sigma_1] + D_m[\sigma_2],\\
E_m[\overline{\eta}] &\le E_m[\eta_1] + E_m[\eta_2].&&
\end{alignat*}

We could proceed by imitating the proof of Theorem \ref{discreteenergyestimate} to get a bound for the energy $E_2[\varepsilon]$; however it's simpler to use some alternative techniques 
to get a bound for $E_1[\varepsilon]$. The reason this works is that we can separate all the estimates into
low-derivative norms of $\varepsilon$ by compensating with high-derivative norms of $\overline{\eta}$.

First we note that since $\delta(t,0)=0$, we have by the Cauchy-Schwarz inequality that $\delta(t,s)^2 \le s \int_0^1 \deltas(t,s)^2 \, ds$. In the notation of Definition \ref{weightednorms}, we conclude that 
\begin{equation}\label{deltamax}
\supnorm{\delta}_{-1, 0} \le \lVert \delta\rVert_{0,1}.
\end{equation}
Since $\delta(t,0)=0$ and $\deltas(t,1)=0$, we have using \eqref{delta} and Corollary \ref{supproductboundcorollary} that
\begin{align*}
\lVert \delta\rVert^2_{0,1} &= \textstyle \int_0^1 \deltas^2 \, ds = -\int_0^1 \delta \deltass \, ds \\
&\le \textstyle 2 \int_0^1 \overline{\sigma}\delta \lvert \epsilondiffss\rvert \lvert \etameandiffss \rvert \, ds + 2\int_0^1 \delta \lvert \epsilondiffst\rvert \lvert \etameandiffst \rvert \, ds \\
&\le 2\supnorm{\overline{\sigma}}_{-2,0} \supnorm{\delta}_{-1,0} \lVert \varepsilon\rVert_{2,2} \lVert \overline{\eta}\rVert_{1,2} + 2\supnorm{\delta}_{-1,0} \lVert \epsilondifft \rVert_{1,1} \lVert \overline{\eta}\rVert_{0,1}.
\end{align*}
Using the fact that $\supnorm{\overline{\sigma}}_{-2,0} = A[\overline{\sigma}]$, along with the inequality \eqref{deltamax}, we conclude
\begin{equation}\label{rootF}
\begin{split}
\lVert \delta\rVert_{0,1} &\le 2A[\overline{\sigma}] \lVert \varepsilon\rVert_{2,2} \lVert \overline{\eta}\rVert_{1,2} + 2\lVert \epsilondifft\rVert_{1,1} \lVert \etameandifft \rVert_{0,1} \\
 &\lesssim (1+A[\overline{\sigma}]) \sqrt{E_1[\varepsilon]} \sqrt{E_2[\overline{\eta}]} \lesssim (E_2[\overline{\eta}])^{3/2} \sqrt{E_1[\varepsilon]},
\end{split}
\end{equation}
using Lemma \ref{normtoenergy} and Lemma \ref{acboundslemma}.

We compute the energies of $\varepsilon$ using the same technique as in Theorem \ref{discreteenergyestimate}: we try to bound $\frac{d\tilde{E}_m}{dt}$ in terms of $E_m$ for $m=0$ and $m=1$. The lowest one is easy: we have by \eqref{varepsilon} that
\begin{align*}
\frac{d\tilde{E}_0[\varepsilon]}{dt} &= \frac{d}{dt} \int_0^1 \big( \lvert \epsilondifft \rvert^2 + \overline{\sigma} \lvert \epsilondiffs\rvert^2\big) \, ds \\
&= C[\overline{\sigma}] E_0[\varepsilon] 
\textstyle + 2 \int_0^1 \deltas \langle \etameandiffs, \epsilondifft \rangle \, ds + 2\int_0^1 \delta \langle \etameandiffss, \epsilondifft \rangle \, ds,
\end{align*}
since the boundary term 
vanishes. 
Then we can estimate the rest:
\begin{equation}\label{zeroenergyestimate}
\begin{split}
\frac{d\tilde{E}_0[\varepsilon]}{dt} &\le C[\overline{\sigma}] E_0[\varepsilon] + 2 \supnorm{\overline{\eta}}_{0,1} \lVert \delta\rVert_{0,1} \lVert \epsilondifft\rVert_{0,0} + 2 \supnorm{\delta}_{-1,0}
\lVert \overline{\eta}\rVert_{1,2} \lVert \epsilondifft\rVert_{0,0} \\
&\le C[\overline{\sigma}] E_0[\varepsilon] + 2 \lVert \delta\rVert_{0,1} \sqrt{E_0[\varepsilon]} \big( \supnorm{\overline{\eta}}_{0,1} + \lVert \overline{\eta}\rVert_{1,2}\big). 
\end{split}
\end{equation}
Now from Lemma \ref{normtoenergy} we have $\lVert \overline\eta\rVert_{1,2} \lesssim \sqrt{E_2[\overline{\eta}]}$, while the term $\supnorm{\overline{\eta}}_{0,1}$ is a bit more difficult (since we have no weighting on the supremum and can't use \eqref{weightedsobolev}). Instead we use the standard Sobolev inequality \eqref{usualsobolev}:
$$
\supnorm{\overline{\eta}}^2_{0,1} 
\lesssim \lVert \overline{\eta}\rVert^2_{0,1} + \lVert \overline{\eta}\rVert^2_{0,2} \lesssim E_3[\overline{\eta}],$$
using Lemma \ref{normtoenergy}. This yields
\begin{equation}\label{E0varepsilon}
\frac{d\tilde{E}_0[\varepsilon]}{dt} \lesssim E_2[\overline{\eta}]^{3/2} E_3[\overline{\eta}]^{1/2} \sqrt{E_0[\varepsilon]} \sqrt{E_1[\varepsilon]}.
\end{equation}

We thus cannot bound $E_0[\varepsilon]$ without also bounding $E_1[\varepsilon]$; this is not surprising since it's hard to even make sense of equation \eqref{delta} without $\etadiffst$ and $\etadiffss$ both being in $L^2$. Our next step is to perform the same estimates for $\tilde{E}_1[\varepsilon]$, at which point the estimates do close up.

We can compute that
\begin{equation}\label{loudbar}
\begin{split}
&\frac{d}{dt} \textstyle \int_0^1 \big( \overline{\sigma} \lvert \epsilondiffst \rvert^2 + \overline{\sigma}^2 \lvert \epsilondiffss\rvert^2\big) \, ds
 \le C[\overline{\sigma}] \lVert \epsilondifft\rVert^2_{1,1} + 2C[\overline{\sigma}]A[\overline{\sigma}] \lVert \varepsilon\rVert^2_{2,2} \\
&\qquad\qquad\qquad + \textstyle 2\int_0^1 \overline{\sigma} \sigmameandiffss \langle\epsilondiffs, \epsilondiffst \rangle \, ds + 2\int_0^1 \overline{\sigma} \deltass \langle \etameandiffs, \epsilondiffst \rangle \, ds \\
&\qquad\qquad\qquad \textstyle
+ 4\int_0^1 \overline{\sigma} \deltas \langle \etameandiffss, \epsilondiffst\rangle \, ds + 2\int_0^1 \overline{\sigma} \delta \langle \etameandiffsss, \epsilondiffst \rangle \, ds,
\end{split}
\end{equation}
where again the boundary term vanishes.

We also want to use a trick to simplify the estimates a bit---the same trick we used in deriving \eqref{notnotnaive}---to reduce the derivatives on $\varepsilon$ to compensate for the high-derivative term $\deltass$. Since 
$\langle \epsilondiffs, \etameandiffs \rangle
= 0$, we have $\langle \epsilondiffst, \etameandiffs \rangle + \langle \epsilondiffs, \etameandiffst \rangle = 0$, and we use this to write $\int_0^1 \overline{\sigma} \deltass \langle \etameandiffs, \epsilondiffst \rangle \, ds = -\int_0^1 \overline{\sigma} \deltass \langle \etameandiffst, \epsilondiffs \rangle\,ds$. Now \eqref{loudbar} becomes
\begin{multline*}
\textstyle \frac{d}{dt} \int_0^1 \big( \overline{\sigma} \lvert \epsilondiffst \rvert^2 + \overline{\sigma}^2 \lvert \epsilondiffss \rvert^2\big) \, ds \le C[\overline{\sigma}](1+A[\overline{\sigma}]) E_1[\varepsilon] 
 + 2\int_0^1 \overline{\sigma} \lvert \sigmameandiffss \rvert \lvert \epsilondiffs\rvert \lvert \epsilondiffst \rvert \, ds  \\ 
\textstyle + 4\int_0^1 \overline{\sigma} \lvert \deltas \rvert \lvert \etameandiffss \rvert \lvert \epsilondiffst \rvert \, ds + 2
\int_0^1 \overline{\sigma} \lvert \delta\rvert \lvert \etameandiffsss \rvert \lvert \epsilondiffst\rvert \, ds
+ 2 \int_0^1 \overline{\sigma} \lvert \deltass \rvert \lvert \etameandiffst \rvert \lvert \epsilondiffs\rvert \, ds.
\end{multline*}
Using Corollary \ref{supproductboundcorollary}, we can easily bound all but the last term:
\begin{multline}\label{E1varepsilonbound1}
\frac{d}{dt} \int_0^1 \big( \overline{\sigma} \lvert \epsilondiffst\rvert^2 + \overline{\sigma}^2 \lvert \epsilondiffss\rvert^2\big) \, ds
\lesssim C[\overline{\sigma}](1+A[\overline{\sigma}]) E_1[\varepsilon]  \\
+ 2A[\overline{\sigma}] \lVert \epsilondifft\rVert_{1,1} \Big( 
\sqrt{D_3[\overline{\sigma}] E_1[\varepsilon]} + \sqrt{E_3[\overline{\eta}]} \lVert \delta\rVert_{0,1} \Big)
+\textstyle  2 A[\overline{\sigma}] \int_0^1 s \lvert \deltass\rvert \lvert \etameandiffst \rvert \lvert \epsilondiffs \rvert \, ds.
\end{multline}
We can bound $D_3[\overline{\sigma}] \lesssim E_3[\overline{\eta}]^3$ using Lemma \ref{sigmasobolevlemma}, while \eqref{rootF} bounds $\lVert \delta\rVert_{0,1}$.

Finally we deal with the last term of \eqref{E1varepsilonbound1} by plugging in \eqref{delta} and using the fact that $\lvert \etameandiffss \rvert^2 + \lvert \epsilondiffss\rvert^2 = \frac{1}{2} (\lvert \partial_s^2\eta_1\rvert^2 + \lvert \partial_s^2\eta_2\rvert^2)$ in order to eliminate the seemingly nonlinear dependence on $\varepsilon$. After a series of computations as above, we obtain
$$
\int_0^1 s\lvert \deltass\rvert \lvert \etameandiffst \rvert \lvert \epsilondiffs\rvert \, ds 
\lesssim E_3[\overline{\eta}]^3 E_1[\varepsilon].
$$

Therefore \eqref{E1varepsilonbound1} becomes
$$
\frac{d}{dt} \int_0^1 \big( \overline{\sigma} \lvert \epsilondiffst \rvert^2 + \overline{\sigma}^2 \lvert \epsilondiffss \rvert^2\big) \, ds
\lesssim E_3[\overline{\eta}]^3 E_1[\varepsilon],$$
and combining this with \eqref{E0varepsilon}, we obtain
$$\frac{d\tilde{E}_1[\varepsilon]}{dt} \lesssim E_3[\overline{\eta}]^3 E_1[\varepsilon].$$
Now using the inequality $s\le B[\overline{\sigma}]\overline{\sigma}(s)$ for all $s$, we bound $E_1[\varepsilon]$ in terms of $\tilde{E}_1[\varepsilon]$:
$$ E_1[\varepsilon] \le \left(1+B[\overline{\sigma}]\right)^2 \tilde{E}_1[\varepsilon].$$ 
Using $B[\overline{\sigma}]\le B[\sigma_1]+B[\sigma_2]$ and the bound \eqref{Bbound} for $B[\sigma_1]$ and $B[\sigma_2]$ in terms of $E_2[\eta_1]$ and $E_2[\eta_2]$ respectively, we ultimately find that 
\begin{equation}\label{finaltildeEestimate}
 \frac{d\tilde{E}_1}{dt} \le N(t) \tilde{E}_1(t),
\end{equation}
where $N(t)$ is a function depending only on the energies $E_3[\eta_1]$ and $E_3[\eta_2]$, which 
are uniformly bounded by assumption.

Using Gronwall's inequality, we conclude that if $\tilde{E}_1(0)=0$, then $\tilde{E}_1(t)=0$ for all time. 
In particular we conclude that $\int_0^1 \overline{\sigma}(t,s)\lvert \epsilondiffs(t,s)\rvert^2 \, ds = 0$ for all $t\in[0,T]$, so that $\partial_s\varepsilon(t,s) \equiv 0$ for all $t\in[0,T]$ and $s\in[0,1]$. Since $\varepsilon(t,1)=0$, we must have $\varepsilon(t,s)=0$ for all $t$ and $s$, whence we conclude $\eta_1(t,s) = \eta_2(t,s)$ for all $t$ and $s$. The fact that $\sigma_1=\sigma_2$ follows. 
\end{proof}

\makeatletter \renewcommand{\@biblabel}[1]{\hfill#1.}\makeatother

\end{document}